\documentclass[10pt]{amsart}
\usepackage[cp1251]{inputenc}
\usepackage[english]{babel}
\usepackage{amsmath}
\usepackage{amssymb}
\usepackage{amsfonts}

\usepackage[linktocpage=true, colorlinks=true, linkcolor=blue, citecolor=blue, urlcolor=blue]{hyperref}

\begin{document}
\renewcommand{\refname}{References}

\thispagestyle{empty}

\title[Expansion of Iterated Stratonovich Stochastic Integrals]
{Expansion of Iterated Stratonovich Stochastic Integrals of 
Arbitrary Multiplicity Based on Generalized 
Iterated Fourier Series Converging Pointwise}
\author[D.F. Kuznetsov]{Dmitriy F. Kuznetsov}
\address{Dmitriy Feliksovich Kuznetsov
\newline\hphantom{iii} Peter the Great Saint-Petersburg Polytechnic University,
\newline\hphantom{iii} Polytechnicheskaya ul., 29,
\newline\hphantom{iii} 195251, Saint-Petersburg, Russia}%
\email{sde\_kuznetsov@inbox.ru}

\maketitle {\small
\begin{quote}
\noindent{\sc Abstract.} 
The article is devoted to the 
expansion of iterated Stratonovich stochastic integrals of arbitrary 
multiplicity $k$ $(k\in\mathbb{N})$ based on 
iterated trigonometric Fourier series
converging pointwise. The case of iterated Fourier--Legendre series
is considered in details for $k=2$. 
The obtained expansions provide a possibility to 
represent the iterated Stratonovich sto\-chas\-tic integral 
in the form of iterated series 
of products of standard Gaussian random variables. Convergence
in the mean of degree $2n$ $(n\in \mathbb{N})$ of the 
expansions is proved.
Some recent results on the expansion of iterated Stratonovich 
stochastic integrals of multiplicities 3 to 6 are given.
The results of the article can be applied to the numerical solution
of Ito stochastic differential equations.
\medskip
\end{quote}
}


\setlength{\baselineskip}{2.2em}

\tableofcontents

\setlength{\baselineskip}{1.2em}


\section{Introduction}

\vspace{5mm}

The idea of representing of iterated Ito and Stratonovich stochastic 
integrals
in the form of multiple stochastic integrals from 
specific discontinuous nonrandom functions of several variables and following 
expansion of these functions using generalized 
iterated and multiple Fourier series in order 
to get effective mean-square approximations of the mentioned stochastic 
integrals was proposed and developed in a lot of 
publications of the author 
\cite{1997}-\cite{xxxxx}.
The terms "generalized iterated Fourier series"\ and  "generalized multiple
Fourier series"\ mean
that these series are constructed using various complete orthonormal systems 
of functions in the space $L_2([t, T])$, and not only using the trigonometric 
system of functions. Here $[t, T]$ is an interval of integration 
of iterated Ito and Stratonovich stochastic integrals.
For the first time approach of generalized iterated and multiple 
Fourier series is considered 
in \cite{1997} (1997), \cite{1a} (1998), and \cite{3a} (2006)
(also see references to early publications (1994-1996) in 
\cite{1997}, \cite{1a}, \cite{3a}, \cite{10a}-\cite{12aa-afterxxx}). 
Usage of the Fourier--Legendre series 
for ap\-pro\-xi\-ma\-ti\-on of iterated Ito and Stratonovich stochastic 
integrals took place for 
the first time in \cite{1997} (1997) (also see
\cite{1a}-\cite{xxxxx}).
The results from \cite{1997}-\cite{xxxxx} and this work
convincingly testify
that there is a doubtless relation between 
the multiplier factor $1/2$, 
which is typical for Stratonovich stochastic integral and included 
into the sum connecting Stratonovich and Ito stochastic integrals, 
and the fact that in the point of finite discontinuity of piecewise
smooth function $f(x)$ its generalized Fourier series  
converges to the value $(f(x+0)+f(x-0))/2.$ 
In addition, as it is demonstrated in \cite{1997}-\cite{xxxxx}, 
the final formulas for 
expansions of iterated Stratonovich stochastic integrals
based on the Fourier--Legendre series are essentially simpler 
than its analogues
based on the trigonometric Fourier series.
Note that another approaches to approximation of iterated
Ito and Stratonovich stochastic integrals 
can be found in \cite{01}-\cite{02}.
For example, in \cite{3a}-\cite{art-zero} the method
of expansion of iterated Ito stochastic integrals
based on generalized multiple Fourier series is proposed
and developed. The ideas underlying this method 
are close to the ideas of the method considered in this article.

\vspace{5mm}

\section{Theorem on Expansion of 
Iterated Stratonovich Stochastic Integrals of 
Arbitrary Multiplicity}

\vspace{5mm}

Let $(\Omega,$ ${\rm F},$ ${\sf P})$ be a complete probability space, let 
$\{{\rm F}_t, t\in[0,T]\}$ be a nondecreasing right-continous family of 
$\sigma$-algebras of ${\rm F},$
and let ${\bf f}_t$ be a standard $m$-dimensional 
Wiener stochastic process, which is
${\rm F}_t$-measurable for any $t\in[0, T].$ We assume that the components
${\bf f}_{t}^{(i)}$ $(i=1,\ldots,m)$ of this process are independent. 

Consider the following iterated Stratonovich and Ito stochastic integrals

\begin{equation}
\label{str}
J^{*}[\psi^{(k)}]_{T,t}=
{\int\limits_t^{*}}^T\psi_k(t_k) \ldots {\int\limits_t^{*}}^{t_2}
\psi_1(t_1) d{\bf w}_{t_1}^{(i_1)}\ldots
d{\bf w}_{t_k}^{(i_k)},
\end{equation}

\begin{equation}
\label{ito}
J[\psi^{(k)}]_{T,t}=\int\limits_t^T\psi_k(t_k) \ldots \int\limits_t^{t_{2}}
\psi_1(t_1) d{\bf w}_{t_1}^{(i_1)}\ldots
d{\bf w}_{t_k}^{(i_k)},
\end{equation}

\vspace{3mm}
\noindent
where every $\psi_l(\tau)$ $(l=1,\ldots,k)$ is a 
nonrandom function 
on $[t,T],$ ${\bf w}_{\tau}^{(i)}={\bf f}_{\tau}^{(i)}$
for $i=1,\ldots,m$ and
${\bf w}_{\tau}^{(0)}=\tau,$
$i_1,\ldots,i_k = 0, 1,\ldots,m,$

\vspace{-1mm}
$$
\int\limits^{*}
\hbox{and}\ \int\limits
$$ 

\vspace{2mm}
\noindent
denote Stratonovich and Ito stochastic integrals,
respectively (in this paper, 
we use the definition of the Stratonovich stochastic integral from \cite{1}).

Further, we will 
denote
the complete orthonormal systems of Legendre polynomials and 
tri\-go\-no\-met\-ric functions 
in the space $L_2([t, T])$ as $\{\phi_j(x)\}_{j=0}^{\infty}$.
We will also pay attention on the 
following well-known facts about these two systems of functions.

\vspace{2mm}

{\it Suppose that the function $f(x)$ is 
bounded at the interval $[t, T].$ Moreover, its derivative
$f'(x)$ is continuous function at the interval $[t, T]$ except may be
the finite number of points 
of the finite discontinuity.
Then the generalized Fourier series 

\vspace{-1mm}
$$
\sum\limits_{j=0}^{\infty}
C_j\phi_j(x)
$$

\vspace{2mm}
\noindent
with the Fourier coefficients

\vspace{-1mm}
$$
C_j=\int\limits_t^{T}f(x)\phi_j(x)dx
$$

\vspace{2mm}
\noindent
converges at any internal point $x$ of 
the interval $[t, T]$ to the value 
$\left(f(x+0)+f(x-0)\right)/2$ and converges  
uniformly to $f(x)$ on any closed interval of continuity 
of the function
$f(x)$ laying inside 
$[t, T]$. At the same time the Fourier--Legendre series 
converges 
if $x=t$ and $x=T$ to $f(t+0)$ and $f(T-0)$ 
correspondently, and the trigonometric Fourier series converges if   
$x=t$ and $x=T$ to $\left(f(t+0)+f(T-0)\right)/2$ 
in the case of periodic continuation 
of the function $f(x)$}.

\vspace{2mm}

Define the following function on the hypercube $[t, T]^k$

\vspace{-1mm}
\begin{equation}
\label{ppp}
K(t_1,\ldots,t_k)=
\begin{cases}
\psi_1(t_1)\ldots \psi_k(t_k),\ t_1<\ldots<t_k\\
~\\
~\\
0,\ \hbox{\rm otherwise}
\end{cases}
=\ \
\prod\limits_{l=1}^k
\psi_l(t_l)\ \prod\limits_{l=1}^{k-1}{\bf 1}_{\{t_l<t_{l+1}\}},\ 
\end{equation}

\vspace{3mm}
\noindent
where $t_1,\ldots,t_k\in [t, T]$ $(k\ge 2)$ and 
$K(t_1)\equiv\psi_1(t_1)$ for $t_1\in[t, T].$ Here 
${\bf 1}_A$ denotes the indicator of the set $A$.

Let us formulate the following statement. 

\vspace{2mm}

{\bf Theorem 1} \cite{10a} (Sect.~2.4) (also see \cite{1997} (1997), 
\cite{1a}, \cite{5b}-\cite{2017a}, \cite{2017}, \cite{10},
\cite{10aaaa}-\cite{12aa-afterxxx}, \cite{xxxxx}).
{\it Suppose that every function $\psi_l(\tau)$ $(l=1,\ldots,k)$ is twice continuously
differentiable at the interval
$[t, T]$ and
$\{\phi_j(x)\}_{j=0}^{\infty}$ is a complete
orthonormal system of trigonometric functions in the space $L_2([t, T])$. 
Then, for the iterated Stratonovich stochastic integral 
$J^{*}[\psi^{(k)}]_{T,t}$ defined by {\rm(\ref{str})}
the following
expansion

\vspace{-1mm}
\begin{equation}
\label{1500}
J^{*}[\psi^{(k)}]_{T,t}=
\sum_{j_1=0}^{\infty}\ldots\sum_{j_k=0}^{\infty}
C_{j_k\ldots j_1}
\prod_{l=1}^k
\zeta^{(i_l)}_{j_l}
\end{equation}

\vspace{4mm}
\noindent
conver\-ging 
in the mean of degree $2n$ $(n\in \mathbb{N})$
is valid, where 

\vspace{-1mm}
$$
\zeta_{j}^{(i)}=
\int\limits_t^T \phi_{j}(s) d{\bf w}_s^{(i)}
$$ 

\vspace{3mm}
\noindent
are independent standard Gaussian random variables
for
various
$i$ or $j$ {\rm(}if $i\ne 0${\rm)} and 

\vspace{-1mm}
\begin{equation}
\label{333.40}
C_{j_k\ldots j_1}=\int\limits_{[t,T]^k}
K(t_1,\ldots,t_k)\prod_{l=1}^{k}\phi_{j_l}(t_l)dt_1\ldots dt_k
\end{equation}

\vspace{3mm}
\noindent
is the Fourier coefficient.}

\vspace{2mm}

Note that (\ref{1500})  means the following 

\vspace{-2mm}
\begin{equation}
\label{1500e}
\lim\limits_{p_1\to\infty}
\varlimsup\limits_{p_2\to\infty}
\ldots\varlimsup\limits_{p_k\to\infty}
{\sf M}\left\{\left(J^{*}[\psi^{(k)}]_{T,t}-
\sum_{j_1=0}^{p_1}\ldots\sum_{j_k=0}^{p_k}
C_{j_k\ldots j_1}
\prod_{l=1}^k
\zeta^{(i_l)}_{j_l}\right)^{2n}\right\}=0,
\end{equation}

\vspace{3mm}
\noindent
where $\varlimsup$ means $\limsup$.

\vspace{2mm}

{\bf Proof.} Let us condider several lemmas.
Define the function $K^{*}(t_1,\ldots,t_k)$ on the hypercube 
$[t,T]^k$ as follows

\vspace{2mm}
$$
K^{*}(t_1,\ldots,t_k)=\prod\limits_{l=1}^k\psi_l(t_l)
\prod_{l=1}^{k-1}\biggl({\bf 1}_{\{t_l<t_{l+1}\}}+
\frac{1}{2}{\bf 1}_{\{t_l=t_{l+1}\}}\biggr)=
$$

\vspace{3mm}
\begin{equation}
\label{1999.1}
=\prod_{l=1}^k \psi_l(t_l)\left(\prod_{l=1}^{k-1}
{\bf 1}_{\{t_l<t_{l+1}\}}+
\sum_{r=1}^{k-1}\frac{1}{2^r}
\sum_{\stackrel{s_r,\ldots,s_1=1}{{}_{s_r>\ldots>s_1}}}^{k-1}\ 
\prod_{l=1}^r {\bf 1}_{\{t_{s_l}=t_{s_l+1}\}}
\prod_{\stackrel{l=1}{{}_{l\ne s_1,\ldots, s_r}}}^{k-1}
{\bf 1}_{\{t_{l}<t_{l+1}\}}\right)
\end{equation}

\vspace{8mm}
\noindent
for $t_1,\ldots,t_k\in[t, T]$\ $(k\ge 2)$ and 
$K^{*}(t_1)\equiv\psi_1(t_1)$ for $t_1\in[t, T],$ 
where ${\bf 1}_A$ is the indicator of the set $A$.

\vspace{2mm}

{\bf Lemma 1} \cite{1997} (1997), 
\cite{1a}, \cite{5b}-\cite{2017a}, \cite{2017}-\cite{12aa-afterxxx}, \cite{xxxxx}.
{\it Under the conditions of Theorem {\rm 1} the function
$K^{*}(t_1,\ldots,t_k)$
is represented in any internal point of the hypercube   
$[t,T]^k$ by the generalized iterated Fourier series

$$
K^{*}(t_1,\ldots,t_k)=
\lim\limits_{p_1\to\infty}\ldots\lim\limits_{p_k\to\infty}
\sum_{j_1=0}^{p_1}\ldots\sum_{j_k=0}^{p_k}
C_{j_k\ldots j_1}\prod_{l=1}^{k} \phi_{j_l}(t_l)
\stackrel{\sf def}{=}
$$

\vspace{2mm}
\begin{equation}
\label{30.18}
\stackrel{\sf def}{=}
\sum_{j_1=0}^{\infty}\ldots \sum_{j_k=0}^{\infty}
C_{j_k\ldots j_1}\prod_{l=1}^{k} \phi_{j_l}(t_l),\ \ \ 
(t_1,\ldots,t_k)\in (t, T)^k,
\end{equation}

\vspace{6mm}
\noindent
where $C_{j_k\ldots j_1}$ has the form {\rm (\ref{333.40})}.
At that, the iterated series {\rm (\ref{30.18})} converges at the 
boundary
of the hypercube $[t,T]^k$
{\rm (}not necessarily to the function $K^{*}(t_1,\ldots,t_k)${\rm )}.}

\vspace{2mm}

{\bf Proof.} We will perform the proof using induction. 
Consider the case $k=2.$ Let us expand the function 
$K^{*}(t_1,t_2)$ using the variable 
$t_1$, when $t_2$ is fixed, into the generalized Fourier series 
at the interval $(t, T)$

\vspace{-1mm}
\begin{equation}
\label{leto8001}
K^{*}(t_1,t_2)=
\sum_{j_1=0}^{\infty}C_{j_1}(t_2)\phi_{j_1}(t_1)\ \ \ (t_1\ne t, T),
\end{equation}

\vspace{2mm}
\noindent
where

\vspace{-2mm}
$$
C_{j_1}(t_2)=\int\limits_t^T
K^{*}(t_1,t_2)\phi_{j_1}(t_1)dt_1=\int\limits_t^T
K(t_1,t_2)\phi_{j_1}(t_1)dt_1=
$$

$$
=\psi_2(t_2)
\int\limits_t^{t_2}\psi_1(t_1)\phi_{j_1}(t_1)dt_1.
$$

\vspace{3mm}

The equality (\ref{leto8001}) is 
fulfilled
pointwise at each point of the interval $(t, T)$ with respect to
the variable $t_1$, when $t_2\in [t, T]$ is fixed, due to 
the piecewise
smoothness of the function $K^{*}(t_1,t_2)$ with respect to the variable 
$t_1\in [t, T]$ ($t_2$ is fixed). 

Note also that due to the well-known properties of the Fourier series, 
the series (\ref{leto8001}) converges when $t_1=t$ and $t_1=T$ 
{\rm (}not necessarily to the function $K^{*}(t_1,t_2)${\rm )}.

Obtaining (\ref{leto8001}) we also used the fact that the right-hand side 
of (\ref{leto8001}) converges when $t_1=t_2$ (point of a finite
discontinuity
of the function $K(t_1,t_2)$) to the value

\vspace{1mm}
$$
\frac{1}{2}\left(K(t_2-0,t_2)+K(t_2+0,t_2)\right)=
\frac{1}{2}\psi_1(t_2)\psi_2(t_2)=
K^{*}(t_2,t_2).
$$

\vspace{4mm}

The function $C_{j_1}(t_2)$ is a continuously differentiable
one
at the interval $[t, T]$. 
Let us expand it into the generalized Fourier 
series at the interval $(t, T)$

\begin{equation}
\label{leto8002}
C_{j_1}(t_2)=
\sum_{j_2=0}^{\infty}C_{j_2 j_1}\phi_{j_2}(t_2)\ \ \ (t_2\ne t, T),
\end{equation}

\vspace{2mm}
\noindent
where 

\vspace{-2mm}
$$
C_{j_2 j_1}=\int\limits_t^T
C_{j_1}(t_2)\phi_{j_2}(t_2)dt_2=
\int\limits_t^T
\psi_2(t_2)\phi_{j_2}(t_2)\int\limits_t^{t_2}
\psi_1(t_1)\phi_{j_1}(t_1)dt_1 dt_2,
$$

\vspace{3mm}
\noindent
and the equality (\ref{leto8002}) is fulfilled pointwise at any point 
of the interval $(t, T)$. The right-hand side 
of
(\ref{leto8002}) converges 
when $t_2=t$ and $t_2=T$ (not necessarily to $C_{j_1}(t_2)$).

Let us substitute (\ref{leto8002}) into (\ref{leto8001})

\begin{equation}
\label{leto8003}
K^{*}(t_1,t_2)=
\sum_{j_1=0}^{\infty}\sum_{j_2=0}^{\infty}C_{j_2 j_1}
\phi_{j_1}(t_1)\phi_{j_2}(t_2),\ \ \ (t_1, t_2)\in (t, T)^2.
\end{equation}

\vspace{3mm}

Note that 
the series on the right-hand side of (\ref{leto8003}) converges at the 
boundary
of the square  $[t, T]^2$ (not necessarily to $K^{*}(t_1,t_2)$).
Lemma 1 is proved for the case $k=2.$

Note that proving Lemma 1 for the case $k=2$, we get the 
following equality (see (\ref{leto8001}))

\begin{equation}
\label{oop1}
\psi_1(t_1)\left({\bf 1}_{\{t_1<t_2\}}+
\frac{1}{2}{\bf 1}_{\{t_1=t_2\}}\right)=
\sum\limits_{j_1=0}^{\infty}\int\limits_{t}^{t_2}\psi_1(t_1) 
\phi_{j_1}(t_1)dt_1 \cdot \phi_{j_1}(t_1),
\end{equation}

\vspace{4mm}
\noindent
which is 
fulfilled pointwise at the interval $(t, T),$
besides
the series on the right-hand side 
of (\ref{oop1}) converges when $t_1=t$ and $t_1=T.$

Let us introduce the assumption of induction

\vspace{2mm}
$$
\sum\limits_{j_1=0}^{\infty}\sum\limits_{j_2=0}^{\infty}\ldots
\sum\limits_{j_{k-2}=0}^{\infty}\psi_{k-1}(t_{k-1})
\int\limits_t^{t_{k-1}}\psi_{k-2}(t_{k-2})\phi_{j_{k-2}}(t_{k-2})
\ldots
$$

\vspace{2mm}
$$
\ldots\int\limits_t^{t_{2}}
\psi_{1}(t_{1})\phi_{j_{1}}(t_{1})dt_1\ldots dt_{k-2}
\prod_{l=1}^{k-2}\phi_{j_{l}}(t_{l})=
$$

\vspace{1mm}
\begin{equation}
\label{oop22}
=\prod\limits_{l=1}^{k-1}\psi_l(t_l)
\prod_{l=1}^{k-2}\left({\bf 1}_{\{t_l<t_{l+1}\}}+
\frac{1}{2}{\bf 1}_{\{t_l=t_{l+1}\}}\right).
\end{equation}

\vspace{4mm}

Then

\vspace{1mm}
$$
\sum\limits_{j_1=0}^{\infty}\sum\limits_{j_2=0}^{\infty}\ldots
\sum\limits_{j_{k-1}=0}^{\infty}\psi_{k}(t_{k})
\int\limits_t^{t_{k}}\psi_{k-1}(t_{k-1})\phi_{j_{k-1}}(t_{k-1})\ldots
$$

\vspace{2mm}
$$
\ldots
\int\limits_t^{t_{2}}\psi_{1}(t_{1})
\phi_{j_{1}}(t_{1})
dt_1\ldots dt_{k-1}\prod_{l=1}^{k-1}\phi_{j_{l}}(t_{l})=
$$

\vspace{2mm}
$$
=\sum\limits_{j_1=0}^{\infty}\sum\limits_{j_2=0}^{\infty}\ldots
\sum\limits_{j_{k-2}=0}^{\infty}
\psi_k(t_k)\left({\bf 1}_{\{t_{k-1}<t_{k}\}}+
\frac{1}{2}{\bf 1}_{\{t_{k-1}=t_{k}\}}\right)\psi_{k-1}(t_{k-1})\times
$$

\vspace{3mm}
$$
\times
\int\limits_t^{t_{k-1}}\psi_{k-2}(t_{k-2})\phi_{j_{k-2}}(t_{k-2})\ldots
\int\limits_t^{t_{2}}\psi_{1}(t_{1})
\phi_{j_{1}}(t_{1})dt_1\ldots dt_{k-2}
\prod_{l=1}^{k-2}\phi_{j_{l}}(t_{l})=
$$

\vspace{3mm}
$$
=\psi_k(t_k)\left({\bf 1}_{\{t_{k-1}<t_{k}\}}+
\frac{1}{2}{\bf 1}_{\{t_{k-1}=t_{k}\}}\right)
\sum\limits_{j_1=0}^{\infty}\sum\limits_{j_2=0}^{\infty}\ldots
\sum\limits_{j_{k-2}=0}^{\infty}
\psi_{k-1}(t_{k-1})\times
$$

\vspace{3mm}
$$
\times
\int\limits_t^{t_{k-1}}\psi_{k-2}(t_{k-2})\phi_{j_{k-2}}(t_{k-2})\ldots
\int\limits_t^{t_{2}}\psi_{1}(t_{1})
\phi_{j_{1}}(t_{1})dt_1\ldots dt_{k-2}
\prod_{l=1}^{k-2}\phi_{j_{l}}(t_{l})=
$$

\vspace{3mm}
$$
=\psi_k(t_k)\left({\bf 1}_{\{t_{k-1}<t_{k}\}}+
\frac{1}{2}{\bf 1}_{\{t_{k-1}=t_{k}\}}\right)
\prod\limits_{l=1}^{k-1}\psi_l(t_l)
\prod_{l=1}^{k-2}\left({\bf 1}_{\{t_l<t_{l+1}\}}+
\frac{1}{2}{\bf 1}_{\{t_l=t_{l+1}\}}\right)=
$$

\vspace{3mm}
\begin{equation}
\label{oop30}
=\prod\limits_{l=1}^{k}\psi_l(t_l)
\prod_{l=1}^{k-1}\left({\bf 1}_{\{t_l<t_{l+1}\}}+
\frac{1}{2}{\bf 1}_{\{t_l=t_{l+1}\}}\right).
\end{equation}

\vspace{7mm}

On the other hand, the left-hand side 
of (\ref{oop30}) can be represented 
in the following form

\vspace{1mm}
$$
\sum_{j_1=0}^{\infty}\ldots \sum_{j_k=0}^{\infty}
C_{j_k\ldots j_1}\prod_{l=1}^{k} \phi_{j_l}(t_l)
$$

\vspace{3mm}
\noindent
by 
expanding the function

\vspace{1mm}
$$
\psi_{k}(t_{k})
\int\limits_t^{t_{k}}\psi_{k-1}(t_{k-1})\phi_{j_{k-1}}(t_{k-1})\ldots
\int\limits_t^{t_{2}}\psi_{1}(t_{1})
\phi_{j_{1}}(t_{1})
dt_1\ldots dt_{k-1}
$$

\vspace{3mm}
\noindent
into the generalized Fourier series at the interval $(t, T)$ 
using the variable 
$t_k$.
Lemma 1 is proved.

Let us introduce the following notations

\vspace{2mm}

$$
J[\psi^{(k)}]_{T,t}^{s_l,\ldots,s_1}\ \ \stackrel{\rm def}{=}\ \
\prod_{p=1}^l {\bf 1}_{\{i_{s_p}=
i_{s_{p}+1}\ne 0\}}\ \times
$$

\vspace{1mm}
$$
\times\
\int\limits_t^T\psi_k(t_k)\ldots \int\limits_t^{t_{s_l+3}}
\psi_{s_l+2}(t_{s_l+2})
\int\limits_t^{t_{s_l+2}}\psi_{s_l}(t_{s_l+1})
\psi_{s_l+1}(t_{s_l+1})\ \times
$$

\vspace{1mm}
$$
\times\
\int\limits_t^{t_{s_l+1}}\psi_{s_l-1}(t_{s_l-1})
\ldots
\int\limits_t^{t_{s_1+3}}\psi_{s_1+2}(t_{s_1+2})
\int\limits_t^{t_{s_1+2}}\psi_{s_1}(t_{s_1+1})
\psi_{s_1+1}(t_{s_1+1})\ \times
$$

\vspace{1mm}
$$
\times\
\int\limits_t^{t_{s_1+1}}\psi_{s_1-1}(t_{s_1-1})
\ldots \int\limits_t^{t_2}\psi_1(t_1)
d{\bf w}_{t_1}^{(i_1)}\ldots d{\bf w}_{t_{s_1-1}}^{(i_{s_1-1})}
dt_{s_1+1}d{\bf w}_{t_{s_1+2}}^{(i_{s_1+2})}\ldots
$$

\vspace{1mm}
\begin{equation}
\label{30.1}
\ldots\
d{\bf w}_{t_{s_l-1}}^{(i_{s_l-1})}
dt_{s_l+1}d{\bf w}_{t_{s_l+2}}^{(i_{s_l+2})}\ldots d{\bf w}_{t_k}^{(i_k)},
\end{equation}

\vspace{4mm}
\noindent
where 
\begin{equation}
\label{30.5550001}
{\rm A}_{k,l}=\left\{(s_l,\ldots,s_1):\
s_l>s_{l-1}+1,\ldots,s_2>s_1+1,\ s_l,\ldots,s_1=1,\ldots,k-1\right\},
\end{equation}

\vspace{-2mm}
$$
(s_l,\ldots,s_1)\in{\rm A}_{k,l},\ \
l=1,\ldots,\left[k/2\right],\ \
i_s=0, 1,\ldots,m,\ \
s=1,\ldots,k,
$$

\vspace{5mm}
\noindent
$[x]$ is an
integer
part of a real number $x,$\
${\bf 1}_A$ is the indicator of the set $A$.

Let us formulate the statement on relation
between
iterated 
Ito and Stratonovich stochastic integrals 
$J^{*}[\psi^{(k)}]_{T,t},$ $J[\psi^{(k)}]_{T,t}$ 
of fixed multiplicity $k$ (see (\ref{str}), (\ref{ito})).

\vspace{2mm}

{\bf Lemma 2}\ \cite{10a} (Sect.~2.4) (also see \cite{1997} (1997), 
\cite{1a}, \cite{5b}-\cite{2017a}, \cite{2017}, \cite{10},
\cite{10aaaa}-\cite{12aa-afterxxx}).\
{\it Suppose that
every $\psi_l(\tau)$ $(l=1,\ldots,k)$ is a continuous
nonrandom
function at the interval $[t, T]$.
Then, the following relation between iterated
Ito and Stra\-to\-no\-vich stochastic integrals is correct

\begin{equation}
\label{30.4}
J^{*}[\psi^{(k)}]_{T,t}=J[\psi^{(k)}]_{T,t}+
\sum_{r=1}^{\left[k/2\right]}\frac{1}{2^r}
\sum_{(s_r,\ldots,s_1)\in {\rm A}_{k,r}}
J[\psi^{(k)}]_{T,t}^{s_r,\ldots,s_1}\ \ \ \hbox{{\rm w.\ p.\ 1}},
\end{equation}

\vspace{4mm}
\noindent
where $\sum\limits_{\emptyset}$ is supposed to be equal to zero{\rm ;}
hereinafter w.\ p.\ {\rm 1}  means  "with probability {\rm 1}".}

\vspace{2mm}

{\bf Proof.} Let us prove the equality (\ref{30.4}) using induction. 
The case $k=1$ is obvious.
If $k=2,$ then from (\ref{30.4}) we get 

\vspace{-3mm}
\begin{equation}
\label{30.6}
J^{*}[\psi^{(2)}]_{T,t}=J[\psi^{(2)}]_{T,t}+
\frac{1}{2}J[\psi^{(2)}]_{T,t}^{1}\ \ \ \hbox{w.\ p.\ 1}.
\end{equation}

\vspace{4mm}

Let us demonstrate that the equality (\ref{30.6}) is correct 
w.~p.~1. In order to do it let us consider the 
function $F(x,\tau)=x\psi_2(\tau)$ and the
process $F(\eta_{\tau,t},\tau),$ where
$\eta_{\tau,t}=J[\psi^{(1)}]_{\tau,t},$ 
$\tau\in[t, T]$. Then

\vspace{1mm}
\begin{equation}
\label{30.7}
\frac{\partial F}{\partial x}(x,\tau)=\psi_2(\tau),\ \ \ 
d\eta_{\tau,t}=\psi_1(\tau)d{\bf w}_{\tau}^{(i_1)}.
\end{equation}

\vspace{4mm}

From (\ref{30.7}) we obtain that the diffusion 
coefficient of the process $\eta_{\tau,t},$ $\tau \in [t, T]$ equals 
to 
${\bf 1}_{\{i_1\ne 0\}}\psi_1(\tau).$
Further, using the standard relations between 
Stratonovich and Ito stochastic integrals \cite{1}
(also see \cite{10a} (Sect.~2.4)), 
we obtain the relation (\ref{30.6}). 
Thus, the statement
of Lemma 2 is proved for $k=1$ and $k=2.$

Assume that the statement of Lemma 2 is correct
for some integer $k$ $(k>2),$ and let us prove its correctness when 
the value $k$ is greater per unit. Using the assumption of induction, 
we obtain w.~p.~1

$$
J^{*}[\psi^{(k+1)}]_{T,t}=
$$

\vspace{2mm}

$$
=
{\int\limits_t^{*}}^T \psi_{k+1}(\tau)
\left(J[\psi^{k}]_{\tau,t}+
\sum_{r=1}^{\left[k/2\right]}\frac{1}{2^r}
\sum_{(s_r,\ldots,s_1)\in {\rm A}_{k,r}}
J[\psi^{(k)}]_{\tau,t}^{s_r,\ldots,s_1}
\right)
d{\bf w}_{\tau}^{(i_{k+1})}=
$$

\vspace{2mm}
$$
={\int\limits_t^{*}}^T \psi_{k+1}(\tau)
J[\psi^{(k)}]_{\tau,t}d{\bf w}_{\tau}^{(i_{k+1})}+
$$

\vspace{2mm}

\begin{equation}
\label{30.8}
+
\sum_{r=1}^{\left[k/2\right]}\frac{1}{2^r}
\sum_{(s_r,\ldots,s_1)\in {\rm A}_{k,r}}
{\int\limits_t^{*}}^T\psi_{k+1}(\tau)
J[\psi^{(k)}]_{\tau,t}^{s_r,\ldots,s_1}d{\bf w}_{\tau}^{(i_{k+1})}.
\end{equation}

\vspace{7mm}

Applying the Ito formula and the standard relation
between 
Stratonovich and Ito stochastic integrals, we get w. p. 1

\vspace{1mm}
\begin{equation}
\label{30.9}
{\int\limits_t^{*}}^T \psi_{k+1}(\tau)
J[\psi^{(k)}]_{\tau,t}d{\bf w}_{\tau}^{(i_{k+1})}=J[\psi^{(k+1)}]_{T,t}+
\frac{1}{2}J[\psi^{(k+1)}]_{T,t}^{k},
\end{equation}

\vspace{2mm}

$$
{\int\limits_t^{*}}^T\psi_{k+1}(\tau)
J[\psi^{(k)}]_{\tau,t}^{s_r,\ldots,s_1}d{\bf w}_{\tau}^{(i_{k+1})}=
$$

\vspace{2mm}
\begin{equation}
\label{30.10}
=\left\{
\begin{matrix}
J[\psi^{(k+1)}]_{T,t}^{s_r,\ldots,s_1}\ &\hbox{if}\ \ s_r=k-1\cr\cr\cr\cr
J[\psi^{(k+1)}]_{T,t}^{s_r,\ldots,s_1}+
J[\psi^{(k+1)}]_{T,t}^{k,s_r,
\ldots,s_1}/2\
&\hbox{if}\ \ s_r<k-1
\end{matrix}\ .\right.
\end{equation}

\vspace{7mm}

After substituting (\ref{30.9}) and (\ref{30.10}) into 
(\ref{30.8}) and 
regrouping 
of summands we pass to the following relations, which are valid
w.~p.~1

\vspace{1mm}
\begin{equation}
\label{30.11}
J^{*}[\psi^{(k+1)}]_{T,t}=J[\psi^{(k+1)}]_{T,t}+
\sum_{r=1}^{\left[k/2\right]}\frac{1}{2^r}
\sum_{(s_r,\ldots,s_1)\in {\rm A}_{k+1,r}}
J[\psi^{(k+1)}]_{T,t}^{s_r,\ldots,s_1} 
\end{equation}

\vspace{4mm}
\noindent
when $k$ is even and

\vspace{1mm}
\begin{equation}
\label{30.12}
J^{*}[\psi^{(k'+1)}]_{T,t}=J[\psi^{(k'+1)}]_{T,t}+
\sum_{r=1}^{\left[k'/2\right]+1}\frac{1}{2^r}
\sum_{(s_r,\ldots,s_1)\in {\rm A}_{k'+1,r}}
J[\psi^{(k'+1)}]_{T,t}^{s_r,\ldots,s_1} 
\end{equation}

\vspace{4mm}
\noindent
when $k'=k+1$ is 
uneven.

From (\ref{30.11}) and (\ref{30.12})
we have w. p. 1 

\vspace{1mm}
\begin{equation}
\label{30.13}
J^{*}[\psi^{(k+1)}]_{T,t}=J[\psi^{(k+1)}]_{T,t}+
\sum_{r=1}^{\left[(k+1)/2\right]}\frac{1}{2^r}
\sum_{(s_r,\ldots,s_1)\in {\rm A}_{k+1,r}}
J[\psi^{(k+1)}]_{T,t}^{s_r,\ldots,s_1}.
\end{equation}

\vspace{4mm}

Lemma 2 is proved. 

Consider the partition $\{\tau_j\}_{j=0}^N$ of the interval 
$[t,T]$ such that

\vspace{1mm}
\begin{equation}
\label{1111}
t=\tau_0<\ldots <\tau_N=T,\ \ \
\Delta_N=
\hbox{\vtop{\offinterlineskip\halign{
\hfil#\hfil\cr
{\rm max}\cr
$\stackrel{}{{}_{0\le j\le N-1}}$\cr
}} }\Delta\tau_j\to 0\ \ \hbox{if}\ \ N\to \infty,\ \ \ 
\Delta\tau_j=\tau_{j+1}-\tau_j.
\end{equation}

\vspace{4mm}

{\bf Lemma 3.} {\it Suppose that
every $\psi_l(\tau)$ $(l=1,\ldots, k)$ is a continuous nonrandom function on 
$[t, T]$. Then

\vspace{-3mm}
\begin{equation}
\label{30.30}
J[\psi^{(k)}]_{T,t}=
\hbox{\vtop{\offinterlineskip\halign{
\hfil#\hfil\cr
{\rm l.i.m.}\cr
$\stackrel{}{{}_{N\to \infty}}$\cr
}} }
\sum_{j_k=0}^{N-1}
\ldots \sum_{j_1=0}^{j_{2}-1}
\prod_{l=1}^k \psi_l(\tau_{j_l})\Delta{\bf w}_{\tau_
{j_l}}^{(i_l)}\ \ \ \hbox{\rm w.\ p.\ 1},
\end{equation}

\vspace{5mm}
\noindent
where $J[\psi^{(k)}]_{T,t}$ is the iterated Ito
stochastic integral {\rm (\ref{ito}),} $\Delta{\bf w}_{\tau_{j}}^{(i)}=
{\bf w}_{\tau_{j+1}}^{(i)}-{\bf w}_{\tau_{j}}^{(i)}$
$(i=0, 1,\ldots,m)$,
$\left\{\tau_{j}\right\}_{j=0}^{N}$ is the partition 
of the interval $[t,T]$ satisfying the condition {\rm (\ref{1111})}.
}

\vspace{2mm}

{\bf Proof.}\ It is easy to notice that using the 
additive property of stochastic integrals we can write the following

\vspace{-1mm}
\begin{equation}
\label{toto}
J[\psi^{(k)}]_{T,t}=
\sum_{j_k=0}^{N-1}\ldots
\sum_{j_{1}=0}^{j_{2}-1}\prod_{l=1}^{k}
J[\psi_l]_{\tau_{j_l+1},\tau_{j_l}}+
\varepsilon_N\ \ \ \ \ \hbox{w.\ p.\ 1},
\end{equation}

\vspace{2mm}
\noindent
where

\vspace{-1mm}
$$
\varepsilon_N = 
\sum_{j_k=0}^{N-1}\int\limits_{\tau_{j_k}}^{\tau_{j_k+1}}
\psi_k(s)\int\limits_{\tau_{j_k}}^{s}
\psi_{k-1}(\tau)J[\psi^{(k-2)}]_{\tau,t}d{\bf w}_\tau^{(i_{k-1})}
d{\bf w}_s^{(i_k)} +
$$

\vspace{1mm}
$$
+ \sum_{r=1}^{k-3}
G[\psi_{k-r+1}^{(k)}]_{N}
\sum_{j_{k-r}=0}^{j_{k-r+1}-1}\int\limits_{\tau_{j_{k-r}}}^{\tau_{j_{k-r}+1}}
\psi_{k-r}(s)\int\limits_{\tau_{j_{k-r}}}^{s}
\psi_{k-r-1}(\tau)J[\psi^{(k-r-2)}]_{\tau,t}
d{\bf w}_\tau^{(i_{k-r-1})}
d{\bf w}_s^{(i_{k-r})}+
$$

\vspace{1mm}
$$
+ G[\psi_3^{(k)}]_{N}
\sum_{j_{2}=0}^{j_{3}-1}J[\psi^{(2)}]_{\tau_{j_{2}+1},
\tau_{j_{2}}},
$$

\vspace{5mm}
$$
G[\psi_m^{(k)}]_{N}=\sum_{j_k=0}^{N-1}\sum_{j_{k-1}=0}^{j_k-1}
\ldots \sum_{j_{m}=0}^{j_{m+1}-1}
\prod_{l=m}^{k}J[\psi_l]_{\tau_{j_l+1},\tau_{j_l}},
$$

\vspace{3mm}
$$
J[\psi_l]_{s,\theta}=\int\limits_{\theta}^s\psi_l(\tau)d{\bf w}_{\tau}
^{(i_l)},
$$

\vspace{3mm}
$$
(\psi_m,\psi_{m+1},\ldots,\psi_{k})\stackrel{\rm def}{=}\psi_m^{(k)},\ \
(\psi_1,\ldots,\psi_{k})\stackrel{\rm def}{=}\psi_1^{(k)}=\psi^{(k)}.
$$

\vspace{6mm}

Using the standard estimates (\ref{99.010}), (\ref{99.010a}) 
for the moments of stochastic 
integrals, we obtain w. p. 1

\begin{equation}
\label{999.0001}
\hbox{\vtop{\offinterlineskip\halign{
\hfil#\hfil\cr
{\rm l.i.m.}\cr
$\stackrel{}{{}_{N\to \infty}}$\cr
}} }\varepsilon_N =0.
\end{equation}

\vspace{3mm}

Comparing (\ref{toto}) and (\ref{999.0001}), we get

\begin{equation}
\label{toto1}
J[\psi^{(k)}]_{T,t}=
\hbox{\vtop{\offinterlineskip\halign{
\hfil#\hfil\cr
{\rm l.i.m.}\cr
$\stackrel{}{{}_{N\to \infty}}
$\cr
}} }
\sum_{j_k=0}^{N-1}\ldots
\sum_{j_{1}=0}^{j_{2}-1}\prod_{l=1}^{k}
J[\psi_l]_{\tau_{j_l+1},\tau_{j_l}}\ \ \ \hbox{w.\ p.\ 1}.
\end{equation}

\vspace{3mm}

Let us write $J[\psi_l]_{\tau_{j_l+1},\tau_{j_l}}$ in the form 

$$
J[\psi_l]_{\tau_{j_l+1},\tau_{j_l}}=
\psi_l(\tau_{j_l})\Delta{\bf w}_{\tau_{j_l}}^{(i_l)}+
\int\limits_{\tau_{j_l}}^{\tau_{j_l+1}}
(\psi_l(\tau)-\psi_l(\tau_{j_l}))d{\bf w}_{\tau}^{(i_l)}\ \ \ 
\hbox{w.\ p.\ 1}
$$

\vspace{2mm}
\noindent
and substitute it into
(\ref{toto1}).
Then, due to the moment properties of stochastic integrals and
continuity (which means uniform continuity) of the functions 
$\psi_l(s)$ ($l=1,\ldots,k$)
it is easy to see that the 
prelimit
expression on the right-hand side of (\ref{toto1}) is a sum of 
the prelimit
expression on the right-hand side of (\ref{30.30}) and the value which 
tends to zero in the mean-square sense if 
$N\to\infty.$ Lemma 3 is proved. 

\vspace{2mm}

{\bf Remark 1.} {\it It is easy to see that if
$\Delta{\bf w}_{\tau_{j_l}}^{(i_l)}$ in {\rm (\ref{30.30})}
for some $l\in\{1,\ldots,k\}$ is replaced with 
$\left(\Delta{\bf w}_{\tau_{j_l}}^{(i_l)}\right)^p$ $(p=2,$
$i_l\ne 0),$ then
the differential $d{\bf w}_{t_{l}}^{(i_l)}$
in the integral $J[\psi^{(k)}]_{T,t}$
will be replaced with $dt_l$.
If $p=3, 4,\ldots,$ then the
right-hand side 
of the formula {\rm (\ref{30.30})}
will become zero w.~p.~{\rm 1}.
If we replace $\Delta{\bf w}_{\tau_{j_l}}^{(i_l)}$ in {\rm (\ref{30.30})}
for some $l\in\{1,\ldots,k\}$
with $\left(\Delta \tau_{j_l}\right)^p$ $(p=2, 3,\ldots),$
then the right-hand side of the formula
{\rm (\ref{30.30})} also 
will be equal to zero w.~p.~{\rm 1}.}

\vspace{2mm}

Let us define the following
multiple stochastic integral

\vspace{-1mm}
\begin{equation}
\label{30.34}
\hbox{\vtop{\offinterlineskip\halign{
\hfil#\hfil\cr
{\rm l.i.m.}\cr
$\stackrel{}{{}_{N\to \infty}}$\cr
}} }\sum_{j_1,\ldots,j_k=0}^{N-1}
\Phi\left(\tau_{j_1},\ldots,\tau_{j_k}\right)
\prod\limits_{l=1}^k\Delta{\bf w}_{\tau_{j_l}}^{(i_l)}
\stackrel{\rm def}{=}J[\Phi]_{T,t}^{(k)},
\end{equation}

\vspace{3mm}
\noindent
where $\Phi(t_1,\ldots,t_k):\ [t, T]^k\to\mathbb{R}$ is a nonrandom function 
(the properties of this function
will be specified further).

Denote

\vspace{-3mm}
\begin{equation}
\label{dom1}
D_k=\{(t_1,\ldots,t_k):\ t\le t_1<\ldots <t_k\le T\}.
\end{equation}

\vspace{4mm}

We will use the same symbol $D_k$ to denote the open and closed 
domains corresponding to the domain $D_k$ defined by (\ref{dom1}).
However, we always specify what domain we consider (open or closed). 
Also we will write $\Phi(t_1,\ldots,t_k)\in C(D_k)$
if $\Phi(t_1,\ldots,t_k)$ is a continuous nonrandom function of $k$ variables
in the closed domain $D_k$.

Let us consider the iterated Ito stochastic integral

\begin{equation}
\label{rrr29}
I[\Phi]_{T,t}^{(k)}\stackrel{\rm def}{=}
\int\limits_t^T\ldots \int\limits_t^{t_2}
\Phi(t_1,\ldots,t_k)d{\bf w}_{t_1}^{(i_1)}\ldots
d{\bf w}_{t_k}^{(i_k)},
\end{equation}

\vspace{3mm}
\noindent
where $\Phi(t_1,\ldots,t_k)\in C(D_k).$

Using the arguments which similar to the arguments used in the 
proof of Lemma 3
it is easy to demonstrate that if
$\Phi(t_1,\ldots,t_k)\in C(D_k),$ then the following equality is fulfilled

\begin{equation}
\label{30.52}
I[\Phi]_{T,t}^{(k)}=\hbox{\vtop{\offinterlineskip\halign{
\hfil#\hfil\cr
{\rm l.i.m.}\cr
$\stackrel{}{{}_{N\to \infty}}$\cr
}} }
\sum_{j_k=0}^{N-1}
\ldots \sum_{j_1=0}^{j_{2}-1}
\Phi(\tau_{j_1},\ldots,\tau_{j_k})
\prod\limits_{l=1}^k\Delta {\bf w}_{\tau_{j_l}}^{(i_l)}\ \ \ \hbox{w.\ p.\ 1}.
\end{equation}

\vspace{3mm}

In order to explain this, let us check the correctness of the equality 
(\ref{30.52}) when $k=3$.
For definiteness we will suppose that 
$i_1,i_2,i_3=1,\ldots,m.$ We have

$$
I[\Phi]_{T,t}^{(3)}\stackrel{\rm def}{=}
\int\limits_t^T\int\limits_t^{t_3}\int\limits_t^{t_2}
\Phi(t_1,t_2,t_3)d{\bf w}_{t_1}^{(i_1)}d{\bf w}_{t_2}^{(i_2)}
d{\bf w}_{t_3}^{(i_3)}=
$$

\vspace{1mm}
$$
=\hbox{\vtop{\offinterlineskip\halign{
\hfil#\hfil\cr
{\rm l.i.m.}\cr
$\stackrel{}{{}_{N\to \infty}}$\cr
}} }
\sum_{j_3=0}^{N-1}
\int\limits_{t}^{\tau_{j_3}}\int\limits_t^{t_2}
\Phi(t_1,t_2,\tau_{j_3})d{\bf w}_{t_1}^{(i_1)}d{\bf w}_{t_2}^{(i_2)}
\Delta{\bf w}_{\tau_{j_3}}^{(i_3)}=
$$

\vspace{1mm}
$$
=\hbox{\vtop{\offinterlineskip\halign{
\hfil#\hfil\cr
{\rm l.i.m.}\cr
$\stackrel{}{{}_{N\to \infty}}$\cr
}} }
\sum_{j_3=0}^{N-1}\sum_{j_2=0}^{j_3-1}
\int\limits_{\tau_{j_2}}^{\tau_{j_2+1}}\int\limits_t^{t_2}
\Phi(t_1,t_2,\tau_{j_3})d{\bf w}_{t_1}^{(i_1)}d{\bf w}_{t_2}^{(i_2)}
\Delta{\bf w}_{\tau_{j_3}}^{(i_3)}=
$$

\vspace{1mm}
$$
=\hbox{\vtop{\offinterlineskip\halign{
\hfil#\hfil\cr
{\rm l.i.m.}\cr
$\stackrel{}{{}_{N\to \infty}}$\cr
}} }
\sum_{j_3=0}^{N-1}\sum_{j_2=0}^{j_3-1}
\int\limits_{\tau_{j_2}}^{\tau_{j_2+1}}
\left(\ \int\limits_t^{\tau_{j_2}}\
+\ \int\limits_{\tau_{j_2}}^{t_2}\ \right)
\Phi(t_1,t_2,\tau_{j_3})d{\bf w}_{t_1}^{(i_1)}d{\bf w}_{t_2}^{(i_2)}
\Delta{\bf w}_{\tau_{j_3}}^{(i_3)}=
$$

\vspace{1mm}
$$
=\hbox{\vtop{\offinterlineskip\halign{
\hfil#\hfil\cr
{\rm l.i.m.}\cr
$\stackrel{}{{}_{N\to \infty}}$\cr
}} }
\sum_{j_3=0}^{N-1}\sum_{j_2=0}^{j_3-1}\sum_{j_1=0}^{j_2-1}
\int\limits_{\tau_{j_2}}^{\tau_{j_2+1}}\int\limits_{\tau_{j_1}}^{\tau_{j_1+1}}
\Phi(t_1,t_2,\tau_{j_3})d{\bf w}_{t_1}^{(i_1)}d{\bf w}_{t_2}^{(i_2)}
\Delta{\bf w}_{\tau_{j_3}}^{(i_3)}+
$$

\vspace{1mm}
\begin{equation}
\label{44444.25}
+\hbox{\vtop{\offinterlineskip\halign{
\hfil#\hfil\cr
{\rm l.i.m.}\cr
$\stackrel{}{{}_{N\to \infty}}$\cr
}} }
\sum_{j_3=0}^{N-1}\sum_{j_2=0}^{j_3-1}
\int\limits_{\tau_{j_2}}^{\tau_{j_2+1}}\int\limits_{\tau_{j_2}}^{t_2}
\Phi(t_1,t_2,\tau_{j_3})d{\bf w}_{t_1}^{(i_1)}d{\bf w}_{t_2}^{(i_2)}
\Delta{\bf w}_{\tau_{j_3}}^{(i_3)}.
\end{equation}

\vspace{5mm}

Let us demonstrate that the second limit on the right-hand side 
of (\ref{44444.25}) 
equals to zero.
Actually, for the second moment of its 
prelimit
expression we get

$$
\sum_{j_3=0}^{N-1}\sum_{j_2=0}^{j_3-1}
\int\limits_{\tau_{j_2}}^{\tau_{j_2+1}}\int\limits_{\tau_{j_2}}^{t_2}
\Phi^2(t_1,t_2,\tau_{j_3})dt_1 dt_2
\Delta\tau_{j_3}
\le M^2  \sum_{j_3=0}^{N-1}\sum_{j_2=0}^{j_3-1}
\frac{1}{2}\left(\Delta\tau_{j_2}\right)^2\Delta\tau_{j_3}\to 0
$$

\vspace{3mm}
\noindent
when $N\to\infty.$
Here $M$ is a constant, which restricts the module of the
function
$\Phi(t_1,t_2,t_3)$ due to its continuity, $\Delta\tau_j=
\tau_{j+1}-\tau_j.$

Considering the obtained conclusions, we have

$$
I[\Phi]_{T,t}^{(3)}\stackrel{\rm def}{=}
\int\limits_t^T\int\limits_t^{t_3}\int\limits_t^{t_2}
\Phi(t_1,t_2,t_3)d{\bf w}_{t_1}^{(i_1)}d{\bf w}_{t_2}^{(i_2)}
d{\bf w}_{t_3}^{(i_3)}=
$$

\vspace{1mm}
$$
=\hbox{\vtop{\offinterlineskip\halign{
\hfil#\hfil\cr
{\rm l.i.m.}\cr
$\stackrel{}{{}_{N\to \infty}}$\cr
}} }
\sum_{j_3=0}^{N-1}\sum_{j_2=0}^{j_3-1}\sum_{j_1=0}^{j_2-1}
\int\limits_{\tau_{j_2}}^{\tau_{j_2+1}}\int\limits_{\tau_{j_1}}^{\tau_{j_1+1}}
\Phi(t_1,t_2,\tau_{j_3})d{\bf w}_{t_1}^{(i_1)}d{\bf w}_{t_2}^{(i_2)}
\Delta{\bf w}_{\tau_{j_3}}^{(i_3)}=
$$

\vspace{1mm}
$$
=\hbox{\vtop{\offinterlineskip\halign{
\hfil#\hfil\cr
{\rm l.i.m.}\cr
$\stackrel{}{{}_{N\to \infty}}$\cr
}} }
\sum_{j_3=0}^{N-1}\sum_{j_2=0}^{j_3-1}\sum_{j_1=0}^{j_2-1}
\int\limits_{\tau_{j_2}}^{\tau_{j_2+1}}\int\limits_{\tau_{j_1}}^{\tau_{j_1+1}}
\left(\Phi(t_1,t_2,\tau_{j_3})-\Phi(t_1,\tau_{j_2},\tau_{j_3})\right)
d{\bf w}_{t_1}^{(i_1)}d{\bf w}_{t_2}^{(i_2)}
\Delta{\bf w}_{\tau_{j_3}}^{(i_3)}+
$$

\vspace{1mm}
$$
+\hbox{\vtop{\offinterlineskip\halign{
\hfil#\hfil\cr
{\rm l.i.m.}\cr
$\stackrel{}{{}_{N\to \infty}}$\cr
}} }
\sum_{j_3=0}^{N-1}\sum_{j_2=0}^{j_3-1}\sum_{j_1=0}^{j_2-1}
\int\limits_{\tau_{j_2}}^{\tau_{j_2+1}}\int\limits_{\tau_{j_1}}^{\tau_{j_1+1}}
\left(\Phi(t_1,\tau_{j_2},\tau_{j_3})-
\Phi(\tau_{j_1},\tau_{j_2},\tau_{j_3})\right)
d{\bf w}_{t_1}^{(i_1)}d{\bf w}_{t_2}^{(i_2)}
\Delta{\bf w}_{\tau_{j_3}}^{(i_3)}+
$$

\vspace{1mm}
\begin{equation}
\label{4444.1}
+\hbox{\vtop{\offinterlineskip\halign{
\hfil#\hfil\cr
{\rm l.i.m.}\cr
$\stackrel{}{{}_{N\to \infty}}$\cr
}} }
\sum_{j_3=0}^{N-1}\sum_{j_2=0}^{j_3-1}\sum_{j_1=0}^{j_2-1}
\Phi(\tau_{j_1},\tau_{j_2},\tau_{j_3})
\Delta{\bf w}_{\tau_{j_1}}^{(i_1)}
\Delta{\bf w}_{\tau_{j_2}}^{(i_2)}
\Delta{\bf w}_{\tau_{j_3}}^{(i_3)}.
\end{equation}

\vspace{5mm}

In order to get the sought result, we just have to demonstrate that 
the first
two limits on the right-hand side of (\ref{4444.1}) equal to zero. 
Let us prove 
that the first one of them equals to zero (proof for the second limit 
is similar).
   
The second moment of prelimit expression of the first limit on the 
right-hand side of (\ref{4444.1}) equals to the following expression

\begin{equation}
\label{4444.01}
\sum_{j_3=0}^{N-1}\sum_{j_2=0}^{j_3-1}\sum_{j_1=0}^{j_2-1}
\int\limits_{\tau_{j_2}}^{\tau_{j_2+1}}\int\limits_{\tau_{j_1}}^{\tau_{j_1+1}}
\left(\Phi(t_1,t_2,\tau_{j_3})-\Phi(t_1,\tau_{j_2},\tau_{j_3})\right)^2
dt_1 dt_2
\Delta\tau_{j_3}.
\end{equation}

\vspace{3mm}

Since the function $\Phi(t_1,t_2,t_3)$ is continuous in 
the closed bo\-un\-ded domain
$D_3,$ then
it is uniformly continuous in this domain. Therefore, if the 
distance between two points of the domain $D_3$ is less than 
$\delta(\varepsilon)$ ($\delta(\varepsilon)>0$ 
exists
for any $\varepsilon>0$ and it does not depend 
on mentioned points), then the cor\-res\-pond\-ing oscillation of the function 
$\Phi(t_1,t_2,t_3)$ for these two points of the domain $D_3$ is less than
$\varepsilon.$

If we assume that $\Delta\tau_j<\delta(\varepsilon)$ ($j=0, 1,\ldots,N-1$),
then the distance between points 
$(t_1,t_2,\tau_{j_3})$,\ $(t_1,\tau_{j_2},\tau_{j_3})$
is obviously less than $\delta(\varepsilon).$ In this case 

\vspace{-1mm}
$$
|\Phi(t_1,t_2,\tau_{j_3})-\Phi(t_1,\tau_{j_2},\tau_{j_3})|<\varepsilon.
$$

\vspace{3mm}

Consequently, when $\Delta\tau_j<\delta(\varepsilon)$ ($j=0,\ 1,\ldots,N-1$)
the expression (\ref{4444.01})  
is estimated by the following value

\vspace{-1mm}
$$
\varepsilon^2
\sum_{j_3=0}^{N-1}\sum_{j_2=0}^{j_3-1}\sum_{j_1=0}^{j_2-1}
\Delta\tau_{j_1}\Delta\tau_{j_2}\Delta\tau_{j_3}<
\varepsilon^2\frac{(T-t)^3}{6}.
$$ 

\vspace{4mm}

Therefore, the first limit on the right-hand side 
of (\ref{4444.1}) equals to zero.
Similarly, we can prove that the second limit on the right-hand
side
of (\ref{4444.1}) equals to zero.

Consequently, the equality (\ref{30.52}) is proved for $k=3$. 
The cases $k=2$ and $k>3$ 
are analyzed absolutely similarly.

It is necessary to note that the proof of 
correctness of (\ref{30.52}) 
is similar when the nonrandom function $\Phi(t_1,\ldots,t_k)$ is 
continuous in 
the open domain $D_k$ and bounded at its boundary.

Let us consider the class ${\rm M}_2([0, T])$ of functions
$\xi: [0,T]\times\Omega\rightarrow \mathbb{R},$ which are 
measurable with respect to the variables
$(t,\omega)$ and
${\rm F}_t$-measurable 
for all $t\in[0,T].$ Moreover, $\xi(\tau,\omega)$ is independent 
with increments ${\bf f}_{t+\Delta}-{\bf f}_{t}$
for $t\ge \tau $ $(\Delta>0),$ 

\vspace{-1mm}
$$
\int\limits_0^T{\sf M}\left\{\xi^2(t,\omega)\right\}dt
<\infty,
$$

\vspace{2mm}
\noindent
and ${\sf M}\left\{\xi^2(t,\omega)\right\}<\infty$
for all $t\in[0,T].$

It is well-known \cite{1}, \cite{100} that the Ito stochastic integral
exists in the mean-square sense for any 
$\xi\in{\rm M}_2([0, T]).$  Further, we will denote  
$\xi(\tau,\omega)$ as $\xi_{\tau}.$

\vspace{2mm}

{\bf Lemma 4.}\ {\it Suppose that $\Phi(t_1,\ldots,t_k)\in C(D_k)$ 
or $\Phi(t_1,\ldots,t_k)$ 
is a continuous nonrandom function in the open domain $D_k$ and bounded at its boundary.
Then

$$
{\sf M}\left\{\biggl|I[\Phi]_{T,t}^{(k)}\biggr|^{2n}\right\}
\le C_{k}
\int\limits_t^T\ldots \int\limits_t^{t_2}
\Phi^{2n}(t_1,\ldots,t_k)dt_1\ldots dt_k,\ \ \
C_{k}<\infty,
$$

\vspace{3mm}
\noindent
where $I[\Phi]_{T,t}^{(k)}$ is defined by the formula {\rm (\ref{rrr29})}.}

\vspace{2mm}

{\bf Proof.}\
Using standard estimates for moments of stochastic integrals,
we have \cite{100}

\begin{equation}
\label{99.010}
{\sf M}\left\{\left|\int\limits_{t}^T \xi_\tau
df_\tau\right|^{2n}\right\} \le (T-t)^{n-1}\left(n(2n-1)\right)^n
\int\limits_{t}^T {\sf M}\left\{\left|\xi_\tau \right|^{2n}\right\}d\tau,
\end{equation}
\begin{equation}
\label{99.010a}
{\sf M}\left\{\left|\int\limits_{t}^T \xi_\tau
d\tau\right|^{2n}\right\} \le (T-t)^{2n-1}
\int\limits_{t}^T {\sf M}\left\{\left|\xi_\tau\right|^{2n}\right\}d\tau,
\end{equation}

\vspace{4mm}
\noindent
where the process $\xi_{\tau}$ is such that
$\left(\xi_{\tau}\right)^n\in{\rm M}_2
([t,T])$ and $f_t$ is a scalar standard Wiener 
process,\
$n=1, 2,\ldots$ 

\vspace{3mm}

Let us denote

\vspace{-2mm}
$$
\xi[\Phi]_{t_{l+1},\ldots,t_k,t}^{(l)}=
\int\limits_t^{t_{l+1}}\ldots \int\limits_t^{t_2}
\Phi(t_1,\ldots,t_k)
d{\bf w}_{t_1}^{(i_1)}\ldots
d{\bf w}_{t_{l}}^{(i_{l})},
$$ 

\vspace{3mm}
\noindent
where $l=1,\ldots,$ $k-1$ and
$\xi[\Phi]_{t_{1},\ldots,t_k,t}^{(0)}\stackrel{\rm def}{=}
\Phi(t_1,\ldots,t_k).$

By induction it is easy to demonstrate that
$\left(\xi[\Phi]_{t_{l+1},\ldots,t_k,t}^{(l)}\right)^n\in{\rm M}_2([t,T])$
with respect to the variable $t_{l+1}.$
Further, using the estimates (\ref{99.010}) and (\ref{99.010a})
repeatedly we obtain the statement of Lemma 4. Lemma 4 is proved.

\vspace{2mm}

{\bf Lemma 5} \cite{1997} (1997), 
\cite{1a}, \cite{5b}-\cite{2017a}, \cite{2017}-\cite{12aa-afterxxx}.
{\it Suppose that every $\varphi_l(s)$
$(l=1,\ldots,k)$ is a continuous nonrandom function on $[t, T]$.
Then

\vspace{-1mm}
\begin{equation}
\label{30.39}
\prod_{l=1}^k 
J[\varphi_l]_{T,t}=J[\Phi]_{T,t}^{(k)}\ \ \ \hbox{\rm w.\ p.\ 1},
\end{equation}

\vspace{3mm}
\noindent
where 
$$
J[\varphi_l]_{T,t}
=\int\limits_t^T \varphi_l(s) d{\bf w}_{s}^{(i_l)},\ \ \ 
\Phi(t_1,\ldots,t_k)=\prod\limits_{l=1}^k\varphi_l(t_l),
$$

\vspace{3mm}
\noindent
and the integral $J[\Phi]_{T,t}^{(k)}$ 
is defined
by the equality
{\rm (\ref{30.34})}.
}

\vspace{2mm}

{\bf Proof.}\
Let at first $i_l\ne 0$ $(l=1,\ldots,k).$
Denote

\vspace{2mm}
$$
J[\varphi_l]_{N}\stackrel{\rm def}{=}\sum\limits_{j=0}^{N-1}
\varphi_l(\tau_j)\Delta{\bf w}_{\tau_j}^{(i_l)}.
$$

\vspace{2mm}

Since
$$
\prod_{l=1}^k J[\varphi_l]_{N}-\prod_{l=1}^k J[\varphi_l]_{T,t}
=
$$

\vspace{1mm}
$$
=\sum_{l=1}^k \left(\prod_{g=1}^{l-1} J[\varphi_g]_{T,t}\right)
\biggl(J[\varphi_l]_{N}-J[\varphi_l]_{T,t}
\biggr)\left(\prod_{g=l+1}^k J[\varphi_g]_{N}\right),
$$

\vspace{6mm}
\noindent
then due to the Minkowski inequality and the inequality 
of Cauchy-Bunyakovsky we obtain

\vspace{-1mm}
\begin{equation}
\label{30.42}
\left({\sf M}\left\{\left|\prod_{l=1}^k J[\varphi_l]_{N}
-\prod_{l=1}^k J[\varphi_l]_{T,t}\right|^2
\right\}\right)^{1/2}\le C_k
\sum_{l=1}^k
\left({\sf M}\left\{
\biggl|J[\varphi_l]_{N}-J[\varphi_l]_{T,t}\biggr|^4\right\}\right)
^{1/4},
\end{equation}

\vspace{5mm}
\noindent 
where $C_k$ is a constant.

Note that

$$
J[\varphi_l]_{N}-J[\varphi_l]_{T,t}=\sum\limits_{j=0}^{N-1}
J[\Delta\varphi_l]_{\tau_{j+1},\tau_j},\ \ \
J[\Delta\varphi_l]_{\tau_{j+1},\tau_j}
=\int\limits_{\tau_j}^{\tau_{j+1}}\left(
\varphi_l(\tau_j)-\varphi_l(s)\right)d{\bf w}_{s}^{(i_l)}.
$$

\vspace{3mm}

Since $J[\Delta\varphi_l]_{\tau_{j+1},\tau_j}$
are independent for various $j,$ then \cite{Sc}

\vspace{1mm}
$$
{\sf M}\left\{\left|\sum_{j=0}^{N-1}J[\Delta\varphi_l]_{\tau_{j+1},
\tau_j}\right|^4
\right\}=
\sum_{j=0}^{N-1}{\sf M}\left\{\biggl|J[\Delta\varphi_l]_{\tau_{j+1},
\tau_j}\biggr|^4
\right\}+ 
$$

\vspace{1mm}
\begin{equation}
\label{30.43}
+6 \sum_{j=0}^{N-1}{\sf M}
\left\{\biggl|J[\Delta\varphi_l]_{\tau_{j+1},\tau_j}\biggr|^2
\right\}
\sum_{q=0}^{j-1}{\sf M}\left\{\biggl|
J[\Delta\varphi_l]_{\tau_{q+1},\tau_q}\biggr|^2
\right\}.
\end{equation}

\vspace{6mm}

Moreover, since
$J[\Delta\varphi_l]_{\tau_{j+1},\tau_j}$ is a Gaussian random variable,
we have

\vspace{2mm}
$$
{\sf M}\left\{\biggl|J[\Delta\varphi_l]_{\tau_{j+1},\tau_j}\biggr|^2\right\}=
\int\limits_{\tau_j}^{\tau_{j+1}}(\varphi_l(\tau_j)-\varphi_l(s))^2ds,
$$

\vspace{2mm}
$$
{\sf M}\left\{\biggl|J[\Delta\varphi_l]_{\tau_{j+1},\tau_j}\biggr|^4\right\}=
3\left(\int\limits_{\tau_j}^{\tau_{j+1}}(\varphi_l(\tau_j)-\varphi_l(s))^2ds
\right)^2.
$$

\vspace{6mm}

Using these relations and continuity (which means uniform continuity)
of the functions $\varphi_l(s),$ we obtain

\vspace{-2mm}
$$
{\sf M}\left\{\left|\sum_{j=0}^{N-1}J[\Delta\varphi_l]_{\tau_{j+1},
\tau_j}\right|^4
\right\}\le 
$$

\vspace{2mm}
$$
\le \varepsilon^4\left(
3 \sum_{j=0}^{N-1}(\Delta\tau_{j})^2+
6 \sum_{j=0}^{N-1}\Delta\tau_{j}
\sum_{q=0}^{j-1}\Delta\tau_{q}\right)
<3\varepsilon^4\left(\delta(\varepsilon) (T-t)+(T-t)^2\right),
$$

\vspace{6mm}
\noindent
where $\Delta\tau_{j}<\delta(\varepsilon),$ $j=0,1,\ldots,N-1$ ($\forall$
$\varepsilon>0$\ $\exists$ 
$\delta(\varepsilon)>0$ which does not
depend on points of the interval $[t, T]$ and such that 
$|\varphi_l(\tau_j)-\varphi_l(s)|<\varepsilon,$ $s\in [\tau_j, \tau_{j+1}]$).
Then the right-hand side of the formula 
(\ref{30.43}) tends to zero when $N\to \infty.$ 

Taking into account this fact as well 
as (\ref{30.42}), we obtain (\ref{30.39}). 
If ${\bf w}_{t_l}^{(i_l)}=t_l$ for some $l\in\{1,\ldots,k\},$
then
the proof of Lemma 5 becomes obviously simpler and  
it is performed similarly. Lemma 5 is proved.

Using Lemma 2 and (\ref{30.52}), we obtain w. p. 1 

\begin{equation}
\label{pp1}
J^{*}[\psi^{(k)}]_{T,t}=
J[\psi^{(k)}]_{T,t}+
\sum_{r=1}^{\left[k/2\right]}\frac{1}{2^r}
\sum_{(s_r,\ldots,s_1)\in {\rm A}_{k,r}}
J[\psi^{(k)}]_{T,t}^{s_r,\ldots,s_1}=
J[{K^{*}}]_{T,t}^{(k)},
\end{equation}

\vspace{3mm}
\noindent
where the stochastic integral
$J[K^{*}]_{T,t}^{(k)}$
is defined in accordance with (\ref{30.34}).

Let us subsitute the relation

$$
K^{*}(t_1,\ldots,t_k)=
$$

\vspace{1mm}
$$
=
\sum_{j_1=0}^{p_1}\ldots\sum_{j_k=0}^{p_k}
C_{j_k\ldots j_1} \prod_{l=1}^{k} \phi_{j_l}(t_l)
+K^{*}(t_1,\ldots,t_k)-\sum_{j_1=0}^{p_1}\ldots\sum_{j_k=0}^{p_k}
C_{j_k\ldots j_1} \prod_{l=1}^{k} \phi_{j_l}(t_l)
$$

\vspace{6mm}
\noindent
into (\ref{pp1}) (here we suppose that $p_1,\ldots,p_k<\infty).$

Then using Lemma 5, we obtain

\vspace{1mm}
\begin{equation}
\label{proof1}
J^{*}[\psi^{(k)}]_{T,t}=
\sum_{j_1=0}^{p_1}\ldots\sum_{j_k=0}^{p_k}
C_{j_k\ldots j_1}
\prod_{l=1}^k \zeta_{j_l}^{(i_l)}+
J[R_{p_1\ldots p_k}]_{T,t}^{(k)}\ \ \ \hbox{w. p. 1,}
\end{equation}

\vspace{4mm}
\noindent
where the stochastic integral
$J[R_{p_1\ldots p_k}]_{T,t}^{(k)}$
is defined in accordance with (\ref{30.34}) and

\vspace{1mm}
\begin{equation}
\label{30.46}
R_{p_1\ldots p_k}(t_1,\ldots,t_k)=
K^{*}(t_1,\ldots,t_k)-
\sum_{j_1=0}^{p_1}\ldots\sum_{j_k=0}^{p_k}
C_{j_k\ldots j_1} \prod_{l=1}^{k} \phi_{j_l}(t_l),
\end{equation}

\vspace{2mm}

$$
\zeta_{j_l}^{(i_l)}=\int\limits_t^T \phi_{j_l}(s) d{\bf w}_s^{(i_l)}.
$$

\vspace{4mm}

According to Lemma 1, we obtain

\begin{equation}
\label{410}
\hbox{\vtop{\offinterlineskip\halign{
\hfil#\hfil\cr
{\rm lim}\cr
$\stackrel{}{{}_{p_1\to \infty}}$\cr
}} }\ldots \hbox{\vtop{\offinterlineskip\halign{
\hfil#\hfil\cr
{\rm lim}\cr
$\stackrel{}{{}_{p_k\to \infty}}$\cr
}} }R_{p_1\ldots p_k}(t_1,\ldots,t_k)=0\ \ \ \hbox{when}\ \ \ (t_1,\ldots,t_k)\in (t,T)^k,
\end{equation}

\vspace{4mm}
\noindent
where the left-hand side of (\ref{410})
is bounded on $[t, T]^k.$

\vspace{2mm}

{\bf Lemma 6.} {\it Under the conditions of Theorem {\rm 1} the following equality is correct

\vspace{1mm}
$$
\lim\limits_{p_1\to \infty}
\varlimsup\limits_{p_2\to \infty}
\ldots
\varlimsup\limits_{p_k\to \infty}
{\sf M}\left\{\left|J[R_{p_1\ldots p_k}]_{T,t}^{(k)}\right|^{2n}
\right\}=0,\ \ \ n\in \mathbb{N}.
$$
}

\vspace{4mm}

{\bf Proof.} At first let us analize in detail the cases $k=2, 3, 4.$ 
Using (\ref{best1}) (see below), we have w.~p.~1

\vspace{3mm}
$$
J[R_{p_1p_2}]_{T,t}^{(2)}
=\hbox{\vtop{\offinterlineskip\halign{
\hfil#\hfil\cr
{\rm l.i.m.}\cr
$\stackrel{}{{}_{N\to \infty}}$\cr
}} }\sum_{l_2=0}^{N-1}
\sum_{l_1=0}^{N-1}
R_{p_1 p_2}(\tau_{l_1},\tau_{l_2})
\Delta{\bf w}_{\tau_{l_1}}^{(i_1)}
\Delta{\bf w}_{\tau_{l_2}}^{(i_2)}=
$$

\vspace{1mm}
$$
=\hbox{\vtop{\offinterlineskip\halign{
\hfil#\hfil\cr
{\rm l.i.m.}\cr
$\stackrel{}{{}_{N\to \infty}}$\cr
}} }\sum_{l_2=0}^{N-1}
\sum_{l_1=0}^{l_2-1}
R_{p_1 p_2}(\tau_{l_1},\tau_{l_2})
\Delta{\bf w}_{\tau_{l_1}}^{(i_1)}
\Delta{\bf w}_{\tau_{l_2}}^{(i_2)}+
\hbox{\vtop{\offinterlineskip\halign{
\hfil#\hfil\cr
{\rm l.i.m.}\cr
$\stackrel{}{{}_{N\to \infty}}$\cr
}} }\sum_{l_1=0}^{N-1}
\sum_{l_2=0}^{l_1-1}
R_{p_1 p_2}(\tau_{l_1},\tau_{l_2})
\Delta{\bf w}_{\tau_{l_1}}^{(i_1)}
\Delta{\bf w}_{\tau_{l_2}}^{(i_2)}+
$$

\vspace{1mm}
$$
+
\hbox{\vtop{\offinterlineskip\halign{
\hfil#\hfil\cr
{\rm l.i.m.}\cr
$\stackrel{}{{}_{N\to \infty}}$\cr
}} }\sum_{l_1=0}^{N-1}
R_{p_1 p_2}(\tau_{l_1},\tau_{l_1})
\Delta{\bf w}_{\tau_{l_1}}^{(i_1)}
\Delta{\bf w}_{\tau_{l_1}}^{(i_2)}=
$$

\vspace{1mm}
$$
=\int\limits_t^T\int\limits_t^{t_2}
R_{p_1p_2}(t_1,t_2)d{\bf w}_{t_1}^{(i_1)}d{\bf w}_{t_2}^{(i_2)}
+\int\limits_t^T\int\limits_t^{t_1}
R_{p_1p_2}(t_1,t_2)d{\bf w}_{t_2}^{(i_2)}d{\bf w}_{t_1}^{(i_1)}+
$$

\begin{equation}
\label{nov800}
+{\bf 1}_{\{i_1=i_2\ne 0\}}
\int\limits_t^T R_{p_1p_2}(t_1,t_1)dt_1,
\end{equation}

\vspace{3mm}
\noindent
where

\vspace{-3mm}
\begin{equation}
\label{zajaz}
R_{p_1 p_2}(t_1,t_2)=
K^{*}(t_1,t_2)-
\sum_{j_1=0}^{p_1}\sum_{j_2=0}^{p_2}C_{j_2 j_1}
\phi_{j_1}(t_1)\phi_{j_2}(t_2),\ \ \ p_1,\ p_2<\infty.
\end{equation}

\vspace{4mm}

Using Lemma 4, we obtain 

\vspace{1mm}
$$
{\sf M}\left\{\left|J[R_{p_1p_2}]_{T,t}^{(2)}\right|^{2n}
\right\}\le C_n\left(\int\limits_t^T\int\limits_t^{t_2}
\left(R_{p_1p_2}(t_1,t_2)\right)^{2n}dt_1 dt_2
+\right.
$$

\vspace{1mm}
\begin{equation}
\label{leto80010}
\left.+
\int\limits_t^T\int\limits_t^{t_1}
\left(R_{p_1p_2}(t_1,t_2)\right)^{2n}dt_2 dt_1+
{\bf 1}_{\{i_1=i_2\ne 0\}}
\int\limits_t^T \left(R_{p_1p_2}(t_1,t_1)\right)^{2n}dt_1\right),
\end{equation}

\vspace{4mm}
\noindent
where constant $C_n<\infty$ depends on 
$n$ and $T-t$ $(n=1, 2,\ldots).$

Further, we have

\vspace{1mm}
$$
\int\limits_t^T\int\limits_t^{t_2}
\left(R_{p_1p_2}(t_1,t_2)\right)^{2n}dt_1 dt_2
+
\int\limits_t^T\int\limits_t^{t_1}
\left(R_{p_1p_2}(t_1,t_2)\right)^{2n}dt_2 dt_1=
$$

\vspace{1mm}
\begin{equation}
\label{leb1}
=
\int\limits_t^T\int\limits_t^{t_2}
\left(R_{p_1p_2}(t_1,t_2)\right)^{2n}dt_1 dt_2
+
\int\limits_t^T\int\limits_{t_2}^{T}
\left(R_{p_1p_2}(t_1,t_2)\right)^{2n}dt_1 dt_2
=
\int\limits_{[t, T]^2}
\left(R_{p_1p_2}(t_1,t_2)\right)^{2n}dt_1 dt_2.
\end{equation}

\vspace{5mm}

Combining (\ref{leto80010}) and (\ref{leb1}), we obtain

\vspace{1mm}
$$
{\sf M}\left\{\left|J[R_{p_1p_2}]_{T,t}^{(2)}\right|^{2n}
\right\}\le 
$$

\vspace{1mm}

\begin{equation}
\label{leto80010aaa}
\le C_n\left(
\int\limits_{[t, T]^2}
\left(R_{p_1p_2}(t_1,t_2)\right)^{2n}dt_1 dt_2
+{\bf 1}_{\{i_1=i_2\ne 0\}}
\int\limits_t^T \left(R_{p_1p_2}(t_1,t_1)\right)^{2n}dt_1\right),
\end{equation}

\vspace{4mm}
\noindent
where constant $C_n<\infty$ depends on 
$n$ and $T-t$ $(n=1, 2,\ldots).$

Since the integrals on the right-hand side of (\ref{leto80010aaa}) 
exist as Riemann integrals, then they are equal to the 
corresponding Lebesgue integrals. 
Moreover,

\vspace{1mm}
\begin{equation}
\label{strange101}
\lim\limits_{p_1\to\infty}\lim\limits_{p_2\to\infty}
\left(R_{p_1p_2}(t_1,t_2)\right)^{2n}=0\ \ \ \hbox{when}\ \ \ (t_1,t_2)\in (t, T)^2,
\end{equation}

\vspace{4mm}
\noindent
where $n\in \mathbb{N}$, the left-hand side
is bounded on $[t, T]^2.$

According to (\ref{leto8001})--(\ref{leto8003}) and (\ref{zajaz}), we obtain

$$
R_{p_1p_2}(t_1,t_2)=\left(K^{*}(t_1,t_2)-\sum\limits_{j_1=0}^{p_1}
C_{j_1}(t_2)\phi_{j_1}(t_1)\right)+
$$

\vspace{1mm}
\begin{equation}
\label{d2020}
+\left(
\sum\limits_{j_1=0}^{p_1}\left(C_{j_1}(t_2)-
\sum\limits_{j_2=0}^{p_2}
C_{j_2j_1}\phi_{j_2}(t_2)\right)
\phi_{j_1}(t_1)\right).
\end{equation}

\vspace{5mm}

Then, applying two times (we mean an iterated passage to the limit
$\lim\limits_{p_1\to\infty}\varlimsup\limits_{p_2\to\infty}$)
the Lebesgue's Dominated Convergence Theorem, 
we get

\begin{equation}
\label{leb2}
\lim\limits_{p_1\to\infty}\varlimsup\limits_{p_2\to\infty}
\int\limits_{[t, T]^2}
\left(R_{p_1p_2}(t_1,t_2)\right)^{2n}dt_1 dt_2=0,\ \ \ \ \ 
\lim\limits_{p_1\to\infty}\varlimsup\limits_{p_2\to\infty}
\int\limits_t^T
\left(R_{p_1p_2}(t_1,t_1)\right)^{2n}dt_1=0.
\end{equation}

\vspace{5mm}

We will discuss the choice of integrable majorants
when applying Lebesgue's 
Dominated Convergence Theorem when we consider the case 
of arbitrary $k\in\mathbb{N}$ later in this section.

From 
(\ref{leto80010aaa}) and (\ref{leb2}) we obtain

\vspace{1mm}
$$
\lim\limits_{p_1\to \infty}
\varlimsup\limits_{p_2\to \infty}
{\sf M}\left\{\left|J[R_{p_1p_2}]_{T,t}^{(2)}\right|^{2n}
\right\}=0,\ \ \ n\in \mathbb{N}.
$$

\vspace{5mm}

Let us consider the case $k=3.$ Using (\ref{oop12}) (see below),
we have w. p. 1 

\vspace{4mm}
$$
J[R_{p_1p_2p_3}]_{T,t}^{(3)}=
\hbox{\vtop{\offinterlineskip\halign{
\hfil#\hfil\cr
{\rm l.i.m.}\cr
$\stackrel{}{{}_{N\to \infty}}$\cr
}} }\sum_{l_3=0}^{N-1}\sum_{l_2=0}^{N-1}
\sum_{l_1=0}^{N-1}
R_{p_1 p_2 p_3}(\tau_{l_1},\tau_{l_2},\tau_{l_3})
\Delta{\bf w}_{\tau_{l_1}}^{(i_1)}
\Delta{\bf w}_{\tau_{l_2}}^{(i_2)}
\Delta{\bf w}_{\tau_{l_3}}^{(i_3)}=
$$

\vspace{2mm}
$$
=\hbox{\vtop{\offinterlineskip\halign{
\hfil#\hfil\cr
{\rm l.i.m.}\cr
$\stackrel{}{{}_{N\to \infty}}$\cr
}} }\sum_{l_3=0}^{N-1}\sum_{l_2=0}^{l_3-1}
\sum_{l_1=0}^{l_2-1}\Biggl(
R_{p_1 p_2 p_3}(\tau_{l_1},\tau_{l_2},\tau_{l_3})
\Delta{\bf w}_{\tau_{l_1}}^{(i_1)}
\Delta{\bf w}_{\tau_{l_2}}^{(i_2)}
\Delta{\bf w}_{\tau_{l_3}}^{(i_3)}+\Biggr.
$$

\vspace{2mm}
$$
+
R_{p_1 p_2 p_3}(\tau_{l_1},\tau_{l_3},\tau_{l_2})
\Delta{\bf w}_{\tau_{l_1}}^{(i_1)}
\Delta{\bf w}_{\tau_{l_3}}^{(i_2)}
\Delta{\bf w}_{\tau_{l_2}}^{(i_3)}+
R_{p_1 p_2 p_3}(\tau_{l_2},\tau_{l_1},\tau_{l_3})
\Delta{\bf w}_{\tau_{l_2}}^{(i_1)}
\Delta{\bf w}_{\tau_{l_1}}^{(i_2)}
\Delta{\bf w}_{\tau_{l_3}}^{(i_3)}+
$$

\vspace{2mm}
$$
+
R_{p_1 p_2 p_3}(\tau_{l_2},\tau_{l_3},\tau_{l_1})
\Delta{\bf w}_{\tau_{l_2}}^{(i_1)}
\Delta{\bf w}_{\tau_{l_3}}^{(i_2)}
\Delta{\bf w}_{\tau_{l_1}}^{(i_3)}+
R_{p_1 p_2 p_3}(\tau_{l_3},\tau_{l_2},\tau_{l_1})
\Delta{\bf w}_{\tau_{l_3}}^{(i_1)}
\Delta{\bf w}_{\tau_{l_2}}^{(i_2)}
\Delta{\bf w}_{\tau_{l_1}}^{(i_3)}+
$$

\vspace{2mm}
$$
\Biggl.
+R_{p_1 p_2 p_3}
(\tau_{l_3},\tau_{l_1},\tau_{l_2})\Delta{\bf w}_{\tau_{l_3}}^{(i_1)}
\Delta{\bf w}_{\tau_{l_1}}^{(i_2)}
\Delta{\bf w}_{\tau_{l_2}}^{(i_3)}\Biggr)+
$$

\vspace{2mm}
$$
+
\hbox{\vtop{\offinterlineskip\halign{
\hfil#\hfil\cr
{\rm l.i.m.}\cr
$\stackrel{}{{}_{N\to \infty}}$\cr
}} }\sum_{l_3=0}^{N-1}\sum_{l_2=0}^{l_3-1}\Biggl(
R_{p_1 p_2 p_3}(\tau_{l_2},\tau_{l_2},\tau_{l_3})
\Delta{\bf w}_{\tau_{l_2}}^{(i_1)}
\Delta{\bf w}_{\tau_{l_2}}^{(i_2)}
\Delta{\bf w}_{\tau_{l_3}}^{(i_3)}+\Biggr.
$$

\vspace{2mm}
$$
+R_{p_1 p_2 p_3}(\tau_{l_2},\tau_{l_3},\tau_{l_2})
\Delta{\bf w}_{\tau_{l_2}}^{(i_1)}
\Delta{\bf w}_{\tau_{l_3}}^{(i_2)}
\Delta{\bf w}_{\tau_{l_2}}^{(i_3)}
\Biggl.+R_{p_1 p_2 p_3}(\tau_{l_3},\tau_{l_2},\tau_{l_2})
\Delta{\bf w}_{\tau_{l_3}}^{(i_1)}
\Delta{\bf w}_{\tau_{l_2}}^{(i_2)}
\Delta{\bf w}_{\tau_{l_2}}^{(i_3)}\Biggr)+
$$

\vspace{2mm}
$$
+
\hbox{\vtop{\offinterlineskip\halign{
\hfil#\hfil\cr
{\rm l.i.m.}\cr
$\stackrel{}{{}_{N\to \infty}}$\cr
}} }\sum_{l_3=0}^{N-1}\sum_{l_1=0}^{l_3-1}\Biggl(
R_{p_1 p_2 p_3}(\tau_{l_1},\tau_{l_3},\tau_{l_3})
\Delta{\bf w}_{\tau_{l_1}}^{(i_1)}
\Delta{\bf w}_{\tau_{l_3}}^{(i_2)}
\Delta{\bf w}_{\tau_{l_3}}^{(i_3)}+\Biggr.
$$

\vspace{2mm}
$$
+R_{p_1 p_2 p_3}(\tau_{l_3},\tau_{l_1},\tau_{l_3})
\Delta{\bf w}_{\tau_{l_3}}^{(i_1)}
\Delta{\bf w}_{\tau_{l_1}}^{(i_2)}
\Delta{\bf w}_{\tau_{l_3}}^{(i_3)}
\Biggl.+R_{p_1 p_2 p_3}(\tau_{l_3},\tau_{l_3},\tau_{l_1})
\Delta{\bf w}_{\tau_{l_3}}^{(i_1)}
\Delta{\bf w}_{\tau_{l_3}}^{(i_2)}
\Delta{\bf w}_{\tau_{l_1}}^{(i_3)}\Biggr)+
$$

\vspace{2mm}
$$
+\hbox{\vtop{\offinterlineskip\halign{
\hfil#\hfil\cr
{\rm l.i.m.}\cr
$\stackrel{}{{}_{N\to \infty}}$\cr
}} }\sum_{l_3=0}^{N-1}
R_{p_1 p_2 p_3}(\tau_{l_3},\tau_{l_3},\tau_{l_3})
\Delta{\bf w}_{\tau_{l_3}}^{(i_1)}
\Delta{\bf w}_{\tau_{l_3}}^{(i_2)}
\Delta{\bf w}_{\tau_{l_3}}^{(i_3)}=
$$

\vspace{2mm}
$$
=
\int\limits_t^T\int\limits_t^{t_3}\int\limits_t^{t_2}
R_{p_1 p_2 p_3}(t_1,t_2,t_3)
d{\bf w}_{t_1}^{(i_1)}
d{\bf w}_{t_2}^{(i_2)}
d{\bf w}_{t_3}^{(i_3)}+
\int\limits_t^T\int\limits_t^{t_3}\int\limits_t^{t_2}
R_{p_1 p_2 p_3}(t_1,t_3,t_2)
d{\bf w}_{t_1}^{(i_1)}
d{\bf w}_{t_2}^{(i_3)}
d{\bf w}_{t_3}^{(i_2)}+
$$

\vspace{1mm}
$$
+
\int\limits_t^T\int\limits_t^{t_3}\int\limits_t^{t_2}
R_{p_1 p_2 p_3}(t_2,t_1,t_3)
d{\bf w}_{t_1}^{(i_2)}
d{\bf w}_{t_2}^{(i_1)}
d{\bf w}_{t_3}^{(i_3)}+
\int\limits_t^T\int\limits_t^{t_3}\int\limits_t^{t_2}
R_{p_1 p_2 p_3}(t_2,t_3,t_1)
d{\bf w}_{t_1}^{(i_3)}
d{\bf w}_{t_2}^{(i_1)}
d{\bf w}_{t_3}^{(i_2)}+
$$

\vspace{1mm}
$$
+
\int\limits_t^T\int\limits_t^{t_3}\int\limits_t^{t_2}
R_{p_1 p_2 p_3}(t_3,t_2,t_1)
d{\bf w}_{t_1}^{(i_3)}
d{\bf w}_{t_2}^{(i_2)}
d{\bf w}_{t_3}^{(i_1)}+
\int\limits_t^T\int\limits_t^{t_3}\int\limits_t^{t_2}
R_{p_1 p_2 p_3}(t_3,t_1,t_2)
d{\bf w}_{t_1}^{(i_2)}
d{\bf w}_{t_2}^{(i_3)}
d{\bf w}_{t_3}^{(i_1)}+
$$

\vspace{1mm}
$$
+{\bf 1}_{\{i_1=i_2\ne 0\}}
\int\limits_t^T\int\limits_t^{t_3}
R_{p_1 p_2 p_3}(t_2,t_2,t_3)
dt_2
d{\bf w}_{t_3}^{(i_3)}+{\bf 1}_{\{i_1=i_3\ne 0\}}
\int\limits_t^T\int\limits_t^{t_3}
R_{p_1 p_2 p_3}(t_2,t_3,t_2)
dt_2
d{\bf w}_{t_3}^{(i_2)}+
$$

\vspace{1mm}
$$
+{\bf 1}_{\{i_2=i_3\ne 0\}}
\int\limits_t^T\int\limits_t^{t_3}
R_{p_1 p_2 p_3}(t_3,t_2,t_2)
dt_2
d{\bf w}_{t_3}^{(i_1)}+{\bf 1}_{\{i_2=i_3\ne 0\}}
\int\limits_t^T\int\limits_t^{t_3}
R_{p_1 p_2 p_3}(t_1,t_3,t_3)
d{\bf w}_{t_1}^{(i_1)}dt_3+
$$

\vspace{1mm}
\begin{equation}
\label{s1s}
+{\bf 1}_{\{i_1=i_3\ne 0\}}
\int\limits_t^T\int\limits_t^{t_3}
R_{p_1 p_2 p_3}(t_3,t_1,t_3)
d{\bf w}_{t_1}^{(i_2)}dt_3
+{\bf 1}_{\{i_1=i_2\ne 0\}}
\int\limits_t^T\int\limits_t^{t_3}
R_{p_1 p_2 p_3}(t_3,t_3,t_1)
d{\bf w}_{t_1}^{(i_3)}dt_3.
\end{equation}

\vspace{7mm}

Applying Lemma 4, we obtain

\vspace{3mm}
$$
{\sf M}\left\{\left|J[R_{p_1p_2p_3}]_{T,t}^{(3)}\right|^{2n}
\right\}\le 
C_n\Biggl(\int\limits_t^T\int\limits_t^{t_3}\int\limits_t^{t_2}
\Biggl(
\left(R_{p_1 p_2 p_3}(t_1,t_2,t_3)\right)^{2n}+
\left(R_{p_1 p_2 p_3}(t_1,t_3,t_2)\right)^{2n}+\Biggr.\Biggr.
$$

\vspace{1mm}
$$
+\left(R_{p_1 p_2 p_3}(t_2,t_1,t_3)\right)^{2n}+
\left(R_{p_1 p_2 p_3}(t_2,t_3,t_1)\right)^{2n}+
\left(R_{p_1 p_2 p_3}(t_3,t_2,t_1)\right)^{2n}+
$$

\vspace{1mm}
$$
\Biggl.
+\left(R_{p_1 p_2 p_3}(t_3,t_1,t_2)\right)^{2n}\Biggr)dt_1dt_2dt_3+
$$

\vspace{1mm}
$$
+
\int\limits_t^T\int\limits_t^{t_3}\Biggl(
{\bf 1}_{\{i_1=i_2\ne 0\}}\Biggl(
\left(R_{p_1 p_2 p_3}(t_2,t_2,t_3)\right)^{2n}+
\left(R_{p_1 p_2 p_3}(t_3,t_3,t_2)\right)^{2n}\Biggr)+\Biggr.
$$

\vspace{1mm}
$$
+{\bf 1}_{\{i_1=i_3\ne 0\}}\Biggl(
\left(R_{p_1 p_2 p_3}(t_2,t_3,t_2)\right)^{2n}+
\left(R_{p_1 p_2 p_3}(t_3,t_2,t_3)\right)^{2n}\Biggr)+
$$

\vspace{1mm}
\begin{equation}
\label{oop16}
\Biggl.
+{\bf 1}_{\{i_2=i_3\ne 0\}}\Biggl(
\left(R_{p_1 p_2 p_3}(t_3,t_2,t_2)\right)^{2n}+
\left(R_{p_1 p_2 p_3}(t_2,t_3,t_3)\right)^{2n}\Biggr)dt_2dt_3\Biggr),\ \ \
C_n<\infty.
\end{equation}

\vspace{6mm}

Further, we have

\vspace{2mm}
$$
\int\limits_t^T\int\limits_t^{t_3}\int\limits_t^{t_2}
\Biggl(
\left(R_{p_1 p_2 p_3}(t_1,t_2,t_3)\right)^{2n}+
\left(R_{p_1 p_2 p_3}(t_1,t_3,t_2)\right)^{2n}+
\left(R_{p_1 p_2 p_3}(t_2,t_1,t_3)\right)^{2n}+\Biggr.
$$

\vspace{1mm}
$$
+
\left(R_{p_1 p_2 p_3}(t_2,t_3,t_1)\right)^{2n}+
\left(R_{p_1 p_2 p_3}(t_3,t_2,t_1)\right)^{2n}
+\left(R_{p_1 p_2 p_3}(t_3,t_1,t_2)\right)^{2n}\Biggr)dt_1dt_2dt_3=
$$

\vspace{1mm}
\begin{equation}
\label{zero1}
=
\int\limits_{[t, T]^3}
\left(R_{p_1 p_2 p_3}(t_1,t_2,t_3)\right)^{2n}dt_1dt_2dt_3,
\end{equation}

\vspace{5mm}

$$
\int\limits_t^T\int\limits_t^{t_3}\biggl(
\left(R_{p_1 p_2 p_3}(t_2,t_2,t_3)\right)^{2n}+
\left(R_{p_1 p_2 p_3}(t_3,t_3,t_2)\right)^{2n}\biggr)dt_2dt_3=
$$

\vspace{1mm}
$$
=\int\limits_t^T\int\limits_t^{t_3}
\left(R_{p_1 p_2 p_3}(t_2,t_2,t_3)\right)^{2n}dt_2dt_3+
\int\limits_t^T\int\limits_{t_3}^{T}
\left(R_{p_1 p_2 p_3}(t_2,t_2,t_3)\right)^{2n}dt_2dt_3=
$$

\vspace{1mm}
\begin{equation}
\label{zero2}
=
\int\limits_{[t, T]^2}
\left(R_{p_1 p_2 p_3}(t_2,t_2,t_3)\right)^{2n}dt_2dt_3,
\end{equation}

\vspace{5mm}

$$
\int\limits_t^T\int\limits_t^{t_3}\biggl(
\left(R_{p_1 p_2 p_3}(t_2,t_3,t_2)\right)^{2n}+
\left(R_{p_1 p_2 p_3}(t_3,t_2,t_3)\right)^{2n}\biggr)dt_2dt_3=
$$

\vspace{1mm}
$$
=\int\limits_t^T\int\limits_t^{t_3}
\left(R_{p_1 p_2 p_3}(t_2,t_3,t_2)\right)^{2n}dt_2dt_3+
\int\limits_t^T\int\limits_{t_3}^{T}
\left(R_{p_1 p_2 p_3}(t_2,t_3,t_2)\right)^{2n}dt_2dt_3=
$$

\vspace{1mm}
\begin{equation}
\label{zero3}
=
\int\limits_{[t, T]^2}
\left(R_{p_1 p_2 p_3}(t_2,t_3,t_2)\right)^{2n}dt_2dt_3,
\end{equation}

\vspace{5mm}

$$
\int\limits_t^T\int\limits_t^{t_3}\biggl(
\left(R_{p_1 p_2 p_3}(t_3,t_2,t_2)\right)^{2n}+
\left(R_{p_1 p_2 p_3}(t_2,t_3,t_3)\right)^{2n}\biggr)dt_2dt_3=
$$

\vspace{1mm}
$$
=\int\limits_t^T\int\limits_t^{t_3}
\left(R_{p_1 p_2 p_3}(t_3,t_2,t_2)\right)^{2n}dt_2dt_3+
\int\limits_t^T\int\limits_{t_3}^{T}
\left(R_{p_1 p_2 p_3}(t_3,t_2,t_2)\right)^{2n}dt_2dt_3=
$$

\vspace{1mm}
\begin{equation}
\label{zero4}
=
\int\limits_{[t, T]^2}
\left(R_{p_1 p_2 p_3}(t_3,t_2,t_2)\right)^{2n}dt_2dt_3.
\end{equation}

\vspace{6mm}

Combining (\ref{oop16}) and (\ref{zero1})--(\ref{zero4}),
we get

\vspace{1mm}
$$
{\sf M}\left\{\left|J[R_{p_1p_2p_3}]_{T,t}^{(3)}\right|^{2n}
\right\}\le 
C_n\left(
\int\limits_{[t, T]^3}
\left(R_{p_1 p_2 p_3}(t_1,t_2,t_3)\right)^{2n}dt_1dt_2dt_3+\right.
$$

\vspace{2mm}
$$
+
{\bf 1}_{\{i_1=i_2\ne 0\}}
\int\limits_{[t, T]^2}
\left(R_{p_1 p_2 p_3}(t_2,t_2,t_3)\right)^{2n}dt_2dt_3+
$$

\vspace{1mm}
$$
+{\bf 1}_{\{i_1=i_3\ne 0\}}
\int\limits_{[t, T]^2}
\left(R_{p_1 p_2 p_3}(t_2,t_3,t_2)\right)^{2n}dt_2dt_3+
$$

\vspace{-1mm}
\begin{equation}
\label{oop16xxx}
\left.
+{\bf 1}_{\{i_2=i_3\ne 0\}}
\int\limits_{[t, T]^2}
\left(R_{p_1 p_2 p_3}(t_3,t_2,t_2)\right)^{2n}dt_2dt_3\right),\ \ \ 
C_n<\infty.
\end{equation}

\vspace{5mm}

Since the integrals on the right-hand side of (\ref{oop16xxx}) 
exist as Riemann integrals, then they are equal to the 
corresponding Lebesgue integrals. 
Moreover, 

\vspace{1mm}
$$
\lim\limits_{p_1\to\infty}\lim\limits_{p_2\to\infty}\lim\limits_{p_3\to\infty}
R_{p_1 p_2 p_3}(t_1,t_2,t_3)=0\ \ \ \hbox{when}\ \ \ (t_1,t_2,t_3)\in (t, T)^3,
$$

\vspace{4mm}
\noindent
where the left-hand side 
is bounded on $[t, T]^3.$

According to the proof of Lemma 1 and (\ref{30.46}) for $k=3$, we have

\vspace{2mm}
$$
R_{p_1p_2p_3}(t_1,t_2,t_3)=\left(K^{*}(t_1,t_2,t_3)-\sum\limits_{j_1=0}^{p_1}
C_{j_1}(t_2,t_3)\phi_{j_1}(t_1)\right)+
$$

\vspace{2mm}
$$
+\left(
\sum\limits_{j_1=0}^{p_1}\left(C_{j_1}(t_2,t_3)-
\sum\limits_{j_2=0}^{p_2}
C_{j_2j_1}(t_3)\phi_{j_2}(t_2)\right)
\phi_{j_1}(t_1)\right)+
$$

\vspace{2mm}
\begin{equation}
\label{lab11}
+\left(
\sum\limits_{j_1=0}^{p_1}\sum\limits_{j_2=0}^{p_2}\left(C_{j_2j_1}(t_3)-
\sum\limits_{j_3=0}^{p_3}
C_{j_3j_2j_1}\phi_{j_3}(t_3)\right)
\phi_{j_2}(t_2)\phi_{j_1}(t_1)\right),
\end{equation}

\vspace{5mm}
\noindent
where

\vspace{-2mm}
$$
C_{j_1}(t_2,t_3)=\int\limits_t^T
K^{*}(t_1,t_2,t_3)\phi_{j_1}(t_1)dt_1,\ \ \ 
C_{j_2j_1}(t_3)=\int\limits_{[t, T]^2}
K^{*}(t_1,t_2,t_3)\phi_{j_1}(t_1)\phi_{j_2}(t_2)dt_1 dt_2.
$$

\vspace{4mm}

Then, applying three times (we mean an iterated passage to the limit
$\lim\limits_{p_1\to\infty}\varlimsup\limits_{p_2\to\infty}
\varlimsup\limits_{p_3\to\infty}$)
the Lebesgue's Dominated Convergence Theorem,
we obtain

\begin{equation}
\label{leb21}
\lim\limits_{p_1\to\infty}\varlimsup\limits_{p_2\to\infty}
\varlimsup\limits_{p_3\to\infty}
\int\limits_{[t, T]^3}
\left(R_{p_1p_2p_3}(t_1,t_2,t_3)\right)^{2n}dt_1 dt_2 dt_3=0,
\end{equation}

\begin{equation}
\label{leb22}
\lim\limits_{p_1\to\infty}\varlimsup\limits_{p_2\to\infty}
\varlimsup\limits_{p_3\to\infty}
\int\limits_{[t, T]^2}
\left(R_{p_1 p_2 p_3}(t_2,t_2,t_3)\right)^{2n}dt_2dt_3=0,
\end{equation}

\begin{equation}
\label{leb23}
\lim\limits_{p_1\to\infty}\varlimsup\limits_{p_2\to\infty}
\varlimsup\limits_{p_3\to\infty}
\int\limits_{[t, T]^2}
\left(R_{p_1 p_2 p_3}(t_2,t_3,t_2)\right)^{2n}dt_2dt_3=0,
\end{equation}

\begin{equation}
\label{leb24}
\lim\limits_{p_1\to\infty}\varlimsup\limits_{p_2\to\infty}
\varlimsup\limits_{p_3\to\infty}
\int\limits_{[t, T]^2}
\left(R_{p_1 p_2 p_3}(t_3,t_2,t_2)\right)^{2n}dt_2dt_3=0.
\end{equation}

\vspace{5mm}

From 
(\ref{oop16xxx})--(\ref{leb24}) we get

\vspace{1mm}
$$
\lim\limits_{p_1\to\infty}
\varlimsup\limits_{p_2\to\infty}\varlimsup\limits_{p_3\to\infty}
{\sf M}\left\{\left|J[R_{p_1p_2p_3}]_{T,t}^{(3)}\right|^{2n}
\right\}=0,\ \ \ n\in \mathbb{N}.
$$

\vspace{5mm}

Let us consider the case $k=4.$ Using (\ref{huh}) (see below),
we have w. p. 1 

\vspace{2mm}
$$
J[R_{p_1p_2p_3p_4}]_{T,t}^{(4)}=
$$

\vspace{2mm}

$$
=
\hbox{\vtop{\offinterlineskip\halign{
\hfil#\hfil\cr
{\rm l.i.m.}\cr
$\stackrel{}{{}_{N\to \infty}}$\cr
}} }\sum_{l_4=0}^{N-1}\sum_{l_3=0}^{N-1}\sum_{l_2=0}^{N-1}
\sum_{l_1=0}^{N-1}
R_{p_1 p_2 p_3 p_4}(\tau_{l_1},\tau_{l_2},\tau_{l_3}, \tau_{l_4})
\Delta{\bf w}_{\tau_{l_1}}^{(i_1)}
\Delta{\bf w}_{\tau_{l_2}}^{(i_2)}
\Delta{\bf w}_{\tau_{l_3}}^{(i_3)}\Delta{\bf w}_{\tau_{l_4}}^{(i_4)}=
$$

\vspace{1.5mm}
$$
=\hbox{\vtop{\offinterlineskip\halign{
\hfil#\hfil\cr
{\rm l.i.m.}\cr
$\stackrel{}{{}_{N\to \infty}}$\cr
}} }\sum_{l_4=0}^{N-1}\sum_{l_3=0}^{l_4-1}\sum_{l_2=0}^{l_3-1}
\sum_{l_1=0}^{l_2-1}
\sum\limits_{(l_1,l_2,l_3,l_4)}
\biggl(R_{p_1 p_2 p_3 p_4}(\tau_{l_1},\tau_{l_2},\tau_{l_3}, \tau_{l_4})
\Delta{\bf w}_{\tau_{l_1}}^{(i_1)}
\Delta{\bf w}_{\tau_{l_2}}^{(i_2)}
\Delta{\bf w}_{\tau_{l_3}}^{(i_3)}\Delta{\bf w}_{\tau_{l_4}}^{(i_4)}\biggr)+
$$

\vspace{1.5mm}
$$
+\hbox{\vtop{\offinterlineskip\halign{
\hfil#\hfil\cr
{\rm l.i.m.}\cr
$\stackrel{}{{}_{N\to \infty}}$\cr
}} }\sum_{l_4=0}^{N-1}\sum_{l_3=0}^{l_4-1}\sum_{l_2=0}^{l_3-1}
\sum\limits_{(l_2,l_2,l_3,l_4)}
\biggl(R_{p_1 p_2 p_3 p_4}(\tau_{l_2},\tau_{l_2},\tau_{l_3}, \tau_{l_4})
\Delta{\bf w}_{\tau_{l_2}}^{(i_1)}
\Delta{\bf w}_{\tau_{l_2}}^{(i_2)}
\Delta{\bf w}_{\tau_{l_3}}^{(i_3)}\Delta{\bf w}_{\tau_{l_4}}^{(i_4)}\biggr)+
$$

\vspace{1.5mm}
$$
+\hbox{\vtop{\offinterlineskip\halign{
\hfil#\hfil\cr
{\rm l.i.m.}\cr
$\stackrel{}{{}_{N\to \infty}}$\cr
}} }\sum_{l_4=0}^{N-1}\sum_{l_3=0}^{l_4-1}\sum_{l_1=0}^{l_3-1}
\sum\limits_{(l_1,l_3,l_3,l_4)}
\biggl(R_{p_1 p_2 p_3 p_4}(\tau_{l_1},\tau_{l_3},\tau_{l_3}, \tau_{l_4})
\Delta{\bf w}_{\tau_{l_1}}^{(i_1)}
\Delta{\bf w}_{\tau_{l_3}}^{(i_2)}
\Delta{\bf w}_{\tau_{l_3}}^{(i_3)}\Delta{\bf w}_{\tau_{l_4}}^{(i_4)}\biggr)+
$$

\vspace{1.5mm}
$$
+\hbox{\vtop{\offinterlineskip\halign{
\hfil#\hfil\cr
{\rm l.i.m.}\cr
$\stackrel{}{{}_{N\to \infty}}$\cr
}} }\sum_{l_4=0}^{N-1}\sum_{l_2=0}^{l_4-1}\sum_{l_1=0}^{l_2-1}
\sum\limits_{(l_1,l_2,l_4,l_4)}
\biggl(R_{p_1 p_2 p_3 p_4}(\tau_{l_1},\tau_{l_2},\tau_{l_4}, \tau_{l_4})
\Delta{\bf w}_{\tau_{l_1}}^{(i_1)}
\Delta{\bf w}_{\tau_{l_2}}^{(i_2)}
\Delta{\bf w}_{\tau_{l_4}}^{(i_3)}\Delta{\bf w}_{\tau_{l_4}}^{(i_4)}\biggr)+
$$

\vspace{1.5mm}
$$
+\hbox{\vtop{\offinterlineskip\halign{
\hfil#\hfil\cr
{\rm l.i.m.}\cr
$\stackrel{}{{}_{N\to \infty}}$\cr
}} }\sum_{l_4=0}^{N-1}\sum_{l_3=0}^{l_4-1}
\sum\limits_{(l_3,l_3,l_3,l_4)}
\biggl(R_{p_1 p_2 p_3 p_4}(\tau_{l_3},\tau_{l_3},\tau_{l_3}, \tau_{l_4})
\Delta{\bf w}_{\tau_{l_3}}^{(i_1)}
\Delta{\bf w}_{\tau_{l_3}}^{(i_2)}
\Delta{\bf w}_{\tau_{l_3}}^{(i_3)}\Delta{\bf w}_{\tau_{l_4}}^{(i_4)}\biggr)+
$$

\vspace{1.5mm}
$$
+\hbox{\vtop{\offinterlineskip\halign{
\hfil#\hfil\cr
{\rm l.i.m.}\cr
$\stackrel{}{{}_{N\to \infty}}$\cr
}} }\sum_{l_4=0}^{N-1}\sum_{l_2=0}^{l_4-1}
\sum\limits_{(l_2,l_2,l_4,l_4)}
\biggl(R_{p_1 p_2 p_3 p_4}(\tau_{l_2},\tau_{l_2},\tau_{l_4}, \tau_{l_4})
\Delta{\bf w}_{\tau_{l_2}}^{(i_1)}
\Delta{\bf w}_{\tau_{l_2}}^{(i_2)}
\Delta{\bf w}_{\tau_{l_4}}^{(i_3)}\Delta{\bf w}_{\tau_{l_4}}^{(i_4)}\biggr)+
$$

\vspace{1.5mm}
$$
+\hbox{\vtop{\offinterlineskip\halign{
\hfil#\hfil\cr
{\rm l.i.m.}\cr
$\stackrel{}{{}_{N\to \infty}}$\cr
}} }\sum_{l_4=0}^{N-1}\sum_{l_1=0}^{l_4-1}
\sum\limits_{(l_1,l_4,l_4,l_4)}
\biggl(R_{p_1 p_2 p_3 p_4}(\tau_{l_1},\tau_{l_4},\tau_{l_4}, \tau_{l_4})
\Delta{\bf w}_{\tau_{l_1}}^{(i_1)}
\Delta{\bf w}_{\tau_{l_4}}^{(i_2)}
\Delta{\bf w}_{\tau_{l_4}}^{(i_3)}\Delta{\bf w}_{\tau_{l_4}}^{(i_4)}\biggr)+
$$

\vspace{1.5mm}
$$
+\hbox{\vtop{\offinterlineskip\halign{
\hfil#\hfil\cr
{\rm l.i.m.}\cr
$\stackrel{}{{}_{N\to \infty}}$\cr
}} }\sum_{l_4=0}^{N-1}
R_{p_1 p_2 p_3 p_4}(\tau_{l_4},\tau_{l_4},\tau_{l_4}, \tau_{l_4})
\Delta{\bf w}_{\tau_{l_4}}^{(i_1)}
\Delta{\bf w}_{\tau_{l_4}}^{(i_2)}
\Delta{\bf w}_{\tau_{l_4}}^{(i_3)}\Delta{\bf w}_{\tau_{l_4}}^{(i_4)}=
$$

\vspace{1.5mm}

$$
=
\int\limits_t^T\int\limits_t^{t_4}\int\limits_t^{t_3}\int\limits_t^{t_2}
\sum\limits_{(t_1,t_2,t_3,t_4)}\biggl(R_{p_1 p_2 p_3 p_4}(t_1,t_2,t_3,t_4)
d{\bf w}_{t_1}^{(i_1)}
d{\bf w}_{t_2}^{(i_2)}
d{\bf w}_{t_3}^{(i_3)}
d{\bf w}_{t_4}^{(i_4)}\biggr)+
$$

\vspace{2mm}

$$
+
{\bf 1}_{\{i_1=i_2\ne 0\}}
\int\limits_t^T\int\limits_t^{t_4}\int\limits_t^{t_3}
\sum\limits_{(t_1,t_3,t_4)}\biggl(R_{p_1 p_2 p_3 p_4}(t_1,t_1,t_3,t_4)
dt_1
d{\bf w}_{t_3}^{(i_3)}
d{\bf w}_{t_4}^{(i_4)}\biggr)+
$$

\vspace{2mm}

$$
+
{\bf 1}_{\{i_1=i_3\ne 0\}}
\int\limits_t^T\int\limits_t^{t_4}\int\limits_t^{t_2}
\sum\limits_{(t_1,t_2,t_4)}\biggl(
R_{p_1 p_2 p_3 p_4}(t_1,t_2,t_1,t_4)
dt_1
d{\bf w}_{t_2}^{(i_2)}
d{\bf w}_{t_4}^{(i_4)}\biggr)+
$$

\vspace{2mm}

$$
+
{\bf 1}_{\{i_1=i_4\ne 0\}}
\int\limits_t^T\int\limits_t^{t_3}\int\limits_t^{t_2}
\sum\limits_{(t_1,t_2,t_3)}
\biggl(R_{p_1 p_2 p_3 p_4}(t_1,t_2,t_3,t_1)
dt_1 d{\bf w}_{t_2}^{(i_2)}
d{\bf w}_{t_3}^{(i_3)}\biggr)+
$$

\vspace{2mm}

$$
+
{\bf 1}_{\{i_2=i_3\ne 0\}}
\int\limits_t^T\int\limits_t^{t_4}\int\limits_t^{t_2}
\sum\limits_{(t_1,t_2,t_4)}
\biggl(R_{p_1 p_2 p_3 p_4}(t_1,t_2,t_2,t_4)
d{\bf w}_{t_1}^{(i_1)}
dt_2
d{\bf w}_{t_4}^{(i_4)}\biggr)+
$$

\vspace{2mm}

$$
+
{\bf 1}_{\{i_2=i_4\ne 0\}}
\int\limits_t^T\int\limits_t^{t_3}\int\limits_t^{t_2}
\sum\limits_{(t_1,t_2,t_3)}
\biggl(R_{p_1 p_2 p_3 p_4}(t_1,t_2,t_3,t_2)
d{\bf w}_{t_1}^{(i_1)}
dt_2
d{\bf w}_{t_3}^{(i_3)}\biggr)+
$$

\vspace{2mm}

$$
+
{\bf 1}_{\{i_3=i_4\ne 0\}}
\int\limits_t^T\int\limits_t^{t_3}\int\limits_t^{t_2}
\sum\limits_{(t_1,t_2,t_3)}
\biggl(R_{p_1 p_2 p_3 p_4}(t_1,t_2,t_3,t_3)
d{\bf w}_{t_1}^{(i_1)}
d{\bf w}_{t_2}^{(i_2)}dt_3\biggr)+
$$

\vspace{2mm}

$$
+
{\bf 1}_{\{i_1=i_2\ne 0\}}{\bf 1}_{\{i_3=i_4\ne 0\}}
\left(\int\limits_t^T\int\limits_t^{t_4}
R_{p_1 p_2 p_3 p_4}(t_2,t_2,t_4,t_4)dt_2 dt_4
+\right.
$$

\vspace{1mm}
$$
\left.+ \int\limits_t^T\int\limits_t^{t_4}
R_{p_1 p_2 p_3 p_4}(t_4,t_4,t_2,t_2)dt_2 dt_4\right)+
$$

\vspace{2mm}

$$
+
{\bf 1}_{\{i_1=i_3\ne 0\}}{\bf 1}_{\{i_2=i_4\ne 0\}}
\left(\int\limits_t^T\int\limits_t^{t_4}
R_{p_1 p_2 p_3 p_4}(t_2,t_4,t_2,t_4)dt_2 dt_4
+\right.
$$

\vspace{1mm}
$$
\left.+\int\limits_t^T\int\limits_t^{t_4}
R_{p_1 p_2 p_3 p_4}(t_4,t_2,t_4,t_2)dt_2 dt_4\right)+
$$

\vspace{2mm}

$$
+
{\bf 1}_{\{i_1=i_4\ne 0\}}{\bf 1}_{\{i_2=i_3\ne 0\}}
\left(\int\limits_t^T\int\limits_t^{t_4}
R_{p_1 p_2 p_3 p_4}(t_2,t_4,t_4,t_2)dt_2 dt_4
+\right.
$$

\vspace{1mm}
\begin{equation}
\label{sogl}
\left.+\int\limits_t^T\int\limits_t^{t_4}
R_{p_1 p_2 p_3 p_4}(t_4,t_2,t_2,t_4)dt_2 dt_4\right),
\end{equation}

\vspace{5mm}
\noindent
where
the expression
$$
\sum\limits_{(a_1, \ldots, a_k)}
$$

\vspace{3mm}
\noindent
means the sum with respect to
all possible  
permutations $(a_1, \ldots, a_k)$. 
Note that an analogue of (\ref{sogl}) was 
obtained in \cite{17}, Sect.~6 (also see \cite{10a}-\cite{12aa-afterxxx})
with using a different approach.

By analogy with
(\ref{oop16xxx}) we obtain

\vspace{1mm}
$$
{\sf M}\left\{\left|J[R_{p_1p_2p_3p_4}]_{T,t}^{(4)}\right|^{2n}
\right\}\le 
C_n\left(
\int\limits_{[t, T]^4}
\left(R_{p_1 p_2 p_3 p_4}(t_1,t_2,t_3,t_4)\right)^{2n}dt_1dt_2dt_3dt_4+\right.
$$

\vspace{1.5mm}
$$
+
{\bf 1}_{\{i_1=i_2\ne 0\}}
\int\limits_{[t, T]^3}
\left(R_{p_1 p_2 p_3 p_4}(t_2,t_2,t_3,t_4)\right)^{2n}dt_2dt_3dt_4+
$$

\vspace{1.5mm}
$$
+
{\bf 1}_{\{i_1=i_3\ne 0\}}
\int\limits_{[t, T]^3}
\left(R_{p_1 p_2 p_3 p_4}(t_2,t_3,t_2,t_4)\right)^{2n}dt_2dt_3dt_4+
$$

\vspace{1.5mm}
$$
+
{\bf 1}_{\{i_1=i_4\ne 0\}}
\int\limits_{[t, T]^3}
\left(R_{p_1 p_2 p_3 p_4}(t_2,t_3,t_4,t_2)\right)^{2n}dt_2dt_3dt_4+
$$

\vspace{1.5mm}
$$
+
{\bf 1}_{\{i_2=i_3\ne 0\}}
\int\limits_{[t, T]^3}
\left(R_{p_1 p_2 p_3 p_4}(t_3,t_2,t_2,t_4)\right)^{2n}dt_2dt_3dt_4+
$$

\vspace{1.5mm}
$$
+
{\bf 1}_{\{i_2=i_4\ne 0\}}
\int\limits_{[t, T]^3}
\left(R_{p_1 p_2 p_3 p_4}(t_3,t_2,t_4,t_2)\right)^{2n}dt_2dt_3dt_4+
$$

\vspace{1.5mm}
$$
+
{\bf 1}_{\{i_3=i_4\ne 0\}}
\int\limits_{[t, T]^3}
\left(R_{p_1 p_2 p_3 p_4}(t_3,t_4,t_2,t_2)\right)^{2n}dt_2dt_3dt_4+
$$

\vspace{1.5mm}
$$
+
{\bf 1}_{\{i_1=i_2\ne 0\}}{\bf 1}_{\{i_3=i_4\ne 0\}}
\int\limits_{[t, T]^2}
\left(R_{p_1 p_2 p_3 p_4}(t_2,t_2,t_4,t_4)\right)^{2n}dt_2dt_4+
$$

\vspace{1.5mm}
$$
+
{\bf 1}_{\{i_1=i_3\ne 0\}}{\bf 1}_{\{i_2=i_4\ne 0\}}
\int\limits_{[t, T]^2}
\left(R_{p_1 p_2 p_3 p_4}(t_2,t_4,t_2,t_4)\right)^{2n}dt_2dt_4+
$$

\vspace{1.5mm}
\begin{equation}
\label{udar}
\left.+
{\bf 1}_{\{i_1=i_4\ne 0\}}{\bf 1}_{\{i_2=i_3\ne 0\}}
\int\limits_{[t, T]^2}
\left(R_{p_1 p_2 p_3 p_4}(t_2,t_4,t_4,t_2)\right)^{2n}dt_2dt_4\right),\ \ \
C_n<\infty.
\end{equation}

\vspace{6mm}

Since the integrals on the right-hand side of (\ref{udar}) 
exist as Riemann integrals, then they are equal to the 
corresponding Lebesgue integrals. 
Moreover, 

\vspace{1mm}
$$
\lim\limits_{p_1\to\infty}\lim\limits_{p_2\to\infty}\lim\limits_{p_3\to\infty}
\lim\limits_{p_4\to\infty}
R_{p_1 p_2 p_3 p_4}(t_1,t_2,t_3,t_4)=0\ \ \ \hbox{when}\ \ \ (t_1,t_2,t_3,t_4)\in (t, T)^4,
$$

\vspace{5mm}
\noindent
where the left-hand side 
is bounded on $[t, T]^4.$

According to the proof of Lemma 1 and (\ref{30.46}) for $k=4$, we have

\vspace{2mm}
$$
R_{p_1p_2p_3p_4}
(t_1,t_2,t_3,t_4)=\left(K^{*}(t_1,t_2,t_3,t_4)-
\sum\limits_{j_1=0}^{p_1}
C_{j_1}(t_2,t_3,t_4)\phi_{j_1}(t_1)\right)+
$$

\vspace{2mm}
$$
+\left(
\sum\limits_{j_1=0}^{p_1}\left(C_{j_1}(t_2,t_3,t_4)-
\sum\limits_{j_2=0}^{p_2}
C_{j_2j_1}(t_3,t_4)\phi_{j_2}(t_2)\right)
\phi_{j_1}(t_1)\right)+
$$

\vspace{2mm}
$$
+\left(
\sum\limits_{j_1=0}^{p_1}\sum\limits_{j_2=0}^{p_2}\left(C_{j_2j_1}(t_3,t_4)-
\sum\limits_{j_3=0}^{p_3}
C_{j_3j_2j_1}(t_4)\phi_{j_3}(t_3)\right)
\phi_{j_2}(t_2)\phi_{j_1}(t_1)\right)+
$$

\vspace{2mm}
$$
+\left(
\sum\limits_{j_1=0}^{p_1}\sum\limits_{j_2=0}^{p_2}\sum\limits_{j_3=0}^{p_3}
\left(C_{j_3j_2j_1}(t_4)-
\sum\limits_{j_4=0}^{p_4}
C_{j_4j_3j_2j_1}\phi_{j_4}(t_4)\right)
\phi_{j_3}(t_3)\phi_{j_2}(t_2)\phi_{j_1}(t_1)\right),
$$

\vspace{5mm}
\noindent
where
$$
C_{j_1}(t_2,t_3,t_4)=\int\limits_t^T
K^{*}(t_1,t_2,t_3,t_4)\phi_{j_1}(t_1)dt_1,
$$

\vspace{1.5mm}
$$
C_{j_2j_1}(t_3,t_4)=\int\limits_{[t, T]^2}
K^{*}(t_1,t_2,t_3,t_4)\phi_{j_1}(t_1)\phi_{j_2}(t_2)dt_1 dt_2,
$$

\vspace{1.5mm}
$$
C_{j_3j_2j_1}(t_4)=\int\limits_{[t, T]^3}
K^{*}(t_1,t_2,t_3,t_4)\phi_{j_1}(t_1)\phi_{j_2}(t_2)\phi_{j_3}(t_3)
dt_1 dt_2 dt_3.
$$

\vspace{5mm}

Then, applying  four times (we mean an iterated passage to the limit
$\lim\limits_{p_1\to\infty}\varlimsup\limits_{p_2\to\infty}
\varlimsup\limits_{p_3\to\infty}\varlimsup\limits_{p_4\to\infty}$)
the Lebesgue's Dominated Convergence Theorem,
we obtain

\vspace{1mm}

\begin{equation}
\label{final1}
\lim\limits_{p_1\to\infty}\varlimsup\limits_{p_2\to\infty}
\varlimsup\limits_{p_3\to\infty}\varlimsup\limits_{p_4\to\infty}\int\limits_{[t, T]^4}
\left(R_{p_1 p_2 p_3 p_4}(t_1,t_2,t_3,t_4)\right)^{2n}dt_1dt_2dt_3dt_4=0,
\end{equation}

\vspace{1.5mm}

\begin{equation}
\label{final2}
\lim\limits_{p_1\to\infty}\varlimsup\limits_{p_2\to\infty}
\varlimsup\limits_{p_3\to\infty}\varlimsup\limits_{p_4\to\infty}\int\limits_{[t, T]^3}
\left(R_{p_1 p_2 p_3 p_4}(t_2,t_2,t_3,t_4)\right)^{2n}dt_2dt_3dt_4=0,
\end{equation}

\vspace{1.5mm}

\begin{equation}
\label{final3}
\lim\limits_{p_1\to\infty}\varlimsup\limits_{p_2\to\infty}
\varlimsup\limits_{p_3\to\infty}\varlimsup\limits_{p_4\to\infty}\int\limits_{[t, T]^3}
\left(R_{p_1 p_2 p_3 p_4}(t_2,t_3,t_2,t_4)\right)^{2n}dt_2dt_3dt_4=0,
\end{equation}

\vspace{1.5mm}

\begin{equation}
\label{final4}
\lim\limits_{p_1\to\infty}\varlimsup\limits_{p_2\to\infty}
\varlimsup\limits_{p_3\to\infty}\varlimsup\limits_{p_4\to\infty}\int\limits_{[t, T]^3}
\left(R_{p_1 p_2 p_3 p_4}(t_2,t_3,t_4,t_2)\right)^{2n}dt_2dt_3dt_4=0,
\end{equation}

\vspace{1.5mm}

\begin{equation}
\label{final5}
\lim\limits_{p_1\to\infty}\varlimsup\limits_{p_2\to\infty}
\varlimsup\limits_{p_3\to\infty}\varlimsup\limits_{p_4\to\infty}\int\limits_{[t, T]^3}
\left(R_{p_1 p_2 p_3 p_4}(t_3,t_2,t_2,t_4)\right)^{2n}dt_2dt_3dt_4=0,
\end{equation}

\vspace{1.5mm}

\begin{equation}
\label{final6}
\lim\limits_{p_1\to\infty}\varlimsup\limits_{p_2\to\infty}
\varlimsup\limits_{p_3\to\infty}\varlimsup\limits_{p_4\to\infty}\int\limits_{[t, T]^3}
\left(R_{p_1 p_2 p_3 p_4}(t_3,t_2,t_4,t_2)\right)^{2n}dt_2dt_3dt_4=0,
\end{equation}

\vspace{1.5mm}

\begin{equation}
\label{final7}
\lim\limits_{p_1\to\infty}\varlimsup\limits_{p_2\to\infty}
\varlimsup\limits_{p_3\to\infty}\varlimsup\limits_{p_4\to\infty}\int\limits_{[t, T]^3}
\left(R_{p_1 p_2 p_3 p_4}(t_3,t_4,t_2,t_2)\right)^{2n}dt_2dt_3dt_4=0,
\end{equation}

\vspace{1.5mm}

\begin{equation}
\label{final8}
\lim\limits_{p_1\to\infty}\varlimsup\limits_{p_2\to\infty}
\varlimsup\limits_{p_3\to\infty}\varlimsup\limits_{p_4\to\infty}\int\limits_{[t, T]^3}
\left(R_{p_1 p_2 p_3 p_4}(t_2,t_2,t_4,t_4)\right)^{2n}dt_2dt_4=0,
\end{equation}

\begin{equation}
\label{final9}
\lim\limits_{p_1\to\infty}\varlimsup\limits_{p_2\to\infty}
\varlimsup\limits_{p_3\to\infty}\varlimsup\limits_{p_4\to\infty}\int\limits_{[t, T]^3}
\left(R_{p_1 p_2 p_3 p_4}(t_2,t_4,t_2,t_4)\right)^{2n}dt_2dt_4=0,
\end{equation}

\begin{equation}
\label{final10}
\lim\limits_{p_1\to\infty}\varlimsup\limits_{p_2\to\infty}
\varlimsup\limits_{p_3\to\infty}\varlimsup\limits_{p_4\to\infty}\int\limits_{[t, T]^3}
\left(R_{p_1 p_2 p_3 p_4}(t_2,t_4,t_4,t_2)\right)^{2n}dt_2dt_4=0.
\end{equation}

\vspace{5mm}

Combaining
(\ref{udar}) with (\ref{final1})--(\ref{final10}), we get

\vspace{1mm}
$$
\lim\limits_{p_1\to\infty}
\varlimsup\limits_{p_2\to\infty}\varlimsup\limits_{p_3\to\infty}
\varlimsup\limits_{p_4\to\infty}
{\sf M}\left\{\left|J[R_{p_1p_2p_3p_4}]_{T,t}^{(4)}\right|^{2n}
\right\}=0,\ \ \ n\in \mathbb{N}.
$$

\vspace{5mm}

Lemma 6 is proved for the case $k=4$.

Let us consider the case of arbitrary $k$ ($k\in\mathbb{N}$).
Let us analyze the stochastic integral defined by
(\ref{30.34}) and find its representation convenient  
for the following consideration. In order to do it we 
introduce several notations. Suppose that

\vspace{2mm}
$$
S_N^{(k)}(a)=
\sum\limits_{j_k=0}^{N-1}\ldots \sum\limits_{j_1=0}^{j_2-1}\ \
\sum\limits_{(j_1,\ldots,j_k)}a_{(j_1,\ldots,j_k)},
$$

\vspace{7mm}

$$
{\rm C}_{s_r}\ldots {\rm C}_{s_1}S_N^{(k)}(a)=
$$

\vspace{1mm}
$$
=
\sum_{j_k=0}^{N-1}\ldots \sum_{j_{s_r+1}=0}^{j_{s_r+2}-1}
\sum_{j_{s_r-1}=0}^{j_{s_r+1}-1}\ldots \sum_{j_{s_1+1}=0}^{j_{s_1+2}-1}
\sum_{j_{s_1-1}=0}^{j_{s_1+1}-1}\ldots \sum_{j_1=0}^{j_2-1}\ \
\sum\limits_{\prod\limits_{l=1}^r{\bf I}_{j_{s_l},j_{s_l+1}}
(j_1,\ldots,j_k)}\ \
a_{\prod\limits_{l=1}^r{\bf I}_{j_{s_l},j_{s_l+1}}
(j_1,\ldots,j_k)},
$$

\vspace{5mm}
\noindent
where

\vspace{-3mm}
$$
\prod\limits_{l=1}^r{\bf I}_{j_{s_l},j_{s_l+1}}
(j_1,\ldots,j_k)\ \ \stackrel{\rm def}{=}\ \ 
{\bf I}_{j_{s_r},j_{s_r+1}}\ldots
{\bf I}_{j_{s_1},j_{s_1+1}}
(j_1,\ldots,j_k),
$$

\vspace{2mm}
$$
{\rm C}_{s_0}\ldots {\rm C}_{s_1}
S_N^{(k)}(a)=S_N^{(k)}(a),\ \ \ \ \ 
\prod\limits_{l=1}^0{\bf I}_{j_{s_l},j_{s_l+1}}
(j_1,\ldots,j_k)=(j_1,\ldots,j_k),
$$

\vspace{3mm}
$$
{\bf I}_{j_l,j_{l+1}}(j_{q_1},\ldots,j_{q_2},j_l,j_{q_3},\ldots,
j_{q_{k-2}},j_l,j_{q_{k-1}},\ldots,j_{q_{k}})\stackrel{\rm def}{=}
$$

\vspace{2mm}
$$
\stackrel{\rm def}{=}(j_{q_1},\ldots,j_{q_2},j_{l+1},j_{q_3},\ldots,
j_{q_{k-2}},j_{l+1},j_{q_{k-1}},\ldots,j_{g_{k}}),
$$

\vspace{7mm}
\noindent
where
$l\in\mathbb{N},\ \
l\ne q_{1},\ldots,q_2,q_3,\ldots,q_{k-2},q_{k-1},\ldots,q_{k},$\ \
$s_1,\ldots,s_r = 1,\ldots, k-1$,\ \ $s_r>\ldots >s_1,$\ \
$a_{(j_{q_1},\ldots,j_{q_k})}$ is a scalar value,\ \
$q_1,\ldots,q_k=1,\ldots,k,$\ \
the expression

\vspace{1mm}
$$
\sum\limits_{(j_{q_1},\ldots,j_{q_k})}
$$

\vspace{4mm}
\noindent
means the sum with respect to  
all possible  
permutations
$(j_{q_1},\ldots,j_{q_k}).$

Using induction it is possible to prove the following equality

\vspace{1mm}
\begin{equation}
\label{979}
\sum_{j_k=0}^{N-1}\ldots \sum_{j_1=0}^{N-1}
a_{(j_1,\ldots,j_k)}
=\sum_{r=0}^{k-1}\ \
\sum_{\stackrel{\Large{s_r,\ldots,s_1=1}}
{{}_{s_r>\ldots>s_1}}}^{k-1}  
{\rm C}_{s_r}\ldots {\rm C}_{s_1}S_N^{(k)}(a),
\end{equation}

\vspace{3mm}
\noindent
where $k=2, 3,\ldots $

Hereinafter in this section, we will identify the following records
$a_{(j_1,\ldots,j_k)}
=a_{(j_1\ldots j_k)}=a_{j_1\ldots j_k}.$
In particular, from (\ref{979}) for $k=2, 3, 4$ 
we get the following formulas

\vspace{3mm}
$$
\sum\limits_{j_2=0}^{N-1}\sum\limits_{j_1=0}^{N-1}
a_{(j_1,j_2)}=S_N^{(2)}(a)
+{\rm C}_{1}S_N^{(2)}(a)=
$$

\vspace{2mm}
$$
=
\sum\limits_{j_2=0}^{N-1}\sum\limits_{j_1=0}^{j_2-1}
\sum\limits_{(j_1,j_2)}a_{(j_1j_2)}+
\sum\limits_{j_2=0}^{N-1}
a_{(j_2j_2)}
=\sum\limits_{j_2=0}^{N-1}\sum\limits_{j_1=0}^{j_2-1}(a_{j_1j_2}+
a_{j_2j_1})+
$$

\vspace{2mm}
\begin{equation}
\label{best1}
+
\sum\limits_{j_2=0}^{N-1}
a_{j_2j_2},
\end{equation}

\vspace{7mm}

$$
\sum\limits_{j_3=0}^{N-1}\sum\limits_{j_2=0}^{N-1}\sum\limits_{j_1=0}^{N-1}
a_{(j_1,j_2,j_3)}=S_N^{(3)}(a)
+{\rm C}_{1}S_N^{(3)}(a)+
{\rm C}_{2}S_N^{(3)}(a)+{\rm C}_{2}{\rm C}_{1}S_N^{(3)}(a)=
$$

\vspace{2mm}
$$
=\sum\limits_{j_3=0}^{N-1}\sum\limits_{j_2=0}^{j_3-1}
\sum\limits_{j_1=0}^{j_2-1}
\sum\limits_{(j_1,j_2,j_3)}a_{(j_1j_2j_3)}+
\sum\limits_{j_3=0}^{N-1}\sum\limits_{j_2=0}^{j_3-1}
\sum\limits_{(j_2,j_2,j_3)}a_{(j_2j_2j_3)}+
$$

\vspace{2mm}
$$
+\sum\limits_{j_3=0}^{N-1}
\sum\limits_{j_1=0}^{j_3-1}
\sum\limits_{(j_1,j_3,j_3)}a_{(j_1j_3j_3)}+
\sum\limits_{j_3=0}^{N-1}a_{(j_3j_3j_3)}=
$$

\vspace{2mm}
$$
=\sum\limits_{j_3=0}^{N-1}\sum\limits_{j_2=0}^{j_3-1}
\sum\limits_{j_1=0}^{j_2-1}
\left(a_{j_1j_2j_3}+a_{j_1j_3j_2}+a_{j_2j_1j_3}+
a_{j_2j_3j_1}+a_{j_3j_2j_1}+a_{j_3j_1j_2}\right)+
$$

\vspace{2mm}
$$
+\sum\limits_{j_3=0}^{N-1}\sum\limits_{j_2=0}^{j_3-1}
\left(a_{j_2j_2j_3}+a_{j_2j_3j_2}+a_{j_3j_2j_2}\right)+
\sum\limits_{j_3=0}^{N-1}
\sum\limits_{j_1=0}^{j_3-1}
\left(a_{j_1j_3j_3}+a_{j_3j_1j_3}+a_{j_3j_3j_1}\right)+
$$

\vspace{2mm}
\begin{equation}
\label{oop12}
+
\sum\limits_{j_3=0}^{N-1}a_{j_3j_3j_3},
\end{equation}

\vspace{7mm}

$$
\sum\limits_{j_4=0}^{N-1}\sum\limits_{j_3=0}^{N-1}
\sum\limits_{j_2=0}^{N-1}\sum\limits_{j_1=0}^{N-1}
a_{(j_1,j_2,j_3,j_4)}=
S_N^{(4)}(a)
+{\rm C}_{1}S_N^{(4)}(a)+{\rm C}_{2}S_N^{(4)}(a)+
$$

\vspace{1mm}
$$
+{\rm C}_{3}S_N^{(4)}(a)+
{\rm C}_{2}{\rm C}_{1}S_N^{(4)}(a)+
{\rm C}_{3}{\rm C}_{1}S_N^{(4)}(a)+
{\rm C}_{3}{\rm C}_{2}S_N^{(4)}(a)+
{\rm C}_{3}{\rm C}_{2}{\rm C}_{1}S_N^{(4)}(a)=
$$

\vspace{1mm}
$$
=\sum\limits_{j_4=0}^{N-1}\sum\limits_{j_3=0}^{j_4-1}
\sum\limits_{j_2=0}^{j_3-1}\sum\limits_{j_1=0}^{j_2-1}
\sum\limits_{(j_1,j_2,j_3,j_4)}a_{(j_1j_2j_3j_4)}
+\sum\limits_{j_4=0}^{N-1}\sum\limits_{j_3=0}^{j_4-1}
\sum\limits_{j_2=0}^{j_3-1}
\sum\limits_{(j_2,j_2,j_3,j_4)}a_{(j_2j_2j_3j_4)}
$$

\vspace{1mm}
$$
+\sum\limits_{j_4=0}^{N-1}\sum\limits_{j_3=0}^{j_4-1}
\sum\limits_{j_1=0}^{j_3-1}
\sum\limits_{(j_1,j_3,j_3,j_4)}a_{(j_1j_3j_3j_4)}
+\sum\limits_{j_4=0}^{N-1}\sum\limits_{j_2=0}^{j_4-1}
\sum\limits_{j_1=0}^{j_2-1}
\sum\limits_{(j_1,j_2,j_4,j_4)}a_{(j_1j_2j_4j_4)}+
$$

\vspace{1mm}
$$
+\sum\limits_{j_4=0}^{N-1}
\sum\limits_{j_3=0}^{j_4-1}
\sum\limits_{(j_3,j_3,j_3,j_4)}a_{(j_3j_3j_3j_4)}
+\sum\limits_{j_4=0}^{N-1}
\sum\limits_{j_2=0}^{j_4-1}
\sum\limits_{(j_2,j_2,j_4,j_4)}a_{(j_2j_2j_4j_4)}+
$$

\vspace{1mm}
$$
+\sum\limits_{j_4=0}^{N-1}
\sum\limits_{j_1=0}^{j_4-1}
\sum\limits_{(j_1,j_4,j_4,j_4)}a_{(j_1j_4j_4j_4)}+
\sum\limits_{j_4=0}^{N-1}a_{j_4j_4j_4j_4}=
$$

\vspace{1mm}
$$
=\sum\limits_{j_4=0}^{N-1}\sum\limits_{j_3=0}^{j_4-1}
\sum\limits_{j_2=0}^{j_3-1}\sum\limits_{j_1=0}^{j_2-1}
\left(a_{j_1j_2j_3j_4}+a_{j_1j_2j_4j_3}+
a_{j_1j_3j_2j_4}+a_{j_1j_3j_4j_2}+\right.
$$

\vspace{1mm}
$$
+a_{j_1j_4j_3j_2}+a_{j_1j_4j_2j_3}+a_{j_2j_1j_3j_4}+
a_{j_2j_1j_4j_3}+
a_{j_2j_4j_1j_3}+a_{j_2j_4j_3j_1}+a_{j_2j_3j_1j_4}+
$$

\vspace{1mm}
$$
+a_{j_2j_3j_4j_1}+a_{j_3j_1j_2j_4}+a_{j_3j_1j_4j_2}+
a_{j_3j_2j_1j_4}+
a_{j_3j_2j_4j_1}+a_{j_3j_4j_1j_2}+a_{j_3j_4j_2j_1}+
$$

\vspace{1mm}
$$
+a_{j_4j_1j_2j_3}+a_{j_4j_1j_3j_2}+a_{j_4j_2j_1j_3}+
\left.a_{j_4j_2j_3j_1}+
a_{j_4j_3j_1j_2}+a_{j_4j_3j_2j_1}\right)+
$$

\vspace{1mm}
$$
+\sum\limits_{j_4=0}^{N-1}\sum\limits_{j_3=0}^{j_4-1}
\sum\limits_{j_2=0}^{j_3-1}
\left(a_{j_2j_2j_3j_4}+a_{j_2j_2j_4j_3}+a_{j_2j_3j_2j_4}+\right.
a_{j_2j_4j_2j_3}+
a_{j_2j_3j_4j_2}+a_{j_2j_4j_3j_2}+
$$

\vspace{1mm}
$$
+a_{j_3j_2j_2j_4}+a_{j_4j_2j_2j_3}+a_{j_3j_2j_4j_2}
\left.+
a_{j_4j_2j_3j_2}+
a_{j_4j_3j_2j_2}+a_{j_3j_4j_2j_2}\right)+
$$

\vspace{1mm}
$$
+\sum\limits_{j_4=0}^{N-1}\sum\limits_{j_3=0}^{j_4-1}
\sum\limits_{j_1=0}^{j_3-1}
\left(a_{j_3j_3j_1j_4}+a_{j_3j_3j_4j_1}+a_{j_3j_1j_3j_4}+\right.
a_{j_3j_4j_3j_1}+
a_{j_3j_4j_1j_3}+a_{j_3j_1j_4j_3}+
$$

\vspace{1mm}
$$
+a_{j_1j_3j_3j_4}+a_{j_4j_3j_3j_1}+a_{j_4j_3j_1j_3}
\left.+a_{j_1j_3j_4j_3}+
a_{j_1j_4j_3j_3}+a_{j_4j_1j_3j_3}\right)+
$$

\vspace{1mm}
$$
+\sum\limits_{j_4=0}^{N-1}\sum\limits_{j_2=0}^{j_4-1}
\sum\limits_{j_1=0}^{j_2-1}
\left(a_{j_4j_4j_1j_2}+a_{j_4j_4j_2j_1}+a_{j_4j_1j_4j_2}+\right.
a_{j_4j_2j_4j_1}+
a_{j_4j_2j_1j_4}+a_{j_4j_1j_2j_4}+
$$

\vspace{1mm}
$$
+a_{j_1j_4j_4j_2}+a_{j_2j_4j_4j_1}+a_{j_2j_4j_1j_4}+
\left.a_{j_1j_4j_2j_4}+
a_{j_1j_2j_4j_4}+a_{j_2j_1j_4j_4}\right)+
$$

\vspace{1mm}
$$
+\sum\limits_{j_4=0}^{N-1}
\sum\limits_{j_3=0}^{j_4-1}
\left(a_{j_3j_3j_3j_4}+a_{j_3j_3j_4j_3}+a_{j_3j_4j_3j_3}+
a_{j_4j_3j_3j_3}\right)+
$$

\vspace{1mm}
$$
+\sum\limits_{j_4=0}^{N-1}
\sum\limits_{j_2=0}^{j_4-1}
\left(a_{j_2j_2j_4j_4}+a_{j_2j_4j_2j_4}+a_{j_2j_4j_4j_2}+\right.\left.
a_{j_4j_2j_2j_4}+
a_{j_4j_2j_4j_2}+a_{j_4j_4j_2j_2}\right)+
$$

\vspace{1mm}
$$
+\sum\limits_{j_4=0}^{N-1}
\sum\limits_{j_1=0}^{j_4-1}
\left(a_{j_1j_4j_4j_4}+a_{j_4j_1j_4j_4}+a_{j_4j_4j_1j_4}+
a_{j_4j_4j_4j_1}\right)+
$$

\vspace{1mm}
\begin{equation}
\label{huh}
+\sum\limits_{j_4=0}^{N-1}a_{j_4j_4j_4j_4}.
\end{equation}

\vspace{6mm}

Perhaps, the formula (\ref{979}) for any $k$ ($k\in \mathbb{N}$) was 
found
by the author for the first time \cite{1997} (1997).

Assume that 

\vspace{-3mm}
$$
a_{(j_1,\ldots,j_k)}=
\Phi\left(\tau_{j_1},\ldots,\tau_{j_k}\right)
\prod\limits_{l=1}^k\Delta{\bf w}_{\tau_{j_l}}^{(i_l)},
$$ 

\vspace{4mm}
\noindent
where $\Phi\left(t_1,\ldots,t_k\right)$ is a nonrandom
function of $k$ variables.
Then from (\ref{30.34})
and (\ref{979}) we have

\vspace{5mm}
$$
J[\Phi]_{T,t}^{(k)}=\sum_{r=0}^{[k/2]}\ \ 
\sum_{(s_r,\ldots,s_1)\in{\rm A}_{k,r}}  \times
$$

\vspace{5mm}
$$
\times\ \ 
\hbox{\vtop{\offinterlineskip\halign{
\hfil#\hfil\cr
{\rm l.i.m.}\cr
$\stackrel{}{{}_{N\to \infty}}$\cr
}} }
\sum_{j_k=0}^{N-1}\ldots \sum_{j_{s_r+1}=0}^{j_{s_r+2}-1}
\sum_{j_{s_r-1}=0}^{j_{s_r+1}-1}\ldots \sum_{j_{s_1+1}=0}^{j_{s_1+2}-1}
\sum_{j_{s_1-1}=0}^{j_{s_1+1}-1}\ldots \sum_{j_1=0}^{j_2-1}\ \ 
\sum\limits_{\prod\limits_{l=1}^r{\bf I}_{j_{s_l},j_{s_l+1}}
(j_1,\ldots,j_k)}\times
$$

\vspace{3mm}
$$
\times
\Biggl[\Phi\biggl(\tau_{j_1},\ldots,\tau_{j_{s_1-1}},
\tau_{j_{s_1+1}},\tau_{j_{s_1+1}},\tau_{j_{s_1+2}},\ldots
,\tau_{j_{s_r-1}},\tau_{j_{s_r+1}},\tau_{j_{s_r+1}},\tau_{j_{s_r+2}},
\ldots,\tau_{j_k}\biggr)
\times\Biggr.
$$

\vspace{4mm}
$$
\times
\Delta{\bf w}_{\tau_{j_1}}^{(i_1)}
\ldots\Delta{\bf w}_{\tau_{j_{s_1-1}}}^{(i_{s_1-1})}
\Delta{\bf w}_{\tau_{j_{s_1+1}}}^{(i_{s_1})}
\Delta{\bf w}_{\tau_{j_{s_1+1}}}^{(i_{s_1+1})}
\Delta{\bf w}_{\tau_{j_{s_1+2}}}^{(i_{s_1+2})}
\ldots
$$

\vspace{2mm}
$$
\Biggl.
\ldots\Delta{\bf w}_{\tau_{j_{s_r-1}}}^{(i_{s_r-1})}
\Delta{\bf w}_{\tau_{j_{s_r+1}}}^{(i_{s_r})}
\Delta{\bf w}_{\tau_{j_{s_r+1}}}^{(i_{s_r+1})}
\Delta{\bf w}_{\tau_{j_{s_r+2}}}^{(i_{s_r+2})}\ldots
\Delta{\bf w}_{\tau_{j_{k}}}^{(i_{k})}\Biggr]=
$$

\vspace{6mm}
\begin{equation}
\label{30.51}
=\sum_{r=0}^{[k/2]}\sum_{(s_r,\ldots,s_1)\in{\rm A}_{k,r}}
I[\Phi]_{T,t}^{(k)s_1,\ldots,s_r}\ \ \ \hbox{w.\ p.\ 1},
\end{equation}

\vspace{3mm}
\noindent
where

\vspace{3mm}
$$
I[\Phi]_{T,t}^{(k)s_1,\ldots,s_r}=
\int\limits_t^T\ldots\int\limits_t^{t_{s_r+3}}
\int\limits_t^{t_{s_r+2}}\int\limits_t^{t_{s_r}}\ldots
\int\limits_t^{t_{s_1+3}}
\int\limits_t^{t_{s_1+2}}\int\limits_t^{t_{s_1}}\ldots\int\limits_t^{t_2}
\sum\limits_{\prod\limits_{l=1}^r{\bf I}_{t_{s_l},t_{s_l+1}}
(t_1,\ldots,t_k)}\times
$$

\vspace{3mm}
$$
\times
\Biggl[\Phi\biggl(t_{1},\ldots,t_{s_1-1},t_{s_1+1},t_{s_1+1},t_{s_1+2},\ldots
,t_{s_r-1},t_{s_r+1},t_{s_r+1},t_{s_r+2},
\ldots,t_{k}\biggr)\times\Biggr.
$$

\vspace{4mm}
$$
\times
d{\bf w}_{t_1}^{(i_1)}\ldots d{\bf w}_{t_{s_1-1}}^{(i_{s_1-1})}
d{\bf w}_{t_{s_1+1}}^{(i_{s_1})}
d{\bf w}_{t_{s_1+1}}^{(i_{s_1+1})}
d{\bf w}_{t_{s_1+2}}^{(i_{s_1+2})}
\ldots 
$$

\vspace{2mm}
\begin{equation}
\label{99999}
\Biggl.
\ldots
d{\bf w}_{t_{s_r-1}}^{(i_{s_r-1})}
d{\bf w}_{t_{s_r+1}}^{(i_{s_r})}
d{\bf w}_{t_{s_r+1}}^{(i_{s_r+1})}
d{\bf w}_{t_{s_r+2}}^{(i_{s_r+2})}\ldots
d{\bf w}_{t_k}^{(i_k)}\Biggr],
\end{equation}

\vspace{5mm}
\noindent
where 
$k\ge 2,$ the set ${\rm A}_{k,r}$ is defined by (\ref{30.5550001}).
We suppose that
the right-hand side of (\ref{99999}) exists as the Ito stochastic integral.

\vspace{2mm}

{\bf Remark 2.}\ {\it The summands on the right-hand side
of {\rm (\ref{99999})} should be 
understood
as follows{\rm :} for each 
permutation
from the set 

$$
\prod\limits_{l=1}^r{\bf I}_{t_{s_l},t_{s_l+1}}
(t_1,\ldots,t_k)
=
(t_{1},\ldots,t_{s_1-1},t_{s_1+1},t_{s_1+1},t_{s_1+2},\ldots
,t_{s_r-1},t_{s_r+1},t_{s_r+1},t_{s_r+2},
\ldots,t_{k})
$$

\vspace{3mm}
\noindent
it is necessary to perform replacement on the right-hand side 
of {\rm (\ref{99999})} 
of all pairs {\rm (}their number is equals to 
$r${\rm )} of differentials $d{\bf w}_{t_p}^{(i)}d{\bf w}_{t_p}^{(j)}$ 
with similar 
lower indexes 
by the values ${\bf 1}_{\{i=j\ne 0\}}dt_p.$}

\vspace{2mm}

Note that the term in (\ref{30.51}) for $r=0$ should be understood as
follows

\begin{equation}
\label{pobeda}
\int\limits_t^T\ldots \int\limits_t^{t_2}
\sum\limits_{(t_1,\ldots,t_k)}
\biggl(\Phi\left(t_1,\ldots,t_k\right)
d{\bf w}_{t_{1}}^{(i_{1})}\ldots
d{\bf w}_{t_k}^{(i_k)}\biggr),
\end{equation}

\vspace{4mm}
\noindent
where
$$
\sum\limits_{(t_1,\ldots,t_k)}
$$ 

\vspace{2mm}
\noindent
means the sum with respect to all
possible permutations
$(t_1,\ldots,t_k).$ 
At the same time 
permutations $(t_1,\ldots,t_k)$ when summing
are performed in (\ref{pobeda}) only in the expression, which
is enclosed in pa\-ren\-the\-ses (see \cite{10a}, Sect.~1.1.3 for details).

Using Lemma 4, we get

\vspace{1mm}
$$
{\sf M}\left\{\left|J[\Phi]_{T,t}^{(k)}\right|^{2n}\right\}\le
$$

\vspace{3mm}
\begin{equation}
\label{333.225}
\le C_{nk}\ \ \sum_{r=0}^{[k/2]}\sum_{(s_r,\ldots,s_1)\in{\rm A}_{k,r}}
{\sf M}\left\{\left|
I[\Phi]_{T,t}^{(k)s_1,\ldots,s_r}\right|^{2n}\right\},
\end{equation}

\vspace{3mm}
\noindent
where

\vspace{1mm}
$$
{\sf M}\left\{\left|
I[\Phi]_{T,t}^{(k)s_1,\ldots,s_r}\right|^{2n}\right\}\le
$$

\vspace{5mm}
$$
\le C_{nk}^{s_1\ldots s_r}
\int\limits_t^T\ldots\int\limits_t^{t_{s_r+3}}
\int\limits_t^{t_{s_r+2}}\int\limits_t^{t_{s_r}}\ldots
\int\limits_t^{t_{s_1+3}}
\int\limits_t^{t_{s_1+2}}\int\limits_t^{t_{s_1}}\ldots\int\limits_t^{t_2}
\sum\limits_{\prod\limits_{l=1}^r{\bf I}_{t_{s_l},t_{s_l+1}}
(t_1,\ldots,t_k)}\times
$$

\vspace{3mm}
$$
\times
\Phi^{2n}\biggl(t_{1},\ldots,t_{s_1-1},t_{s_1+1},t_{s_1+1},t_{s_1+2},\ldots,
t_{s_r-1},t_{s_r+1},t_{s_r+1},t_{s_r+2},\ldots,t_{k}\biggr)\times
$$

\vspace{2mm}
\begin{equation}
\label{333.226}
\times
dt_1\ldots dt_{s_1-1}
dt_{s_1+1}dt_{s_1+2}
\ldots dt_{s_r-1}dt_{s_r+1}
dt_{s_r+2}\ldots
dt_k,
\end{equation}
\vspace{4mm}

\noindent
where $C_{nk}$ and $C_{nk}^{s_1\ldots s_r}$ are constants
and
permutations
when
summing in
(\ref{333.226}) are performed only in the value

$$
\Phi^{2n}\biggl(t_{1},\ldots,t_{s_1-1},t_{s_1+1},t_{s_1+1},t_{s_1+2},\ldots,
t_{s_r-1},t_{s_r+1},t_{s_r+1},t_{s_r+2},\ldots,t_{k}\biggr).
$$

\vspace{6mm}

Consider (\ref{333.225}), (\ref{333.226}) for
$\Phi(t_1,\ldots,t_k)=R_{p_1\ldots p_k}(t_1,\ldots,t_k)$

\vspace{7mm}

$$
{\sf M}\left\{\left|J[R_{p_1\ldots p_k}]_{T,t}^{(k)}\right|^{2n}\right\}
\le 
$$

\vspace{4mm}
\begin{equation}
\label{333.225e}
\le
C_{nk}\ \ \sum_{r=0}^{[k/2]}\sum_{(s_r,\ldots,s_1)\in{\rm A}_{k,r}}
{\sf M}\left\{\left|
I[R_{p_1\ldots p_k}]_{T,t}^{(k)s_1,\ldots,s_r}\right|^{2n}\right\},
\end{equation}

\vspace{9mm}

$$
{\sf M}\left\{\left|
I[R_{p_1\ldots p_k}]_{T,t}^{(k)s_1,\ldots,s_r}\right|^{2n}\right\}\le
$$

\vspace{5mm}
$$
\le C_{nk}^{s_1\ldots s_r}
\int\limits_t^T\ldots\int\limits_t^{t_{s_r+3}}
\int\limits_t^{t_{s_r+2}}\int\limits_t^{t_{s_r}}\ldots
\int\limits_t^{t_{s_1+3}}
\int\limits_t^{t_{s_1+2}}\int\limits_t^{t_{s_1}}\ldots\int\limits_t^{t_2}
\sum\limits_{\prod\limits_{l=1}^r{\bf I}_{t_{s_l},t_{s_l+1}}
(t_1,\ldots,t_k)}\times
$$

\vspace{4mm}
$$
\times
R_{p_1\ldots p_k}^{2n}
\biggl(t_{1},\ldots,t_{s_1-1},t_{s_1+1},t_{s_1+1},t_{s_1+2},\ldots,
t_{s_r-1},t_{s_r+1},t_{s_r+1},t_{s_r+2},\ldots,t_{k}\biggr)\times
$$

\vspace{2mm}
\begin{equation}
\label{333.226e}
\times
dt_1\ldots dt_{s_1-1}
dt_{s_1+1}dt_{s_1+2}
\ldots dt_{s_r-1}dt_{s_r+1}
dt_{s_r+2}\ldots
dt_k,
\end{equation}
\vspace{5mm}

\noindent
where $C_{nk}$ and $C_{nk}^{s_1\ldots s_r}$ are constants and
permutations
when summing
in
(\ref{333.226e}) are performed only in the value

\vspace{1mm}
$$
R_{p_1\ldots p_k}^{2n}
\biggl(t_{1},\ldots,t_{s_1-1},t_{s_1+1},t_{s_1+1},t_{s_1+2},\ldots,
t_{s_r-1},t_{s_r+1},t_{s_r+1},t_{s_r+2},\ldots,t_{k}\biggr).
$$

\vspace{4mm}

From the other hand, we can consider
the generalization of the formulas (\ref{leto80010aaa}),
(\ref{oop16xxx}), (\ref{udar}) 
for the case of arbitrary $k$ ($k\in\mathbb{N}$).
In order to do this, let us
consider the unordered
set $\{1, 2, \ldots, k\}$ 
and separate it into two parts:
the first part consists of $r$ unordered 
pairs (sequence order of these pairs is also unimportant) and the 
second one consists of the 
remaining $k-2r$ numbers.
So, we have

\vspace{2mm}
\begin{equation}
\label{leto5007}
(\{
\underbrace{\{g_1, g_2\}, \ldots, 
\{g_{2r-1}, g_{2r}\}}_{\small{\hbox{part 1}}}
\},
\{\underbrace{q_1, \ldots, q_{k-2r}}_{\small{\hbox{part 2}}}
\}),
\end{equation}

\vspace{5mm}
\noindent
where 
$$
\{g_1, g_2, \ldots, 
g_{2r-1}, g_{2r}, q_1, \ldots, q_{k-2r}\}=\{1, 2, \ldots, k\},
$$

\vspace{5mm}
\noindent
braces   
mean an unordered 
set, and pa\-ren\-the\-ses mean an ordered set.

We will say that (\ref{leto5007}) is a partition 
and consider the sum with respect to all possible
partitions

\vspace{1mm}
\begin{equation}
\label{leto5008}
\sum_{\stackrel{(\{\{g_1, g_2\}, \ldots, 
\{g_{2r-1}, g_{2r}\}\}, \{q_1, \ldots, q_{k-2r}\})}
{{}_{\{g_1, g_2, \ldots, 
g_{2r-1}, g_{2r}, q_1, \ldots, q_{k-2r}\}=\{1, 2, \ldots, k\}}}}
a_{g_1 g_2, \ldots, 
g_{2r-1} g_{2r}, q_1 \ldots q_{k-2r}},
\end{equation}

\vspace{5mm}
\noindent
where $a_{g_1 g_2, \ldots, 
g_{2r-1} g_{2r}, q_1 \ldots q_{k-2r}}\in \mathbb{R}.$

Below there are several examples of sums in the form (\ref{leto5008})

\vspace{1mm}
$$
\sum_{\stackrel{(\{g_1, g_2\})}{{}_{\{g_1, g_2\}=\{1, 2\}}}}
a_{g_1 g_2}=a_{12},
$$

\vspace{6mm}
\begin{equation}
\label{f110}
\sum_{\stackrel{(\{\{g_1, g_2\}, \{g_3, g_4\}\})}
{{}_{\{g_1, g_2, g_3, g_4\}=\{1, 2, 3, 4\}}}}
a_{g_1 g_2, g_3 g_4}=a_{12,34} + a_{13,24} + a_{23,14},
\end{equation}

\vspace{6mm}
$$
\sum_{\stackrel{(\{g_1, g_2\}, \{q_1, q_{2}\})}
{{}_{\{g_1, g_2, q_1, q_{2}\}=\{1, 2, 3, 4\}}}}
a_{g_1 g_2, q_1 q_{2}}=
$$

\begin{equation}
\label{f111}
=a_{12,34}+a_{13,24}+a_{14,23}
+a_{23,14}+a_{24,13}+a_{34,12},
\end{equation}

\vspace{6mm}
$$
\sum_{\stackrel{(\{g_1, g_2\}, \{q_1, q_{2}, q_3\})}
{{}_{\{g_1, g_2, q_1, q_{2}, q_3\}=\{1, 2, 3, 4, 5\}}}}
a_{g_1 g_2, q_1 q_{2}q_3}
=
$$

$$
=a_{12,345}+a_{13,245}+a_{14,235}
+a_{15,234}+a_{23,145}+a_{24,135}+
$$

$$
+a_{25,134}+a_{34,125}+a_{35,124}+a_{45,123},
$$

\vspace{6mm}
$$
\sum_{\stackrel{(\{\{g_1, g_2\}, \{g_3, g_{4}\}\}, \{q_1\})}
{{}_{\{g_1, g_2, g_3, g_{4}, q_1\}=\{1, 2, 3, 4, 5\}}}}
a_{g_1 g_2, g_3 g_{4},q_1}
=
$$

$$
=
a_{12,34,5}+a_{13,24,5}+a_{14,23,5}+
a_{12,35,4}+a_{13,25,4}+a_{15,23,4}+
$$

$$
+a_{12,54,3}+a_{15,24,3}+a_{14,25,3}+a_{15,34,2}+a_{13,54,2}+a_{14,53,2}+
$$

$$
+
a_{52,34,1}+a_{53,24,1}+a_{54,23,1}.
$$

\vspace{9mm}

Now we can 
generalize the formulas (\ref{leto80010aaa}),
(\ref{oop16xxx}), (\ref{udar}) 
for the case of arbitrary $k$ ($k\in\mathbb{N}$)

\vspace{4mm}

$$
{\sf M}\left\{\left|J[R_{p_1\ldots p_k}]_{T,t}^{(k)}\right|^{2n}\right\}
\le C_{nk}\left(
\int\limits_{[t, T]^k}
\left(R_{p_1\ldots p_k}(t_1,\ldots,t_k)\right)^{2n}dt_1\ldots dt_k+
\right.
$$

\vspace{3mm}

$$
+\sum\limits_{r=1}^{[k/2]}
\sum_{\stackrel{(\{\{g_1, g_2\}, \ldots, 
\{g_{2r-1}, g_{2r}\}\}, \{q_1, \ldots, q_{k-2r}\})}
{{}_{\{g_1, g_2, \ldots, 
g_{2r-1}, g_{2r}, q_1, \ldots, q_{k-2r}\}=\{1, 2, \ldots, k\}}}}
{\bf 1}_{\{i_{g_{{}_{1}}}=i_{g_{{}_{2}}}\ne 0\}}
\ldots 
{\bf 1}_{\{i_{g_{{}_{2r-1}}}=i_{g_{{}_{2r}}}\ne 0\}}\times
$$

\vspace{3mm}

$$
\times
\int\limits_{[t, T]^{k-r}}
\left(R_{p_1\ldots p_k}\biggl(
t_1,\ldots,t_k\biggr)\biggl.\biggr|_{t_{g_{{}_{1}}}=t_{g_{{}_{2}}},\ldots,
t_{g_{{}_{2r-1}}}=t_{g_{{}_{2r}}}}
\right)^{2n}\times
$$

\vspace{3mm}
\begin{equation}
\label{udar1}
\left.\times \biggl(dt_1\ldots dt_k\biggr)\Biggl.
\Biggr|_{\left(dt_{g_{{}_{1}}}dt_{g_{{}_{2}}}\right)
\curvearrowright dt_{g_1},\ldots,
\left(dt_{g_{{}_{2r-1}}}dt_{g_{{}_{2r}}}\right)\curvearrowright
dt_{g_{{}_{2r-1}}}}\right),
\end{equation}

\vspace{10mm}
\noindent
where $C_{nk}$ is a constant,

\vspace{-1mm}
$$
\biggl(
t_1,\ldots,t_k\biggr)\biggl.\biggr|_{t_{g_{{}_{1}}}=t_{g_{{}_{2}}},\ldots,
t_{g_{{}_{2r-1}}}=t_{g_{{}_{2r}}}}
$$

\vspace{5mm}
\noindent
means the ordered set $(t_1,\ldots,t_k)$ where we put
$t_{g_{{}_{1}}}=t_{g_{{}_{2}}},$ $\ldots,$
$t_{g_{{}_{2r-1}}}=t_{g_{{}_{2r}}}.$

Moreover,

\vspace{-1mm}
$$
\biggl(dt_1\ldots dt_k\biggr)\Biggl.
\Biggr|_{\left(dt_{g_{{}_{1}}}dt_{g_{{}_{2}}}\right)
\curvearrowright dt_{g_1},\ldots,
\left(dt_{g_{{}_{2r-1}}}dt_{g_{{}_{2r}}}\right)\curvearrowright
dt_{g_{{}_{2r-1}}}}
$$

\vspace{5mm}
\noindent 
means the product $dt_1\ldots dt_k$ where we replace
all pairs 
$dt_{g_{{}_{1}}}dt_{g_{{}_{2}}},$ $\ldots,$ 
$dt_{g_{{}_{2r-1}}}dt_{g_{{}_{2r}}}$ by 
$dt_{g_1},$ $\ldots,$ $dt_{g_{{}_{2r-1}}}$
correspondingly.

Note that the estimate like (\ref{udar1}), 
where all indicators ${\bf 1}_{\{\cdot\}}$ must be 
replaced with $1$, can be obtained from 
the estimates (\ref{333.225e}), (\ref{333.226e}).
The comparison of (\ref{udar1}) with the relation (5.36)
in \cite{10} (Theorem 5.2, p.~A.273) or with 
the relation (1.54) in \cite{10a} (Theorem~1.2, p.~60) shows
a similar structure of these formulas.

Let us consider the particular case of (\ref{udar1}) for $k=4$

\vspace{2mm}

$$
{\sf M}\left\{\left|J[R_{p_1 p_2 p_3 p_4}]_{T,t}^{(4)}\right|^{2n}\right\}
\le C_{n4}\left(
\int\limits_{\stackrel{~}{[t, T]^4}}
\left(R_{p_1 p_2 p_3 p_4}(t_1,t_2,t_3,t_4)\right)^{2n}dt_1 dt_2 dt_3 dt_4+
\right.
$$

\vspace{3mm}
$$
+
\sum_{\stackrel{(\{g_1, g_2\}, \{q_1, q_{2}\})}
{{}_{\{g_1, g_2, q_1, q_{2}\}=\{1, 2, 3, 4\}}}}
{\bf 1}_{\{i_{g_{{}_{1}}}=i_{g_{{}_{2}}}\ne 0\}}
\int\limits_{\stackrel{~}{[t, T]^3}}
\left(R_{p_1 p_2 p_3 p_4}\biggl(
t_1,t_2,t_3,t_4\biggr)\biggl.\biggr|_{t_{g_{{}_{1}}}=t_{g_{{}_{2}}}}
\right)^{2n}
\times
$$

\vspace{3mm}
$$
\times
\biggl(dt_1 dt_2 dt_3 dt_4\biggr)
\Biggl.\Biggr|_{\left(dt_{g_{{}_{1}}}dt_{g_{{}_{2}}}\right)
\curvearrowright dt_{g_1}}+
$$

\vspace{3mm}
$$
+
\sum_{\stackrel{(\{\{g_1, g_2\}, \{g_3, g_{4}\}\})}
{{}_{\{g_1, g_2, g_3, g_{4}\}=\{1, 2, 3, 4\}}}}
{\bf 1}_{\{i_{g_{{}_{1}}}=i_{g_{{}_{2}}}\ne 0\}}
{\bf 1}_{\{i_{g_{{}_{3}}}=i_{g_{{}_{4}}}\ne 0\}}
\int\limits_{\stackrel{~}{[t, T]^2}}
\left(R_{p_1 p_2 p_3 p_4}\biggl(
t_1,t_2,t_3,t_4\biggr)\biggl.\biggr|_{t_{g_{{}_{1}}}=t_{g_{{}_{2}}},
t_{g_{{}_{3}}}=t_{g_{{}_{4}}}}
\right)^{2n}
\times
$$

\begin{equation}
\label{f112}
\times\left.
\biggl(dt_1 dt_2 dt_3 dt_4\biggr)
\biggl|_{\left(dt_{g_{{}_{1}}}dt_{g_{{}_{2}}}\right)\biggr.
\curvearrowright dt_{g_1},
\left(dt_{g_{{}_{3}}}dt_{g_{{}_{4}}}\right)
\curvearrowright dt_{g_3}}\right).
\end{equation}

\vspace{8mm}

It is not difficult to notice that (\ref{f112}) is consistent with 
(\ref{udar})
(see (\ref{f110}), (\ref{f111})).

According to (\ref{1999.1}), we have the following 
expression for all internal points of the hypercube $[t,T]^k$

\vspace{2mm}
$$
R_{p_1\ldots p_k}(t_1,\ldots,t_k)=
$$

\vspace{2mm}
$$
=
\prod_{l=1}^k \psi_l(t_l)\left(\prod_{l=1}^{k-1}
{\bf 1}_{\{t_l<t_{l+1}\}}+
\sum_{r=1}^{k-1}\frac{1}{2^r}
\sum_{\stackrel{s_r,\ldots,s_1=1}{{}_{s_r>\ldots>s_1}}}^{k-1}
\prod_{l=1}^r {\bf 1}_{\{t_{s_l}=t_{s_l+1}\}}
\prod_{\stackrel{l=1}{{}_{l\ne s_1,\ldots, s_r}}}^{k-1}
{\bf 1}_{\{t_{l}<t_{l+1}\}}\right)-
$$

\vspace{2mm}
\begin{equation}
\label{30.48}
-\sum_{j_1=0}^{p_1}\ldots\sum_{j_k=0}^{p_k}
C_{j_k\ldots j_1} \prod_{l=1}^{k} \phi_{j_l}(t_l).
\end{equation}

\vspace{5mm}

Due to (\ref{30.48}) 
the function $R_{p_1\ldots p_k}(t_1,\ldots,t_k)$
is continuous 
in the open
domains of integration of integrals on the right-hand side
of (\ref{333.226e}) 
and it is bounded
at the 
boundaries  
of these 
domains for $p_1,\ldots,p_k<\infty.$

Let us perform the iterated  passage to the limit

\vspace{-1mm}
$$
\lim\limits_{p_1\to\infty}\varlimsup\limits_{p_2\to\infty}\ldots 
\varlimsup\limits_{p_k\to\infty}
$$

\vspace{3mm}
\noindent
under the  
integral signs in the estimate (\ref{udar1})
(like it was performed for 
the 2-dimensional, 3-dimentional, and 4-dimensional cases (see above)).
Then, taking into account
(\ref{410}), we obtain the required 
result. More precisely,
since the integrals on the right-hand side of (\ref{udar1}) 
exist as Riemann integrals, then they are equal to the 
corresponding Lebesgue integrals. 
Moreover, 

\vspace{1mm}
$$
\lim\limits_{p_1\to\infty}\ldots \lim\limits_{p_k\to\infty}
R_{p_1\ldots p_k}(t_1,\ldots,t_k)=0,\ \ \ \hbox{when}\ \ \ (t_1,\ldots,t_k)\in (t, T)^k,
$$

\vspace{4mm}
\noindent
where the left-hand side 
is bounded on $[t, T]^k.$

According to the proof of Lemma 1 and (\ref{30.46}), we have

\vspace{2mm}
$$
R_{p_1\ldots p_k}(t_1,\ldots,t_k)=
$$

\vspace{2mm}
$$
=
\left(K^{*}(t_1,\ldots,t_k)-\sum\limits_{j_1=0}^{p_1}
C_{j_1}(t_2,\ldots,t_k)\phi_{j_1}(t_1)\right)+
$$

\vspace{2mm}
$$
+\left(
\sum\limits_{j_1=0}^{p_1}\left(C_{j_1}(t_2,\ldots,t_k)-
\sum\limits_{j_2=0}^{p_2}
C_{j_2j_1}(t_3,\ldots,t_k)\phi_{j_2}(t_2)\right)
\phi_{j_1}(t_1)\right)+
$$

$$
\ldots
$$

\begin{equation}
\label{strange1}
+\left(
\sum\limits_{j_1=0}^{p_1}\ldots 
\sum\limits_{j_{k-1}=0}^{p_{k-1}}\left(C_{j_{k-1}\ldots j_1}(t_k)-
\sum\limits_{j_k=0}^{p_k}
C_{j_k\ldots j_1}\phi_{j_k}(t_k)\right)
\phi_{j_{k-1}}(t_{k-1})\ldots \phi_{j_1}(t_1)\right),
\end{equation}

\vspace{6mm}
\noindent
where
$$
C_{j_1}(t_2,\ldots,t_k)=\int\limits_t^T
K^{*}(t_1,\ldots,t_k)\phi_{j_1}(t_1)dt_1,
$$

\vspace{2mm}
$$
C_{j_2j_1}(t_3,\ldots,t_k)=\int\limits_{[t, T]^2}
K^{*}(t_1,\ldots,t_k)\phi_{j_1}(t_1)\phi_{j_2}(t_2)dt_1 dt_2,
$$

$$
\ldots
$$

$$
C_{j_{k-1}\ldots j_1}(t_k)=\int\limits_{[t, T]^{k-1}}
K^{*}(t_1,\ldots,t_k)\prod\limits_{l=1}^{k-1}
\phi_{j_l}(t_l)dt_1\ldots dt_{k-1}.
$$

\vspace{6mm}

Then, applying $k$ times (we mean an iterated passage to the limit
$\lim\limits_{p_1\to\infty}\varlimsup\limits_{p_2\to\infty}\ldots 
\varlimsup\limits_{p_k\to\infty}$)
the Lebesgue's Dominated Convergence Theorem in the integrals
on the right-hand side of (\ref{udar1}),
we obtain

\vspace{-2mm}
$$
\lim\limits_{p_1\to\infty}
\varlimsup\limits_{p_2\to\infty}\ldots \varlimsup\limits_{p_k\to\infty}
{\sf M}\left\{\left|J[R_{p_1\ldots p_k}]_{T,t}^{(k)}\right|^{2n}
\right\}=0,\ \ \ n\in \mathbb{N}.
$$

\vspace{5mm}

Let us discuss the choice of integrable majorants
when applying Lebesgue's 
Dominated Convergence Theorem 
in (\ref{udar1}).

It is well known that \cite{bari}

\vspace{-3mm}
\begin{equation}
\label{strange2}
\left|\sum\limits_{k=1}^N \frac{\sin kx}{k}\right|\le C
\end{equation}

\vspace{5mm}
\noindent
for all $N$ and $x$, where constant $C$ does not depend on $N$ and $x$.

Moreover,

\vspace{-3mm}
\begin{equation}
\label{strange3}
\sum\limits_{j=1}^N \frac{1}{j^2}\le \sum\limits_{j=1}^{\infty} \frac{1}{j^2}=
\frac{\pi^2}{6}.
\end{equation}

\vspace{5mm}

Applying double integration by parts (as in (2.28), Sect.~2.1.1  \cite{10a}), we estimate
the partial sums of one-dimensional trigonometric Fourier series

\vspace{1mm}
$$
\sum\limits_{j_1=0}^{p_1}
C_{j_1}(t_2,\ldots,t_k)\phi_{j_1}(t_1),\ \ \ 
\sum\limits_{j_2=0}^{p_2}
C_{j_2j_1}(t_3,\ldots,t_k)\phi_{j_2}(t_2),\ \ \ 
\ldots\ \ \ 
\sum\limits_{j_k=0}^{p_k}
C_{j_k\ldots j_1}\phi_{j_k}(t_k)
$$ 

\vspace{4mm}
\noindent
in (\ref{strange1}) using (\ref{strange3}) and (see (\ref{strange2}))

\vspace{1mm}
$$
\left|\sum\limits_{k=1}^N \frac{1}{k}\sin\frac{2\pi k (x-y)}{T-t}\right|\le C,\ \ \ 
\left|\sum\limits_{k=1}^N \frac{1}{k}\sin\frac{2\pi k (x-t)}{T-t}\right|\le C
$$

\vspace{5mm}
\noindent 
(here $N\in \mathbb{N}$ and $x, y\in \mathbb{R}$, constant $C$ does not depend on $N$ and $x, y$) as follows

$$
\left|\sum\limits_{j_1=0}^{p_1}
C_{j_1}(t_2,\ldots,t_k)\phi_{j_1}(t_1)\right|\le C_1,\ \ \ 
\left|\sum\limits_{j_1=0}^{p_1}
C_{j_1}(t_2,\ldots,t_k)\phi_{j_1}(t_1)\right|\le C_2,\ \ 
\ldots\ \ 
\left|\sum\limits_{j_k=0}^{p_k}
C_{j_k\ldots j_1}\phi_{j_k}(t_k)\right|\le C_k,
$$

\vspace{4mm}
\noindent
where constant $C_1$ does not depend on $p_1,$ constant $C_2$ does not depend on $p_2,$ etc.

Moreover, 

\vspace{-1mm}
$$
\left|K^{*}(t_1,\ldots,t_k)\right|\le \tilde C_1,\ \ \ 
\left|C_{j_1}(t_2,\ldots,t_k)\right|\le \tilde C_2,\ \ \ 
\ldots\ \ \
\left|C_{j_{k-1}\ldots j_1}(t_k)\right|\le \tilde C_k,
$$

\vspace{3mm}
\noindent
where constant $\tilde C_1$ does not depend on $p_1,$ 
constant $\tilde C_2$ does not depend on $p_2,$ etc.

Further, the construction of 
integrable majorants
when applying Lebesgue's 
Dominated Convergence Theorem 
in (\ref{udar1}) is obvious.

For example, to pass to the limit
$\varlimsup\limits_{p_k\to\infty},$
the integrable majorant has the form
(it is constructed on the base of (\ref{strange1}))

\vspace{-1mm}
$$
\biggl(R_{p_1\ldots p_k}(t_1,\ldots,t_k)\biggr)^{2n}\le 
$$

\vspace{-1mm}
$$
\le\Biggl(\left(\tilde C_1 + C_1\right)+\Biggr.
$$

$$
+
\sum\limits_{j_1=0}^{p_1}\left(\tilde C_2 + C_2\right)
\left|\phi_{j_1}(t_1)\right|+
\ldots\ 
$$

\vspace{-1mm}
$$
\Biggl.
\ldots\ +
\sum\limits_{j_1=0}^{p_1}\ldots 
\sum\limits_{j_{k-1}=0}^{p_{k-1}}\left(\tilde C_k + C_k\right)
\left|\phi_{j_{k-1}}(t_{k-1})\ldots \phi_{j_1}(t_1)\right|\Biggr)^{2n}\le
$$

$$
\le \Biggl(\left(\tilde C_1 + C_1\right)+\Biggr.
$$

$$
+
\sqrt{\frac{2}{T-t}}\ (p_1+1)\left(\tilde C_2 + C_2\right)+
\ldots\
$$

\vspace{-1mm}
\begin{equation}
\label{strange6}
\Biggl.\ldots\ +\Biggl(\sqrt{\frac{2}{T-t}}\Biggr)^{k-1}
(p_1+1)\ldots (p_{k-1}+1)
\left(\tilde C_k + C_k\right)\Biggr)^{2n},
\end{equation}

\vspace{5mm}
\noindent
where $n\in\mathbb{N},$ the numbers $p_1,\ldots,p_{k-1}$ are fixed and the right-hand side of
(\ref{strange6}) is independent of $p_k.$

Theorem 1 is proved.

It easy to notice that if we expand the function
$K^{*}(t_1,\ldots,t_k)$ into the generalized 
Fourier series at the interval $(t, T)$
at first with respect to the variable $t_k$, after that
with respect to the variable $t_{k-1}$, etc., then
we will have the expansion

\begin{equation}
\label{otit3333}
K^{*}(t_1,\ldots,t_k)=\sum_{j_k=0}^{\infty}\ldots \sum_{j_1=0}^{\infty}
C_{j_k\ldots j_1}\prod_{l=1}^{k} \phi_{j_l}(t_l)
\end{equation}

\vspace{4mm}
\noindent
instead of the expansion (\ref{30.18}).

Let us prove the expansion (\ref{otit3333}). Similarly 
with (\ref{oop1}) we have the following equality

\begin{equation}
\label{oop1otit}
\psi_k(t_k)\left({\bf 1}_{\{t_{k-1}<t_k\}}+
\frac{1}{2}{\bf 1}_{\{t_{k-1}=t_k\}}\right)=
\sum\limits_{j_k=0}^{\infty}\ \int\limits_{t_{k-1}}^{T}\psi_k(t_k)
\phi_{j_k}(t_k)dt_k\phi_{j_k}(t_k),
\end{equation}

\vspace{4mm}
\noindent
which is 
fulfilled pointwise at the interval $(t, T),$
besides
the series on the right-hand side
of (\ref{oop1otit}) converges when $t_1=t$ and $t_1=T.$ 

Let us introduce the assumption of induction

\vspace{1mm}
$$
\sum\limits_{j_k=0}^{\infty}\ldots
\sum\limits_{j_{3}=0}^{\infty}\psi_{2}(t_{2})
\int\limits_{t_2}^T\psi_{3}(t_{3})\phi_{j_{3}}(t_{3})\ldots
\int\limits_{t_{k-1}}^T
\psi_{k}(t_{k})\phi_{j_{k}}(t_{k})dt_k\ldots dt_{3}
\prod_{l=3}^{k}\phi_{j_{l}}(t_{l})=
$$

\vspace{1mm}
\begin{equation}
\label{oop22otit}
=\prod\limits_{l=2}^{k}\psi_l(t_l)
\prod_{l=2}^{k-1}\left({\bf 1}_{\{t_l<t_{l+1}\}}+
\frac{1}{2}{\bf 1}_{\{t_l=t_{l+1}\}}\right).
\end{equation}

\vspace{4mm}

Then

\vspace{1mm}
$$
\sum\limits_{j_k=0}^{\infty}\ldots
\sum\limits_{j_{3}=0}^{\infty}\sum\limits_{j_{2}=0}^{\infty}
\psi_{1}(t_{1})
\int\limits_{t_1}^{T}\psi_{2}(t_{2})\phi_{j_{2}}(t_{2})\ldots
\int\limits_{t_{k-1}}^{T}\psi_{k}(t_{k})
\phi_{j_{k}}(t_{k})
dt_k\ldots dt_{2}\prod_{l=2}^{k}\phi_{j_{l}}(t_{l})=
$$

\vspace{3mm}
$$
=\sum\limits_{j_k=0}^{\infty}\ldots
\sum\limits_{j_{3}=0}^{\infty}
\psi_1(t_1)\left({\bf 1}_{\{t_{1}<t_{2}\}}+
\frac{1}{2}{\bf 1}_{\{t_{1}=t_{2}\}}\right)\psi_{2}(t_{2})\times
$$

\vspace{3mm}
$$
\times
\int\limits_{t_2}^{T}\psi_{3}(t_{3})\phi_{j_{3}}(t_{3})\ldots
\int\limits_{t_{k-1}}^T\psi_{k}(t_{k})
\phi_{j_{k}}(t_{k})dt_k\ldots dt_{3}
\prod_{l=3}^{k}\phi_{j_{l}}(t_{l})=
$$

\vspace{3mm}
$$
=\psi_1(t_1)\left({\bf 1}_{\{t_{1}<t_{2}\}}+
\frac{1}{2}{\bf 1}_{\{t_{1}=t_{2}\}}\right)
\sum\limits_{j_k=0}^{\infty}\ldots
\sum\limits_{j_{3}=0}^{\infty}
\psi_{2}(t_{2})\times
$$

\vspace{3mm}
$$
\times
\int\limits_{t_2}^T\psi_{3}(t_{3})\phi_{j_{3}}(t_{3})\ldots
\int\limits_{t_{k-1}}^T\psi_{k}(t_{k})
\phi_{j_{k}}(t_{k})dt_k\ldots dt_{3}
\prod_{l=3}^{k}\phi_{j_{l}}(t_{l})=
$$

\vspace{3mm}
$$
=\psi_1(t_1)\left({\bf 1}_{\{t_{1}<t_{2}\}}+
\frac{1}{2}{\bf 1}_{\{t_{1}=t_{2}\}}\right)
\prod\limits_{l=2}^{k}\psi_l(t_l)
\prod_{l=2}^{k-1}\left({\bf 1}_{\{t_l<t_{l+1}\}}+
\frac{1}{2}{\bf 1}_{\{t_l=t_{l+1}\}}\right)=
$$

\vspace{3mm}
\begin{equation}
\label{oop30otit}
=\prod\limits_{l=1}^{k}\psi_l(t_l)
\prod_{l=1}^{k-1}\left({\bf 1}_{\{t_l<t_{l+1}\}}+
\frac{1}{2}{\bf 1}_{\{t_l=t_{l+1}\}}\right).
\end{equation}

\vspace{7mm}

From the other hand, the left-hand side
of (\ref{oop30otit}) can be represented 
in the following form

\vspace{1mm}
$$
\sum_{j_k=0}^{\infty}\ldots \sum_{j_1=0}^{\infty}
C_{j_k\ldots j_1}\prod_{l=1}^{k} \phi_{j_l}(t_l)
$$

\vspace{5mm}
\noindent
by 
expanding the function

\vspace{1mm}
$$
\psi_{1}(t_{1})
\int\limits_{t_1}^{T}\psi_{2}(t_{2})\phi_{j_{2}}(t_{2})\ldots
\int\limits_{t_{k-1}}^{T}\psi_{k}(t_{k})
\phi_{j_{k}}(t_{k})
dt_k\ldots dt_{2}
$$

\vspace{5mm}
\noindent
into the generalized Fourier series at the interval $(t, T)$ 
using the variable 
$t_1.$
Here we applied the following replacement of 
integration order

\vspace{2mm}
$$
\int\limits_t^T\psi_{1}(t_{1})
\int\limits_{t_1}^{T}\psi_{2}(t_{2})\phi_{j_{2}}(t_{2})\ldots
\int\limits_{t_{k-1}}^{T}\psi_{k}(t_{k})
\phi_{j_{k}}(t_{k})
dt_k\ldots dt_2dt_{1}=
$$

\vspace{2mm}
$$
=\int\limits_t^T\psi_{k}(t_{k})\phi_{j_{k}}(t_{k})
\ldots
\int\limits_t^{t_{3}}\psi_{2}(t_{2})
\phi_{j_{2}}(t_{2})\int\limits_t^{t_{2}}\psi_{1}(t_{1})
\phi_{j_{1}}(t_{1})
dt_1dt_2\ldots dt_k=
$$

\vspace{2mm}
$$
=C_{j_k\ldots j_1}.
$$

\vspace{6mm}

The expansion (\ref{otit3333}) is proved. So, we can formulate the 
following theorem.

\vspace{2mm}
         
{\bf Theorem 2} \cite{10a} (Sect.~2.4) (also see \cite{1997} (1997), 
\cite{2013}, \cite{2017a}, \cite{2017}, \cite{10}, \cite{10aaaa}-\cite{12aa-afterxxx}.\ 
{\it Suppose that the conditions of 
Theorem {\rm 1} are fulfilled.
Then

\begin{equation}
\label{aq1}
J^{*}[\psi^{(k)}]_{T,t}=
\sum_{j_k=0}^{\infty}\ldots\sum_{j_1=0}^{\infty}
C_{j_k\ldots j_1}
\prod_{l=1}^k
\zeta^{(i_l)}_{j_l},
\end{equation}

\vspace{5mm}
\noindent
where notations are the same as in Theorem {\rm 1}.}

\vspace{2mm}

Note that (\ref{aq1})  means the following 

\vspace{2mm}
$$
\lim\limits_{p_k\to\infty}
\varlimsup\limits_{p_{k-1}\to\infty}\ldots
\varlimsup\limits_{p_1\to\infty}
{\sf M}\left\{\left(J^{*}[\psi^{(k)}]_{T,t}-
\sum_{j_k=0}^{p_k}\ldots\sum_{j_1=0}^{p_1}
C_{j_k\ldots j_1}
\prod_{l=1}^k
\zeta^{(i_l)}_{j_l}\right)^{2n}\right\}=0, 
$$

\vspace{5mm}
\noindent
where
$n\in\mathbb{N}.$

Let us make a remark about how one can obtain
an analogue of Theorem~1 for the 
complete orthonormal system of Legendre 
polynomials in the space $L_2([t, T])$
and $n=1$ (the case of mean-square convergence), $k=2.$

First note the well known estimate for Legendre
polynomials \cite{Gob}

\begin{equation}
\label{otit987}
\left|P_j(y)\right| <\frac{K}{\sqrt{j+1}(1-y^2)^{1/4}},\ \ \ 
y\in (-1, 1),\ \ \ j\in \mathbb{N},
\end{equation}

\vspace{4mm}
\noindent
where constant $K$ does not depend on $y$ and $j$.

By analogy with (\ref{leto80010aaa}) we have

$$
{\sf M}\left\{\left(J[R_{p_1p_2}]_{T,t}^{(2)}\right)^{2}
\right\}\le
$$

\vspace{2mm}
\begin{equation}
\label{strange10}
\le
2\int\limits_{[t, T]^2}
\left(R_{p_1p_2}(t_1,t_2)\right)^{2}dt_1 dt_2
+
{\bf 1}_{\{i_1=i_2\ne 0\}}
\left(\int\limits_t^T R_{p_1p_2}(t_1,t_1)dt_1\right)^2.
\end{equation}

\vspace{5mm}

From Remark~1.6, Sect.~1.7.2 \cite{10a} and (1.72), (2.103) \cite{10a} we obtain
for the case of Legendre polynomials

\vspace{-1mm}
$$
\lim\limits_{p_1\to\infty}\varlimsup\limits_{p_2\to\infty}
\int\limits_{[t, T]^2}
\left(R_{p_1p_2}(t_1,t_2)\right)^{2}dt_1 dt_2=0.
$$

\vspace{4mm}

Further, we have (see (\ref{d2020}))

\vspace{-1mm}
$$
R_{p_1p_2}(t_1,t_1)=\left(K^{*}(t_1,t_1)-\sum\limits_{j_1=0}^{p_1}
C_{j_1}(t_1)\phi_{j_1}(t_1)\right)+
$$

\begin{equation}
\label{strange14}
+\left(
\sum\limits_{j_1=0}^{p_1}\left(C_{j_1}(t_1)-
\sum\limits_{j_2=0}^{p_2}
C_{j_2j_1}\phi_{j_2}(t_1)\right)
\phi_{j_1}(t_1)\right).
\end{equation}

\vspace{5mm}

Then, taking into account (\ref{strange101}),
(\ref{strange14}) and applying two times (we mean here an iterated passage to the limit
$\lim\limits_{p_1\to\infty}\varlimsup\limits_{p_2\to\infty}$)
the Lebesgue's 
Dominated Convergence Theorem,
we obtain

\vspace{-1mm}
$$
\lim\limits_{p_1\to\infty}\varlimsup\limits_{p_2\to\infty}
\int\limits_t^T
R_{p_1p_2}(t_1,t_1)dt_1=0.
$$

\vspace{4mm}

Let us discuss the choice of integrable majorants
when applying Lebesgue's 
Dominated Convergence Theorem in our case.

Using double integration by parts (as in (2.22), Sect.~2.1.1 \cite{10a}), we estimate
the partial sums of one-dimensional Fourier--Legendre series

$$
\sum\limits_{j_1=0}^{p_1}
C_{j_1}(t_1)\phi_{j_1}(t_1),\ \ \ 
\sum\limits_{j_2=0}^{p_2}
C_{j_2j_1}\phi_{j_2}(t_1)
$$ 

\vspace{4mm}
\noindent
in (\ref{strange14}) using (\ref{otit987}) and (\ref{strange3}) as follows

\begin{equation}
\label{strange17}
\left|\sum\limits_{j_1=0}^{p_1}
C_{j_1}(t_1)\phi_{j_1}(t_1)\right|
\le 
K_1\Biggl(1+
\frac{1}{\left(1-(z(t_1))^2\right)^{1/2}}+
\frac{1}{\left(1-(z(t_1))^2\right)^{1/4}}\Biggr),
\end{equation}

\begin{equation}
\label{strange18}
\left|\sum\limits_{j_2=0}^{p_2}
C_{j_2j_1}\phi_{j_2}(t_1)
\right|
\le K_2\left(1+\frac{1}{\left(1-(z(t_1))^2\right)^{1/4}}\right),
\end{equation}

\vspace{4mm}
\noindent
where 
$$
z(s)=\left(s-\frac{T+t}{2}\right)\frac{2}{T-t},
$$

\vspace{3mm}
\noindent
constant $K_1$ does not depend on $p_1,$ and constant $K_2$ does not depend on $p_2.$ 

Thus, integrable majorants
in our case can be easily constracted 
using (\ref{strange14}), (\ref{strange17}) and (\ref{strange18}) (see the proof of Theorem~1
for details).

An analogue of Theorem~1 for the 
case of Legendre 
polynomials 
and $n=1$ (the case of mean-square convergence), $k=2$ is obtained.

\vspace{5mm}

\section{Examples. The Case of Legendre Polynomials}

\vspace{5mm}

In this section, we provide some practical material (based on 
an analogue of Theorem~1 for the 
case of Legendre 
polynomials and $k=2,$ $n=1$) on expansions 
of iterated
Stra\-to\-no\-vich stochastic integrals of the following form \cite{10a}-\cite{12aa-afterxxx}

\vspace{-1mm}
\begin{equation}
\label{k1001}
I_{(l_1\ldots l_k)T,t}^{*(i_1\ldots i_k)}
={\int\limits_t^{*}}^T (t-t_k)^{l_k} \ldots {\int\limits_t^{*}}^{t_2}
(t-t_1)^{l_1} d{\bf f}_{t_1}^{(i_1)}\ldots
d{\bf f}_{t_k}^{(i_k)},
\end{equation}

\vspace{4mm}
\noindent
where $i_1,\ldots, i_k=1,\dots,m,$\ \  $l_1,\ldots,l_k=0, 1,\ldots$

The complete orthonormal system of Legendre polynomials in the 
space $L_2([t,T])$ looks as follows

\begin{equation}
\label{4009}
\phi_j(x)=\sqrt{\frac{2j+1}{T-t}}P_j\left(\left(
x-\frac{T+t}{2}\right)\frac{2}{T-t}\right),\ \ \ j=0, 1, 2,\ldots,
\end{equation}

\vspace{3mm}
\noindent
where $P_j(x)$ is the Legendre polynomial. 

Using 
an analogue of Theorem~1 for the 
system of 
functions (\ref{4009}) and $k=2,$ $n=1$, we obtain the following expansions of iterated
Stratonovich stochastic integrals \cite{1997}--\cite{12aa-afterxxx},
\cite{13}, \cite{15}, \cite{15b}, \cite{15d}-\cite{art-zero}

\vspace{2mm}
$$
I_{(0)T,t}^{*(i_1)}=\sqrt{T-t}\zeta_0^{(i_1)},
$$

\vspace{3mm}

\begin{equation}
\label{yyyy1}
I_{(1)T,t}^{*(i_1)}=-\frac{(T-t)^{3/2}}{2}\left(\zeta_0^{(i_1)}+
\frac{1}{\sqrt{3}}\zeta_1^{(i_1)}\right),
\end{equation}

\vspace{3mm}

\begin{equation}
\label{yyyy2}
I_{(2)T,t}^{*(i_1)}=\frac{(T-t)^{5/2}}{3}\left(\zeta_0^{(i_1)}+
\frac{\sqrt{3}}{2}\zeta_1^{(i_1)}+
\frac{1}{2\sqrt{5}}\zeta_2^{(i_1)}\right),
\end{equation}

\vspace{4mm}

\begin{equation}
\label{yyyy3}
I_{(00)T,t}^{*(i_1 i_2)}=
\frac{T-t}{2}\left(\zeta_0^{(i_1)}\zeta_0^{(i_2)}+\sum_{i=1}^{\infty}
\frac{1}{\sqrt{4i^2-1}}\biggl(
\zeta_{i-1}^{(i_1)}\zeta_{i}^{(i_2)}-
\zeta_i^{(i_1)}\zeta_{i-1}^{(i_2)}\biggr)\right),
\end{equation}

\vspace{6mm}

$$
I_{(01)T,t}^{*(i_1 i_2)}=-\frac{T-t}{2}I_{(00)T,t}^{*(i_1 i_2)}
-\frac{(T-t)^2}{4}\left(
\frac{\zeta_0^{(i_1)}\zeta_1^{(i_2)}}{\sqrt{3}}+\right.
$$

\vspace{2mm}
$$
+\left.\sum_{i=0}^{\infty}\left(
\frac{(i+2)\zeta_i^{(i_1)}\zeta_{i+2}^{(i_2)}
-(i+1)\zeta_{i+2}^{(i_1)}\zeta_{i}^{(i_2)}}
{\sqrt{(2i+1)(2i+5)}(2i+3)}-
\frac{\zeta_i^{(i_1)}\zeta_{i}^{(i_2)}}{(2i-1)(2i+3)}\right)\right),
$$

\vspace{6mm}

$$
I_{(10)T,t}^{*(i_1 i_2)}=-\frac{T-t}{2}I_{(00)T,t}^{*(i_1 i_2)}
-\frac{(T-t)^2}{4}\left(
\frac{\zeta_0^{(i_2)}\zeta_1^{(i_1)}}{\sqrt{3}}+\right.
$$

\vspace{2mm}
$$
+\left.\sum_{i=0}^{\infty}\left(
\frac{(i+1)\zeta_{i+2}^{(i_2)}\zeta_{i}^{(i_1)}
-(i+2)\zeta_{i}^{(i_2)}\zeta_{i+2}^{(i_1)}}
{\sqrt{(2i+1)(2i+5)}(2i+3)}+
\frac{\zeta_i^{(i_1)}\zeta_{i}^{(i_2)}}{(2i-1)(2i+3)}\right)\right),
$$

\vspace{9mm}

$$
I_{(02)T,t}^{*(i_1 i_2)}=-\frac{(T-t)^2}{4}I_{(00)T,t}^{*(i_1 i_2)}
-(T-t)I_{(01)T,t}^{*(i_1 i_2)}+
\frac{(T-t)^3}{8}\left(
\frac{2\zeta_2^{(i_2)}\zeta_0^{(i_1)}}{3\sqrt{5}}+\right.
$$

\vspace{2mm}
$$
+\frac{1}{3}\zeta_0^{(i_1)}\zeta_0^{(i_2)}+
\sum_{i=0}^{\infty}\left(
\frac{(i+2)(i+3)\zeta_{i+3}^{(i_2)}\zeta_{i}^{(i_1)}
-(i+1)(i+2)\zeta_{i}^{(i_2)}\zeta_{i+3}^{(i_1)}}
{\sqrt{(2i+1)(2i+7)}(2i+3)(2i+5)}+
\right.
$$

\vspace{2mm}
$$
\left.\left.+
\frac{(i^2+i-3)\zeta_{i+1}^{(i_2)}\zeta_{i}^{(i_1)}
-(i^2+3i-1)\zeta_{i}^{(i_2)}\zeta_{i+1}^{(i_1)}}
{\sqrt{(2i+1)(2i+3)}(2i-1)(2i+5)}\right)\right),
$$

\vspace{9mm}

$$
I_{(20)T,t}^{*(i_1 i_2)}=-\frac{(T-t)^2}{4}I_{(00)T,t}^{*(i_1 i_2)}
-(T-t)I_{(10)T,t}^{*(i_1 i_2)}+
\frac{(T-t)^3}{8}\left(
\frac{2\zeta_0^{(i_2)}\zeta_2^{(i_1)}}{3\sqrt{5}}
+\right.
$$

\vspace{2mm}
$$
+\frac{1}{3}\zeta_0^{(i_1)}\zeta_0^{(i_2)}+
\sum_{i=0}^{\infty}\left(
\frac{(i+1)(i+2)\zeta_{i+3}^{(i_2)}\zeta_{i}^{(i_1)}
-(i+2)(i+3)\zeta_{i}^{(i_2)}\zeta_{i+3}^{(i_1)}}
{\sqrt{(2i+1)(2i+7)}(2i+3)(2i+5)}+
\right.
$$

\vspace{2mm}
$$
\left.\left.+
\frac{(i^2+3i-1)\zeta_{i+1}^{(i_2)}\zeta_{i}^{(i_1)}
-(i^2+i-3)\zeta_{i}^{(i_2)}\zeta_{i+1}^{(i_1)}}
{\sqrt{(2i+1)(2i+3)}(2i-1)(2i+5)}\right)\right),
$$

\vspace{9mm}

$$
I_{(11)T,t}^{*(i_1 i_2)}=-\frac{(T-t)^2}{4}I_{(00)T,t}^{*(i_1 i_2)}
-\frac{(T-t)}{2}\left(I_{(10)T,t}^{*(i_1 i_2)}+
I_{(01)T,t}^{*(i_1 i_2)}\right)+
$$

\vspace{2mm}
$$
+
\frac{(T-t)^3}{8}\left(
\frac{1}{3}\zeta_1^{(i_1)}\zeta_1^{(i_2)}
+
\sum_{i=0}^{\infty}\left(
\frac{(i+1)(i+3)\left(\zeta_{i+3}^{(i_2)}\zeta_{i}^{(i_1)}
-\zeta_{i}^{(i_2)}\zeta_{i+3}^{(i_1)}\right)}
{\sqrt{(2i+1)(2i+7)}(2i+3)(2i+5)}+
\right.\right.
$$

\vspace{2mm}
$$
\left.\left.+
\frac{(i+1)^2\left(\zeta_{i+1}^{(i_2)}\zeta_{i}^{(i_1)}
-\zeta_{i}^{(i_2)}\zeta_{i+1}^{(i_1)}\right)}
{\sqrt{(2i+1)(2i+3)}(2i-1)(2i+5)}\right)\right),
$$

\vspace{6mm}

$$
I_{(3)T,t}^{*(i_1)}=-\frac{(T-t)^{7/2}}{4}\left(\zeta_0^{(i_1)}+
\frac{3\sqrt{3}}{5}\zeta_1^{(i_1)}+
\frac{1}{\sqrt{5}}\zeta_2^{(i_1)}+
\frac{1}{5\sqrt{7}}\zeta_3^{(i_1)}\right),
$$

\vspace{6mm}
\noindent
where 
\begin{equation}
\label{cx1}
\zeta_{j}^{(i)}=
\int\limits_t^T \phi_{j}(s) d{\bf f}_s^{(i)}
\end{equation}

\vspace{5mm}
\noindent
are independent standard Gaussian random variables
for 
various
$i$ or $j$ $(i=1,\ldots,m)$.

\vspace{5mm}

\section{Examples. The Case of Trigonometric Functions}

\vspace{5mm}

Let us consider the Milstein expansions of the 
integrals $I_{(1)T,t}^{(i_1)},$
$I_{(00)T,t}^{*(i_1 i_2)},$ $I_{(2)T,t}^{*(i_1)}$ (see \cite{1}-\cite{3})
based on the trigonometric Fourier expansion of the Brownian
Bridge process
(the version of the so-called Karhunen--Loeve expansion)

\begin{equation}
\label{xxxx1}
I_{(1)T,t}^{*(i_1)}=-\frac{{(T-t)}^{3/2}}{2}
\left(\zeta_0^{(i_1)}-\frac{\sqrt{2}}{\pi}\sum_{r=1}^{\infty}
\frac{1}{r}
\zeta_{2r-1}^{(i_1)}
\right),
\end{equation}

\vspace{1mm}

\begin{equation}
\label{xxxx2}
I_{(2)T,t}^{*(i_1)}=
(T-t)^{5/2}\left(
\frac{1}{3}\zeta_0^{(i_1)}+\frac{1}{\sqrt{2}\pi^2}
\sum_{r=1}^{\infty}\frac{1}{r^2}\zeta_{2r}^{(i_1)}
-\frac{1}{\sqrt{2}\pi}\sum_{r=1}^{\infty}
\frac{1}{r}\zeta_{2r-1}^{(i_1)}\right),
\end{equation}

\vspace{5mm}

$$
I_{(00)T,t}^{*(i_1 i_2)}=\frac{1}{2}(T-t)\Biggl(
\zeta_{0}^{(i_1)}\zeta_{0}^{(i_2)}\Biggr.
+
$$

\vspace{2mm}

\begin{equation}
\label{xxxx3}
\Biggl.+\frac{1}{\pi}
\sum_{r=1}^{\infty}\frac{1}{r}\biggl(
\zeta_{2r}^{(i_1)}\zeta_{2r-1}^{(i_2)}-
\zeta_{2r-1}^{(i_1)}\zeta_{2r}^{(i_2)}
+\sqrt{2}\biggl(\zeta_{2r-1}^{(i_1)}\zeta_{0}^{(i_2)}-
\zeta_{0}^{(i_1)}\zeta_{2r-1}^{(i_2)}\biggr)\biggr)\Biggr),
\end{equation}

\vspace{6mm}
\noindent
where $\zeta_0^{(i)},$ $\zeta_{2r}^{(i)},$
$\zeta_{2r-1}^{(i)}$
($i=1,\ldots,m$) are independent
standard Gaussian random variables defined by the relation (\ref{cx1})
in which $\left\{\phi_j(x)\right\}_{j=0}^{\infty}$
is a complete orthonornal system of trigonometric
functions in $L_2([t, T]).$

It is obviously that at least (\ref{xxxx1})--(\ref{xxxx3}) 
are significantly more complicated 
in comparison with (\ref{yyyy1})--(\ref{yyyy3}).
Note that (\ref{xxxx1})--(\ref{xxxx3}) also can be 
obtained using Theorem 1 \cite{1997}, \cite{1a},
\cite{3a}-\cite{2017a}, \cite{2017}-\cite{art-zero}.

\vspace{5mm}

\section{Further Remarks}

\vspace{5mm}

In this section, we consider some approaches on the base of Theorem 1 
for the case $k=2.$ Moreover,
we explain the potential 
difficulties associated with the use 
of generalized multiple Fourier series converging 
almost everywhere 
in the hypercube 
$[t, T]^k$ in the proof of Theorem 1.

First, we show how iterated series can be replaced by 
multiple one in Theorem 1 (the case $k=2$ and $n=1$) 
and in analogue of Theorem~1 for the 
case of Legendre 
polynomials 
(the case $k=2$ and $n=1$).

We have

$$
\lim\limits_{p\to\infty}
{\sf M}\left\{\left(
J^{*}[\psi^{(2)}]_{T,t}-
\sum\limits_{j_1=0}^{p}\sum\limits_{j_2=0}^{p}
C_{j_2 j_1}\zeta_{j_1}^{(i_1)}
\zeta_{j_2}^{(i_2)}\right)^{2}\right\}=
$$

\vspace{2mm}
$$
=\lim\limits_{p\to\infty}\varlimsup\limits_{q\to\infty}
{\sf M}\left\{\left(
J^{*}[\psi^{(2)}]_{T,t}-
\sum\limits_{j_1=0}^{p}\sum\limits_{j_2=0}^{p}
C_{j_2 j_1}\zeta_{j_1}^{(i_1)}
\zeta_{j_2}^{(i_2)}\right)^{2}\right\}\le
$$

\vspace{2mm}
$$
\le
\lim\limits_{p\to\infty}\varlimsup\limits_{q\to\infty}\left(
2{\sf M}\left\{\left(
J^{*}[\psi^{(2)}]_{T,t}-
\sum\limits_{j_1=0}^{p}\sum\limits_{j_2=0}^{q}
C_{j_2 j_1}\zeta_{j_1}^{(i_1)}
\zeta_{j_2}^{(i_2)}\right)^{2}\right\}+\right.
$$

\vspace{2mm}
$$
\left.+
2{\sf M}\left\{\left(
\sum\limits_{j_1=0}^{p}\sum\limits_{j_2=0}^{q}
C_{j_2 j_1}\zeta_{j_1}^{(i_1)}
\zeta_{j_2}^{(i_2)}
-
\sum\limits_{j_1=0}^{p}\sum\limits_{j_2=0}^{p}
C_{j_2 j_1}\zeta_{j_1}^{(i_1)}
\zeta_{j_2}^{(i_2)}\right)^{2}\right\}\right)=
$$

\vspace{2mm}
$$
=\lim\limits_{p\to\infty}\varlimsup\limits_{q\to\infty}
2{\sf M}\left\{\left(
\sum\limits_{j_1=0}^{p}\sum\limits_{j_2=p+1}^{q}
C_{j_2 j_1}\zeta_{j_1}^{(i_1)}
\zeta_{j_2}^{(i_2)}\right)^{2}\right\}=
$$

\vspace{3mm}

$$
=\lim\limits_{p\to\infty}\varlimsup\limits_{q\to\infty}
2
\sum\limits_{j_1=0}^{p}\sum\limits_{j_1'=0}^{p}
\sum\limits_{j_2=p+1}^{q}\sum\limits_{j_2'=p+1}^{q}
C_{j_2 j_1}C_{j_2' j_1'}{\sf M}\left\{\zeta_{j_1}^{(i_1)}
\zeta_{j_1'}^{(i_1)}\right\}
{\sf M}\left\{\zeta_{j_2}^{(i_2)}\zeta_{j_2'}^{(i_2)}\right\}=
$$

\vspace{3mm}

$$
=2\lim\limits_{p\to\infty}\lim\limits_{q\to\infty}
\sum\limits_{j_1=0}^{p}\sum\limits_{j_2=p+1}^{q}
C_{j_2 j_1}^2=
$$

\vspace{3mm}

\begin{equation}
\label{fd1}
=2\lim\limits_{p\to\infty}\lim\limits_{q\to\infty}
\left(\sum\limits_{j_1=0}^{p}\sum\limits_{j_2=0}^{q}
C_{j_2 j_1}^2-
\sum\limits_{j_1=0}^{p}\sum\limits_{j_2=0}^{p}
C_{j_2 j_1}^2\right)
=
\end{equation}

\vspace{2mm}

\begin{equation}
\label{fd2}
=2\left(\lim\limits_{p,q\to\infty}
\sum\limits_{j_1=0}^{p}\sum\limits_{j_2=0}^{q}
C_{j_2 j_1}^2-
\lim\limits_{p\to\infty}
\sum\limits_{j_1=0}^{p}\sum\limits_{j_2=0}^{p}
C_{j_2 j_1}^2\right)
=
\end{equation}

\vspace{2mm}
\begin{equation}
\label{fd3}
=2\int\limits_{[t,T]^2}K^2 (t_1,t_2)dt_1dt_2-
2\int\limits_{[t,T]^2}K^2 (t_1,t_2)dt_1dt_2=0,
\end{equation}

\vspace{4.5mm}
\noindent 
where the function $K(t_1,t_2)$ is defined by (\ref{ppp}) for $k=2.$

Note that the transition from (\ref{fd1}) to (\ref{fd2})
is based on the theorem on reducing 
of a limit to iterated one.
Moreover, the transition from (\ref{fd2}) to (\ref{fd3})
is based on the Parseval equality.

Thus, we obtain the following Theorem.

\vspace{2mm}

{\bf Theorem 3.}\ {\it Assume that
$\{\phi_j(x)\}_{j=0}^{\infty}$ is a complete orthonormal
system of Legendre polynomials or trigonometric functions
in the space $L_2([t, T])$. 
At the same time $\psi_2(\tau)$ is a continuously dif\-ferentiable 
nonrandom function on $[t, T]$ and $\psi_1(\tau)$ is twice 
continuously differentiable nonrandom function on $[t, T]$. 
Then$,$ for the iterated Stratonovich stochastic integral {\rm (\ref{str})}
of multiplicity $2$

\vspace{1mm}
$$
J^{*}[\psi^{(2)}]_{T,t}=
{\int\limits_t^{*}}^T\psi_2(t_2)
{\int\limits_t^{*}}^{t_2}
\psi_1(t_1) d{\bf w}_{t_1}^{(i_1)}
d{\bf w}_{t_2}^{(i_2)}\ \ \ 
(i_1,i_2=0, 1,\ldots,m)
$$

\vspace{5mm}
\noindent
the following 
expansion 

\vspace{1mm}
$$
J^{*}[\psi^{(2)}]_{T,t}=
\hbox{\vtop{\offinterlineskip\halign{
\hfil#\hfil\cr
{\rm l.i.m.}\cr
$\stackrel{}{{}_{p\to \infty}}$\cr
}} }
\sum\limits_{j_1,j_2=0}^{p}
C_{j_2 j_1}\zeta_{j_1}^{(i_1)}
\zeta_{j_2}^{(i_2)}
$$

\vspace{4mm}
\noindent
converging in the mean-square sense is valid, where the Fourier coefficient 
$C_{j_2 j_1}$ has the form

$$
C_{j_2 j_1}=\int\limits_t^T\psi_2(t_2)\phi_{j_2}(t_2)
\int\limits_t^{t_2}
\psi_1(t_1)\phi_{j_1}(t_1)
dt_1dt_2
$$

\vspace{2mm}
\noindent
and
$$
\zeta_{j}^{(i)}=
\int\limits_t^T \phi_{j}(s) d{\bf w}_s^{(i)}
$$ 

\vspace{4mm}
\noindent
are independent standard Gaussian random variables for various 
$i$ or $j$ {\rm (}if $i\ne 0${\rm )},
${\bf w}_{\tau}^{(i)}={\bf f}_{\tau}^{(i)}$ are independent 
standard Wiener processes
$(i=1,\ldots,m)$ and 
${\bf w}_{\tau}^{(0)}=\tau.$}

\vspace{2mm}

Note that Theorem 3 is a modification (for the 
case $p_1=p_2=p$ of series summation) of Theorem 2.1 \cite{10a}.

From the other hand, Theorem 1 implies the following

\vspace{2mm}
$$
0\le
\left\vert \lim\limits_{p_1\to\infty}
\varlimsup\limits_{p_2\to\infty}\ldots 
\varlimsup\limits_{p_k\to\infty}{\sf M}\left\{
\sum_{j_1=0}^{p_1}\ldots\sum_{j_k=0}^{p_k}
C_{j_k\ldots j_1}
\prod_{l=1}^k
\zeta^{(i_l)}_{j_l}-J^{*}[\psi^{(k)}]_{T,t}\right\}\right\vert \le
$$

\vspace{4mm}
$$
\le \lim\limits_{p_1\to\infty}
\varlimsup\limits_{p_2\to\infty}\ldots 
\varlimsup\limits_{p_k\to\infty} \left\vert {\sf M}\left\{
\sum_{j_1=0}^{p_1}\ldots\sum_{j_k=0}^{p_k}
C_{j_k\ldots j_1}
\prod_{l=1}^k
\zeta^{(i_l)}_{j_l}-J^{*}[\psi^{(k)}]_{T,t}\right\}\right\vert \le
$$

\vspace{4mm}
$$
\le \lim\limits_{p_1\to\infty}
\varlimsup\limits_{p_2\to\infty}\ldots 
\varlimsup\limits_{p_k\to\infty} {\sf M}\left\{\left\vert 
J^{*}[\psi^{(k)}]_{T,t}-\sum_{j_1=0}^{p_1}\ldots\sum_{j_k=0}^{p_k}
C_{j_k\ldots j_1}
\prod_{l=1}^k
\zeta^{(i_l)}_{j_l}\right\vert\right\} \le
$$

\vspace{4mm}
\begin{equation}
\label{het}
\le 
\lim\limits_{p_1\to\infty}
\varlimsup\limits_{p_2\to\infty}\ldots 
\varlimsup\limits_{p_k\to\infty}
\left(\hspace{-1mm}{\sf M}\hspace{-1mm}
\left\{\left(J^{*}[\psi^{(k)}]_{T,t}-
\sum_{j_1=0}^{p_1}\ldots\sum_{j_k=0}^{p_k}
C_{j_k\ldots j_1}
\prod_{l=1}^k
\zeta^{(i_l)}_{j_l}\right)^{2}\right\}\right)^{\hspace{-2mm}1/2}
\hspace{-1mm}=0.
\end{equation}

\vspace{6mm}

Moreover,

$$
\lim\limits_{p_1\to\infty}
\varlimsup\limits_{p_2\to\infty}\ldots 
\varlimsup\limits_{p_k\to\infty} 
\left(
\sum_{j_1=0}^{p_1}\ldots\sum_{j_k=0}^{p_k}
C_{j_k\ldots j_1}
{\sf M}\left\{\prod_{l=1}^k
\zeta^{(i_l)}_{j_l}\right\}-
{\sf M}\left\{J^{*}[\psi^{(k)}]_{T,t}\right\}\right)=
$$

\vspace{4mm}
\begin{equation}
\label{het100}
=
\lim\limits_{p_1\to\infty}
\varlimsup\limits_{p_2\to\infty}\ldots 
\varlimsup\limits_{p_k\to\infty}
\sum_{j_1=0}^{p_1}\ldots\sum_{j_k=0}^{p_k}C_{j_k\ldots j_1}
{\sf M}\left\{\prod_{l=1}^k
\zeta^{(i_l)}_{j_l}\right\}
-{\sf M}\left\{
J^{*}[\psi^{(k)}]_{T,t}\right\}.
\end{equation}

\vspace{7mm}

Combining (\ref{het}) and (\ref{het100}), we obtain

\vspace{1mm}
\begin{equation}
\label{het1}
{\sf M}\left\{
J^{*}[\psi^{(k)}]_{T,t}\right\}=
\lim\limits_{p_1\to\infty}
\varlimsup\limits_{p_2\to\infty}\ldots 
\varlimsup\limits_{p_k\to\infty}
\sum_{j_1=0}^{p_1}\ldots\sum_{j_k=0}^{p_k}
C_{j_k\ldots j_1}
{\sf M}\left\{\prod_{l=1}^k
\zeta^{(i_l)}_{j_l}\right\}.
\end{equation}

\vspace{6mm}

The relation (\ref{het1}) with $k=2$ implies the
following

$$
{\sf M}\left\{
J^{*}[\psi^{(2)}]_{T,t}\right\}= \frac{1}{2}{\bf 1}_{\{i_1=i_2\ne 0\}}
\int\limits_t^T
\psi_1(\tau)\psi_2(\tau)d\tau=
$$

\vspace{2mm}
\begin{equation}
\label{het2}
=\lim\limits_{p_1\to\infty}
\varlimsup\limits_{p_2\to\infty} 
\sum_{j_1=0}^{p_1}\sum_{j_2=0}^{p_2}
C_{j_2j_1}
{\sf M}\left\{
\zeta^{(i_1)}_{j_1}\zeta^{(i_2)}_{j_2}\right\},
\end{equation}

\vspace{6mm}
\noindent
where ${\bf 1}_A$ is the indicator of the set $A.$

Since

\vspace{-1mm}
$$
{\sf M}\left\{
\zeta^{(i_1)}_{j_1}\zeta^{(i_2)}_{j_2}\right\}=
{\bf 1}_{\{i_1=i_2\ne 0\}}{\bf 1}_{\{j_1=j_2\}},
$$

\vspace{5mm}
\noindent
then from (\ref{het2}) we obtain

\vspace{2mm}
$$
{\sf M}\left\{
J^{*}[\psi^{(2)}]_{T,t}\right\}=\lim\limits_{p_1\to\infty}
\varlimsup\limits_{p_2\to\infty} 
\sum_{j_1=0}^{p_1}\sum_{j_2=0}^{p_2}
C_{j_2j_1}{\bf 1}_{\{j_1=j_2\}}{\bf 1}_{\{i_1=i_2\ne 0\}}=
$$

\vspace{3mm}
\begin{equation}
\label{het5}
={\bf 1}_{\{i_1=i_2\ne 0\}}\lim\limits_{p_1\to\infty}
\varlimsup\limits_{p_2\to\infty} 
\sum_{j_1=0}^{{\rm min}\{p_1, p_2\}}
C_{j_1j_1}=
{\bf 1}_{\{i_1=i_2\ne 0\}}\sum_{j_1=0}^{\infty}
C_{j_1j_1},
\end{equation}

\vspace{6mm}
\noindent
where $C_{j_1j_1}$ is defined by (\ref{333.40}) for $k=2$ and $j_1=j_2,$
i.e. 

\vspace{2mm}
$$
C_{j_1j_1}=\int\limits_t^T \psi_2(t_2)\phi_{j_1}(t_2)
\int\limits_t^{t_2} \psi_1(t_1)\phi_{j_1}(t_1)dt_1 dt_2.
$$

\vspace{5mm}

From (\ref{het2}) and (\ref{het5}) we obtain the following relation

\begin{equation}
\label{het7}
\sum_{j_1=0}^{\infty}
C_{j_1j_1}=
\frac{1}{2}
\int\limits_t^T
\psi_1(\tau)\psi_2(\tau)d\tau.
\end{equation}

\vspace{4mm}

Note that the equality (\ref{het7}) and existence of the limit on the left-hand side of
(\ref{het7}) are proved in \cite{10a} (Sect.~2.1.2, 2.1.4), \cite{12} for the polynomial
and trigonometric cases 
($\psi_1(\tau), \psi_2(\tau)$ are continuously differentiable functions on $[t, T]$)
as well as for an arbitrary 
complete orthonormal system of functions in $L_2([t, T])$ 
and $\psi_1(\tau), \psi_2(\tau)\in L_2([t, T]).$

Let us address now to the following theorem on 
expansion of iterated Ito stochastic integrals (\ref{ito}).

\vspace{2mm}

{\bf Theorem 4} \cite{3a} (2006), \cite{3aa}-\cite{art-zero}. 
{\it Suppose that
every $\psi_l(\tau)$ $(l=1,\ldots, k)$ is a continuous 
nonrandom function on
the interval $[t, T]$ and
$\{\phi_j(x)\}_{j=0}^{\infty}$ is a complete orthonormal system  
of continuous func\-ti\-ons in the space $L_2([t,T]).$ 
Then

$$
J[\psi^{(k)}]_{T,t}\  =\ 
\hbox{\vtop{\offinterlineskip\halign{
\hfil#\hfil\cr
{\rm l.i.m.}\cr
$\stackrel{}{{}_{p_1,\ldots,p_k\to \infty}}$\cr
}} }\sum_{j_1=0}^{p_1}\ldots\sum_{j_k=0}^{p_k}
C_{j_k\ldots j_1}\Biggl(
\prod_{l=1}^k\zeta_{j_l}^{(i_l)}\ -
\Biggr.
$$

\vspace{2mm}
\begin{equation}
\label{tyyy}
-\ \Biggl.
\hbox{\vtop{\offinterlineskip\halign{
\hfil#\hfil\cr
{\rm l.i.m.}\cr
$\stackrel{}{{}_{N\to \infty}}$\cr
}} }\sum_{(l_1,\ldots,l_k)\in {\rm G}_k}
\phi_{j_{1}}(\tau_{l_1})
\Delta{\bf w}_{\tau_{l_1}}^{(i_1)}\ldots
\phi_{j_{k}}(\tau_{l_k})
\Delta{\bf w}_{\tau_{l_k}}^{(i_k)}\Biggr),
\end{equation}

\vspace{5mm}
\noindent
where

$$
{\rm G}_k={\rm H}_k\backslash{\rm L}_k,\ \ \
{\rm H}_k=\{(l_1,\ldots,l_k):\ l_1,\ldots,l_k=0,\ 1,\ldots,N-1\},
$$

\vspace{-2mm}
$$
{\rm L}_k=\{(l_1,\ldots,l_k):\ l_1,\ldots,l_k=0,\ 1,\ldots,N-1;\
l_g\ne l_r\ (g\ne r);\ g, r=1,\ldots,k\},
$$

\vspace{5mm}
\noindent
${\rm l.i.m.}$ is a limit in the mean-square sense,
$i_1,\ldots,i_k=0,1,\ldots,m,$ 

\vspace{-1mm}
\begin{equation}
\label{rr23}
\zeta_{j}^{(i)}=
\int\limits_t^T \phi_{j}(s) d{\bf w}_s^{(i)}
\end{equation} 

\vspace{2mm}
\noindent
are independent standard Gaussian random variables
for various
$i$ or $j$ {\rm(}if $i\ne 0${\rm),}
$C_{j_k\ldots j_1}$ is the Fourier coefficient {\rm(\ref{333.40}),}
$\Delta{\bf w}_{\tau_{j}}^{(i)}=
{\bf w}_{\tau_{j+1}}^{(i)}-{\bf w}_{\tau_{j}}^{(i)}$
$(i=0,\ 1,\ldots,m),$
$\left\{\tau_{j}\right\}_{j=0}^{N}$ is the partition of
the interval $[t,T],$ which satisfies the condition {\rm (\ref{1111})}.}

\vspace{2mm}

Consider trans\-for\-med particular cases for 
$k=1,\ldots,5$ of Theorem 4 \cite{3a} (2006), \cite{3aa}-\cite{art-zero}

\vspace{1mm}
\begin{equation}
\label{a1}
J[\psi^{(1)}]_{T,t}
=\hbox{\vtop{\offinterlineskip\halign{
\hfil#\hfil\cr
{\rm l.i.m.}\cr
$\stackrel{}{{}_{p_1\to \infty}}$\cr
}} }\sum_{j_1=0}^{p_1}
C_{j_1}\zeta_{j_1}^{(i_1)},
\end{equation}

\vspace{3mm}
\begin{equation}
\label{leto5001}
J[\psi^{(2)}]_{T,t}
=\hbox{\vtop{\offinterlineskip\halign{
\hfil#\hfil\cr
{\rm l.i.m.}\cr
$\stackrel{}{{}_{p_1,p_2\to \infty}}$\cr
}} }\sum_{j_1=0}^{p_1}\sum_{j_2=0}^{p_2}
C_{j_2j_1}\Biggl(\zeta_{j_1}^{(i_1)}\zeta_{j_2}^{(i_2)}
-{\bf 1}_{\{i_1=i_2\ne 0\}}
{\bf 1}_{\{j_1=j_2\}}\Biggr),
\end{equation}

\vspace{8mm}
$$
J[\psi^{(3)}]_{T,t}=
\hbox{\vtop{\offinterlineskip\halign{
\hfil#\hfil\cr
{\rm l.i.m.}\cr
$\stackrel{}{{}_{p_1,\ldots,p_3\to \infty}}$\cr
}} }\sum_{j_1=0}^{p_1}\sum_{j_2=0}^{p_2}\sum_{j_3=0}^{p_3}
C_{j_3j_2j_1}\Biggl(
\zeta_{j_1}^{(i_1)}\zeta_{j_2}^{(i_2)}\zeta_{j_3}^{(i_3)}
-\Biggr.
$$
\begin{equation}
\label{leto5002}
\Biggl.-{\bf 1}_{\{i_1=i_2\ne 0\}}
{\bf 1}_{\{j_1=j_2\}}
\zeta_{j_3}^{(i_3)}
-{\bf 1}_{\{i_2=i_3\ne 0\}}
{\bf 1}_{\{j_2=j_3\}}
\zeta_{j_1}^{(i_1)}-
{\bf 1}_{\{i_1=i_3\ne 0\}}
{\bf 1}_{\{j_1=j_3\}}
\zeta_{j_2}^{(i_2)}\Biggr),
\end{equation}

\vspace{8mm}
$$
J[\psi^{(4)}]_{T,t}
=
\hbox{\vtop{\offinterlineskip\halign{
\hfil#\hfil\cr
{\rm l.i.m.}\cr
$\stackrel{}{{}_{p_1,\ldots,p_4\to \infty}}$\cr
}} }\sum_{j_1=0}^{p_1}\ldots\sum_{j_4=0}^{p_4}
C_{j_4\ldots j_1}\Biggl(
\prod_{l=1}^4\zeta_{j_l}^{(i_l)}
\Biggr.
-
$$
$$
-
{\bf 1}_{\{i_1=i_2\ne 0\}}
{\bf 1}_{\{j_1=j_2\}}
\zeta_{j_3}^{(i_3)}
\zeta_{j_4}^{(i_4)}
-
{\bf 1}_{\{i_1=i_3\ne 0\}}
{\bf 1}_{\{j_1=j_3\}}
\zeta_{j_2}^{(i_2)}
\zeta_{j_4}^{(i_4)}-
$$
$$
-
{\bf 1}_{\{i_1=i_4\ne 0\}}
{\bf 1}_{\{j_1=j_4\}}
\zeta_{j_2}^{(i_2)}
\zeta_{j_3}^{(i_3)}
-
{\bf 1}_{\{i_2=i_3\ne 0\}}
{\bf 1}_{\{j_2=j_3\}}
\zeta_{j_1}^{(i_1)}
\zeta_{j_4}^{(i_4)}-
$$
$$
-
{\bf 1}_{\{i_2=i_4\ne 0\}}
{\bf 1}_{\{j_2=j_4\}}
\zeta_{j_1}^{(i_1)}
\zeta_{j_3}^{(i_3)}
-
{\bf 1}_{\{i_3=i_4\ne 0\}}
{\bf 1}_{\{j_3=j_4\}}
\zeta_{j_1}^{(i_1)}
\zeta_{j_2}^{(i_2)}+
$$
$$
+
{\bf 1}_{\{i_1=i_2\ne 0\}}
{\bf 1}_{\{j_1=j_2\}}
{\bf 1}_{\{i_3=i_4\ne 0\}}
{\bf 1}_{\{j_3=j_4\}}
+
$$
$$
+
{\bf 1}_{\{i_1=i_3\ne 0\}}
{\bf 1}_{\{j_1=j_3\}}
{\bf 1}_{\{i_2=i_4\ne 0\}}
{\bf 1}_{\{j_2=j_4\}}+
$$
\begin{equation}
\label{leto5003}
+\Biggl.
{\bf 1}_{\{i_1=i_4\ne 0\}}
{\bf 1}_{\{j_1=j_4\}}
{\bf 1}_{\{i_2=i_3\ne 0\}}
{\bf 1}_{\{j_2=j_3\}}\Biggr),
\end{equation}

\vspace{9mm}
$$
J[\psi^{(5)}]_{T,t}
=\hbox{\vtop{\offinterlineskip\halign{
\hfil#\hfil\cr
{\rm l.i.m.}\cr
$\stackrel{}{{}_{p_1,\ldots,p_5\to \infty}}$\cr
}} }\sum_{j_1=0}^{p_1}\ldots\sum_{j_5=0}^{p_5}
C_{j_5\ldots j_1}\Biggl(
\prod_{l=1}^5\zeta_{j_l}^{(i_l)}
-\Biggr.
$$
$$
-
{\bf 1}_{\{i_1=i_2\ne 0\}}
{\bf 1}_{\{j_1=j_2\}}
\zeta_{j_3}^{(i_3)}
\zeta_{j_4}^{(i_4)}
\zeta_{j_5}^{(i_5)}-
{\bf 1}_{\{i_1=i_3\ne 0\}}
{\bf 1}_{\{j_1=j_3\}}
\zeta_{j_2}^{(i_2)}
\zeta_{j_4}^{(i_4)}
\zeta_{j_5}^{(i_5)}-
$$
$$
-
{\bf 1}_{\{i_1=i_4\ne 0\}}
{\bf 1}_{\{j_1=j_4\}}
\zeta_{j_2}^{(i_2)}
\zeta_{j_3}^{(i_3)}
\zeta_{j_5}^{(i_5)}-
{\bf 1}_{\{i_1=i_5\ne 0\}}
{\bf 1}_{\{j_1=j_5\}}
\zeta_{j_2}^{(i_2)}
\zeta_{j_3}^{(i_3)}
\zeta_{j_4}^{(i_4)}-
$$
$$
-
{\bf 1}_{\{i_2=i_3\ne 0\}}
{\bf 1}_{\{j_2=j_3\}}
\zeta_{j_1}^{(i_1)}
\zeta_{j_4}^{(i_4)}
\zeta_{j_5}^{(i_5)}-
{\bf 1}_{\{i_2=i_4\ne 0\}}
{\bf 1}_{\{j_2=j_4\}}
\zeta_{j_1}^{(i_1)}
\zeta_{j_3}^{(i_3)}
\zeta_{j_5}^{(i_5)}-
$$
$$
-
{\bf 1}_{\{i_2=i_5\ne 0\}}
{\bf 1}_{\{j_2=j_5\}}
\zeta_{j_1}^{(i_1)}
\zeta_{j_3}^{(i_3)}
\zeta_{j_4}^{(i_4)}
-{\bf 1}_{\{i_3=i_4\ne 0\}}
{\bf 1}_{\{j_3=j_4\}}
\zeta_{j_1}^{(i_1)}
\zeta_{j_2}^{(i_2)}
\zeta_{j_5}^{(i_5)}-
$$
$$
-
{\bf 1}_{\{i_3=i_5\ne 0\}}
{\bf 1}_{\{j_3=j_5\}}
\zeta_{j_1}^{(i_1)}
\zeta_{j_2}^{(i_2)}
\zeta_{j_4}^{(i_4)}
-{\bf 1}_{\{i_4=i_5\ne 0\}}
{\bf 1}_{\{j_4=j_5\}}
\zeta_{j_1}^{(i_1)}
\zeta_{j_2}^{(i_2)}
\zeta_{j_3}^{(i_3)}+
$$
$$
+
{\bf 1}_{\{i_1=i_2\ne 0\}}
{\bf 1}_{\{j_1=j_2\}}
{\bf 1}_{\{i_3=i_4\ne 0\}}
{\bf 1}_{\{j_3=j_4\}}\zeta_{j_5}^{(i_5)}+
{\bf 1}_{\{i_1=i_2\ne 0\}}
{\bf 1}_{\{j_1=j_2\}}
{\bf 1}_{\{i_3=i_5\ne 0\}}
{\bf 1}_{\{j_3=j_5\}}\zeta_{j_4}^{(i_4)}+
$$
$$
+
{\bf 1}_{\{i_1=i_2\ne 0\}}
{\bf 1}_{\{j_1=j_2\}}
{\bf 1}_{\{i_4=i_5\ne 0\}}
{\bf 1}_{\{j_4=j_5\}}\zeta_{j_3}^{(i_3)}+
{\bf 1}_{\{i_1=i_3\ne 0\}}
{\bf 1}_{\{j_1=j_3\}}
{\bf 1}_{\{i_2=i_4\ne 0\}}
{\bf 1}_{\{j_2=j_4\}}\zeta_{j_5}^{(i_5)}+
$$
$$
+
{\bf 1}_{\{i_1=i_3\ne 0\}}
{\bf 1}_{\{j_1=j_3\}}
{\bf 1}_{\{i_2=i_5\ne 0\}}
{\bf 1}_{\{j_2=j_5\}}\zeta_{j_4}^{(i_4)}+
{\bf 1}_{\{i_1=i_3\ne 0\}}
{\bf 1}_{\{j_1=j_3\}}
{\bf 1}_{\{i_4=i_5\ne 0\}}
{\bf 1}_{\{j_4=j_5\}}\zeta_{j_2}^{(i_2)}+
$$
$$
+
{\bf 1}_{\{i_1=i_4\ne 0\}}
{\bf 1}_{\{j_1=j_4\}}
{\bf 1}_{\{i_2=i_3\ne 0\}}
{\bf 1}_{\{j_2=j_3\}}\zeta_{j_5}^{(i_5)}+
{\bf 1}_{\{i_1=i_4\ne 0\}}
{\bf 1}_{\{j_1=j_4\}}
{\bf 1}_{\{i_2=i_5\ne 0\}}
{\bf 1}_{\{j_2=j_5\}}\zeta_{j_3}^{(i_3)}+
$$
$$
+
{\bf 1}_{\{i_1=i_4\ne 0\}}
{\bf 1}_{\{j_1=j_4\}}
{\bf 1}_{\{i_3=i_5\ne 0\}}
{\bf 1}_{\{j_3=j_5\}}\zeta_{j_2}^{(i_2)}+
{\bf 1}_{\{i_1=i_5\ne 0\}}
{\bf 1}_{\{j_1=j_5\}}
{\bf 1}_{\{i_2=i_3\ne 0\}}
{\bf 1}_{\{j_2=j_3\}}\zeta_{j_4}^{(i_4)}+
$$
$$
+
{\bf 1}_{\{i_1=i_5\ne 0\}}
{\bf 1}_{\{j_1=j_5\}}
{\bf 1}_{\{i_2=i_4\ne 0\}}
{\bf 1}_{\{j_2=j_4\}}\zeta_{j_3}^{(i_3)}+
{\bf 1}_{\{i_1=i_5\ne 0\}}
{\bf 1}_{\{j_1=j_5\}}
{\bf 1}_{\{i_3=i_4\ne 0\}}
{\bf 1}_{\{j_3=j_4\}}\zeta_{j_2}^{(i_2)}+
$$
$$
+
{\bf 1}_{\{i_2=i_3\ne 0\}}
{\bf 1}_{\{j_2=j_3\}}
{\bf 1}_{\{i_4=i_5\ne 0\}}
{\bf 1}_{\{j_4=j_5\}}\zeta_{j_1}^{(i_1)}+
{\bf 1}_{\{i_2=i_4\ne 0\}}
{\bf 1}_{\{j_2=j_4\}}
{\bf 1}_{\{i_3=i_5\ne 0\}}
{\bf 1}_{\{j_3=j_5\}}\zeta_{j_1}^{(i_1)}+
$$
\begin{equation}
\label{a5}
+\Biggl.
{\bf 1}_{\{i_2=i_5\ne 0\}}
{\bf 1}_{\{j_2=j_5\}}
{\bf 1}_{\{i_3=i_4\ne 0\}}
{\bf 1}_{\{j_3=j_4\}}\zeta_{j_1}^{(i_1)}\Biggr),
\end{equation}

\vspace{5mm}
\noindent
where ${\bf 1}_A$ is the indicator of the set $A.$

Note that in \cite{10a}, \cite{11}, \cite{new-2023a}
Theorem~4 is generalized to the case
of an arbitrary 
complete orthonormal system of functions in $L_2([t, T])$ 
and $\psi_1(\tau),\ldots, \psi_k(\tau)\in L_2([t, T])$
(see Theorem~11 below).

From (\ref{leto5001}) 
for the case of an arbitrary 
complete orthonormal system of functions in $L_2([t, T]),$ 
$\psi_1(\tau),\psi_2(\tau)\in L_2([t, T])$
and (\ref{het7}) we obtain

\vspace{2mm}
$$
J[\psi^{(2)}]_{T,t}
=\hbox{\vtop{\offinterlineskip\halign{
\hfil#\hfil\cr
{\rm l.i.m.}\cr
$\stackrel{}{{}_{p_1,p_2\to \infty}}$\cr
}} }\sum_{j_1=0}^{p_1}\sum_{j_2=0}^{p_2}
C_{j_2j_1}\Biggl(\zeta_{j_1}^{(i_1)}\zeta_{j_2}^{(i_2)}
-{\bf 1}_{\{i_1=i_2\ne 0\}}
{\bf 1}_{\{j_1=j_2\}}\Biggr)=
$$

\vspace{4mm}

$$
=
\hbox{\vtop{\offinterlineskip\halign{
\hfil#\hfil\cr
{\rm l.i.m.}\cr
$\stackrel{}{{}_{p_1,p_2\to \infty}}$\cr
}} }\sum_{j_1=0}^{p_1}\sum_{j_2=0}^{p_2}
C_{j_2j_1}\zeta_{j_1}^{(i_1)}\zeta_{j_2}^{(i_2)}
-
{\bf 1}_{\{i_1=i_2\ne 0\}}\sum_{j_1=0}^{\infty}
C_{j_1j_1} =
$$

\vspace{2mm}

\begin{equation}
\label{het10}
=
\hbox{\vtop{\offinterlineskip\halign{
\hfil#\hfil\cr
{\rm l.i.m.}\cr
$\stackrel{}{{}_{p_1,p_2\to \infty}}$\cr
}} }\sum_{j_1=0}^{p_1}\sum_{j_2=0}^{p_2}
C_{j_2j_1}\zeta_{j_1}^{(i_1)}\zeta_{j_2}^{(i_2)}
-
\frac{1}{2}{\bf 1}_{\{i_1=i_2\ne 0\}}
\int\limits_t^T
\psi_1(\tau)\psi_2(\tau)d\tau.
\end{equation}

\vspace{3mm}

Since 
\begin{equation}
\label{uyes1}
J^{*}[\psi^{(2)}]_{T,t}=J[\psi^{(2)}]_{T,t}+
\frac{1}{2}{\bf 1}_{\{i_1=i_2\ne 0\}}
\int\limits_t^T
\psi_1(\tau)\psi_2(\tau)d\tau\ \ \ \hbox{w.\ p.\ 1},
\end{equation}

\vspace{4mm}
\noindent
where $\psi_1(\tau), \psi_2(\tau)$ are continuous functions
on $[t, T]$ (this condition is related to the definition
of the Stratonovich stochastic integral that we use \cite{1} (also see \cite{10a} (Sect.~2.1.1))),
then from (\ref{het10}) we finally get
the following expansion

$$
J^{*}[\psi^{(2)}]_{T,t}
=\hbox{\vtop{\offinterlineskip\halign{
\hfil#\hfil\cr
{\rm l.i.m.}\cr
$\stackrel{}{{}_{p_1,p_2\to \infty}}$\cr
}} }\sum_{j_1=0}^{p_1}\sum_{j_2=0}^{p_2}
C_{j_2j_1}\zeta_{j_1}^{(i_1)}\zeta_{j_2}^{(i_2)}.
$$

\vspace{5mm}

Thus, we obtain the following theorem.

\vspace{2mm}

{\bf Theorem 5} \cite{10a} (Sect.~2.1.4).\ {\it Assume that
$\{\phi_j(x)\}_{j=0}^{\infty}$ is 
an arbitrary 
complete orthonormal system of functions
in the space $L_2([t, T])$. Moreover,
$\psi_1(\tau), \psi_2(\tau)$ are continuous nonrandom functions
on $[t, T]$.
Then, for the iterated Stratonovich stochastic integral {\rm (\ref{str})}
of multiplicity $2$

\vspace{1mm}
$$
J^{*}[\psi^{(2)}]_{T,t}=
{\int\limits_t^{*}}^T \psi_2(t_2)
{\int\limits_t^{*}}^{t_2}
\psi_1(t_1) d{\bf w}_{t_1}^{(i_1)}
d{\bf w}_{t_2}^{(i_2)}\ \ \ 
(i_1,i_2=0, 1,\ldots,m)
$$

\vspace{5mm}
\noindent
the following 
expansion 

\vspace{1mm}
$$
J^{*}[\psi^{(2)}]_{T,t}
=\hbox{\vtop{\offinterlineskip\halign{
\hfil#\hfil\cr
{\rm l.i.m.}\cr
$\stackrel{}{{}_{p_1,p_2\to \infty}}$\cr
}} }\sum_{j_1=0}^{p_1}\sum_{j_2=0}^{p_2}
C_{j_2j_1}\zeta_{j_1}^{(i_1)}\zeta_{j_2}^{(i_2)}
$$

\vspace{4mm}
\noindent
converging in the mean-square sense is valid, where the Fourier coefficient 
$C_{j_2 j_1}$ has the form

$$
C_{j_2 j_1}=\int\limits_t^T\psi_2(t_2)\phi_{j_2}(t_2)
\int\limits_t^{t_2}
\psi_1(t_1)\phi_{j_1}(t_1)
dt_1dt_2
$$

\vspace{2mm}
\noindent
and
$$
\zeta_{j}^{(i)}=
\int\limits_t^T \phi_{j}(s) d{\bf w}_s^{(i)}
$$ 

\vspace{3mm}
\noindent
are independent standard Gaussian random variables for various 
$i$ or $j$ {\rm (}if $i\ne 0${\rm )},
${\bf w}_{\tau}^{(i)}={\bf f}_{\tau}^{(i)}$ are independent 
standard Wiener processes
$(i=1,\ldots,m)$ and 
${\bf w}_{\tau}^{(0)}=\tau.$}

\vspace{2mm}

Note that 
analogues of Theorem 5 for the multiplicities 3 to 6
of the iterated Stratonovich stochastic integrals (\ref{str})
and the systems of Legendre polynomials and trigonometric functions
have been formulated and proved in 
\cite{10a}, \cite{12}, \cite{15a}, \cite{hhh111hhh}, \cite{new-art-1-xxy}
(see Theorems~12--15 below).

We have 

\vspace{-1mm}
$$
J^{*}[\psi^{(2)}]_{T,t}^{p_1,p_2}\stackrel{\sf def}{=}J[\psi^{(2)}]_{T,t}^{p_1,p_2}+
\frac{1}{2}{\bf 1}_{\{i_1=i_2\ne 0\}}
\int\limits_t^T
\psi_1(s)\psi_2(s)ds=
$$

\vspace{1mm}
$$
=\sum_{j_1=0}^{p_1}\sum_{j_2=0}^{p_2}
C_{j_2j_1}\Biggl(\zeta_{j_1}^{(i_1)}\zeta_{j_2}^{(i_2)}
-{\bf 1}_{\{i_1=i_2\ne 0\}}
{\bf 1}_{\{j_1=j_2\}}\Biggr)+\frac{1}{2}{\bf 1}_{\{i_1=i_2\ne 0\}}
\int\limits_t^T
\psi_1(s)\psi_2(s)ds
=
$$

\vspace{1mm}
\begin{equation}
\label{ziko432}
=
\sum_{j_1=0}^{p_1}\sum_{j_2=0}^{p_2}
C_{j_2j_1}\zeta_{j_1}^{(i_1)}\zeta_{j_2}^{(i_2)}+{\bf 1}_{\{i_1=i_2\ne 0\}}
\left(\frac{1}{2}\int\limits_t^T
\psi_1(s)\psi_2(s)ds-\sum_{j_1=0}^{{\rm min}\{p_1,p_2\}}
C_{j_1j_1}\right),
\end{equation}

\vspace{5mm}
\noindent
where 
$$
J[\psi^{(2)}]_{T,t}^{p_1,p_2}=
\sum_{j_1=0}^{p_1}\sum_{j_2=0}^{p_2}
C_{j_2j_1}\Biggl(\zeta_{j_1}^{(i_1)}\zeta_{j_2}^{(i_2)}
-{\bf 1}_{\{i_1=i_2\ne 0\}}
{\bf 1}_{\{j_1=j_2\}}\Biggr)
$$

\vspace{5mm}
\noindent
is the approximation of iterated Ito stochastic integral 
(\ref{ito}) $(k=2)$ based on Theorem~4 (see (\ref{leto5001})).

Moreover, from (\ref{uyes1}) and (\ref{ziko432}) we obtain

\vspace{-1mm}
\begin{equation}
\label{ziko9991}
{\sf M}\left\{\left(J^{*}[\psi^{(2)}]_{T,t}-
J^{*}[\psi^{(2)}]_{T,t}^{p_1,p_2}\right)^{2n}\right\}=
{\sf M}\left\{\left(J[\psi^{(2)}]_{T,t}-
J[\psi^{(2)}]_{T,t}^{p_1,p_2}\right)^{2n}\right\}\ \to\ 0
\end{equation}

\vspace{3mm}
\noindent
if $p_1,p_2\to\infty$ $(n\in \mathbb{N})$, where the relation

\vspace{-1mm}
$$
{\sf M}\left\{\left(J[\psi^{(2)}]_{T,t}-
J[\psi^{(2)}]_{T,t}^{p_1,p_2}\right)^{2n}\right\}\ \to\ 0
$$

\vspace{3mm}
\noindent
if $p_1,p_2\to\infty$ $(n\in \mathbb{N})$ is proved in \cite{10a} (see Sect.~1.1.9, 1.12).

Further (see (\ref{ziko432})), 

\vspace{-1mm}
$$
{\sf M}\left\{\left(J^{*}[\psi^{(2)}]_{T,t}-
\sum\limits_{j_1=0}^{p_1}\sum\limits_{j_2=0}^{p_2}
C_{j_2 j_1}\zeta_{j_1}^{(i_1)}
\zeta_{j_2}^{(i_2)}
\right)^{2n}\right\}=
$$

\vspace{2mm}
$$
={\sf M}\left\{\left(J^{*}[\psi^{(2)}]_{T,t}-J^{*}[\psi^{(2)}]_{T,t}^{p_1,p_2}
+{\bf 1}_{\{i_1=i_2\ne 0\}}
\left(\frac{1}{2}\int\limits_t^T
\psi_1(s)\psi_2(s)ds-\sum_{j_1=0}^{{\rm min}\{p_1,p_2\}}
C_{j_1j_1}\right)\right)^{2n}\right\}\le
$$

\vspace{2mm}
\begin{equation}
\label{ziko9999}
\le K_n \left(\hspace{-0.3mm} 
{\sf M}\left\{\left(J^{*}[\psi^{(2)}]_{T,t}-
J^{*}[\psi^{(2)}]_{T,t}^{p_1,p_2}\right)^{2n}\right\}+
{\bf 1}_{\{i_1=i_2\ne 0\}}
\left(\frac{1}{2}\int\limits_t^T
\psi_1(s)\psi_2(s)ds-\hspace{-2mm}\sum_{j_1=0}^{{\rm min}\{p_1,p_2\}}
\hspace{-1mm}C_{j_1j_1}\right)^{\hspace{-1mm}2n}\right)\hspace{-0.5mm},
\end{equation}

\vspace{4mm}
\noindent
where constant $K_n<\infty$ depends on 
$n$.

Taking into account (\ref{het7}), (\ref{ziko9991}), and (\ref{ziko9999}), we get

\vspace{-1mm}
\begin{equation}
\label{ziko3210}
\lim\limits_{p_1,p_2\to\infty}
{\sf M}\left\{\left(J^{*}[\psi^{(2)}]_{T,t}-
\sum\limits_{j_1=0}^{p_1}\sum\limits_{j_2=0}^{p_2}
C_{j_2 j_1}\zeta_{j_1}^{(i_1)}
\zeta_{j_2}^{(i_2)}
\right)^{2n}\right\}=0.
\end{equation}

\vspace{3mm}
\noindent
Thus, we obtain the following theorem.

\vspace{2mm}

{\bf Theorem 6}\ \cite{10a} (Sect.~2.4.2).\ {\it Suppose that 
$\{\phi_j(x)\}_{j=0}^{\infty}$ is a complete orthonormal system of 
Legendre polynomials or trigonometric functions in the space $L_2([t, T]).$
Moreover,
$\psi_1(\tau), \psi_2(\tau)$ are continuous
nonrandom functions
on $[t,T].$ 
Then, for the iterated Stratonovich stochastic integral {\rm (\ref{str})}
of multiplicity $2$

\vspace{-1mm}
$$
J^{*}[\psi^{(2)}]_{T,t}=
{\int\limits_t^{*}}^T\psi_2(t_2)
{\int\limits_t^{*}}^{t_2}
\psi_1(t_1) d{\bf w}_{t_1}^{(i_1)}
d{\bf w}_{t_2}^{(i_2)}\ \ \
(i_1,i_2=0, 1,\ldots,m)
$$

\vspace{3mm}
\noindent
the following 
expansion 

\vspace{-1mm}
$$
J^{*}[\psi^{(2)}]_{T,t}=
\sum\limits_{j_1,j_2=0}^{\infty}
C_{j_2 j_1}\zeta_{j_1}^{(i_1)}
\zeta_{j_2}^{(i_2)}
$$

\vspace{3mm}
\noindent
that converges in the mean of degree $2n,$ $n\in\mathbb{N}$ {\rm (see (\ref{ziko3210}))}
is valid, where the Fourier coefficient 
$C_{j_2 j_1}$ is defined by {\rm (\ref{333.40})}
and

\vspace{-1mm}
$$
\zeta_{j}^{(i)}=
\int\limits_t^T \phi_{j}(s) d{\bf w}_s^{(i)}
$$ 

\vspace{3mm}
\noindent
are independent standard Gaussian random variables for various 
$i$ or $j$ {\rm (}in the case when $i\ne 0${\rm )},
${\bf w}_{\tau}^{(i)}={\bf f}_{\tau}^{(i)}$ are independent 
standard Wiener processes
$(i=1,\ldots,m)$ and 
${\bf w}_{\tau}^{(0)}=\tau.$}

\vspace{2mm}

Let us consider some other approaches
close to the approaches outlined in this section.
Now we turn to multiple trigonometric Fourier series converging 
almost everywhere. Let us for\-mu\-late the well-known result
from the theory of multiple
trigonometric Fourier series.

\vspace{2mm}

{\bf Theorem 7} \cite{dudu}. {\it Suppose that

\vspace{1mm}
\begin{equation}
\label{vot}
\int\limits_{[0,2\pi]^k}
|f(x_1,\ldots,x_k)|\left({\rm log}^{+}|f(x_1,\ldots,x_k)|\right)^k
{\rm log}^{+}{\rm log}^{+}|f(x_1,\ldots,x_k)|dx_1\ldots dx_k<\infty.
\end{equation}

\vspace{5mm}
\noindent
Then, for the square partial sums

$$
\sum_{j_1=0}^p\ldots \sum_{j_k=0}^p
C_{j_k\ldots j_1}\prod\limits_{l=1}^k\phi_{j_l}(x_l)
$$

\vspace{4mm}
\noindent
of the multiple trigonometric Fourier series we have

$$
\lim\limits_{p\to\infty}
\sum_{j_1=0}^p\ldots \sum_{j_k=0}^p
C_{j_k\ldots j_1}\prod\limits_{l=1}^k\phi_{j_l}(x_l)=
f(x_1,\ldots,x_k)
$$

\vspace{4mm}
\noindent
almost everywhere in $[0, 2\pi]^k,$ where 
$\{\phi_j(x)\}_{j=0}^{\infty}$ is a complete orthonormal
system of trigonometric functions
in the space $L_2([0, 2\pi]),$

$$
C_{j_k\ldots j_1}=
\int\limits_{[0,2\pi]^k}
f(x_1,\ldots,x_k)\prod\limits_{l=1}^k
\phi_{j_l}(x_l)
dx_1\ldots dx_k
$$

\vspace{4mm}
\noindent
is the Fourier coefficient of the function $f(x_1,\ldots,x_k),$ and 
${\rm log}^{+}x={\rm log} \max\{1,\ x\}.$}

\vspace{2mm}

Obviously, Theorem 7 can be reformulated for the hypercube
$[t, T]^k$ instead of the 
hypercube $[0, 2\pi]^k.$

If we tried to apply Theorem 7 in the proof of Theorem 1, 
then we would encounter the following difficulties.
Note that the right-hand side of (\ref{udar1}) contains
multiple integrals 
over hypercubes of various dimensions, namely over hypercubes
$[t, T]^k$, $[t, T]^{k-1},$ etc.
Obviously, the convergence almost everywhere in 
$[t, T]^k$ does not mean the convergence almost everywhere in 
$[t, T]^{k-1}$, $[t, T]^{k-2},$ etc.
This means that we could not apply the Lebesgue's 
Dominated Convergence Theorem 
in the proof of Lemma 6 
and thus could not  complete the proof of Theorem 1.
Although multiple series are more convenient for approximation 
than iterated series as in Theorem 1.

Suppose that $\psi_1(\tau), \psi_2(\tau)$ are continuously
differentiable functions on
$[t, T]$ and
$\{\phi_j(x)\}_{j=0}^{\infty}$ is a complete
orthonormal system of Legendre polynomials
or trigonometric functions in the space $L_2([t, T])$. 
In \cite{10a} (Sect.~2.1.2) it was shown that

\begin{equation}
\label{za1}
\lim\limits_{p_1,p_2\to\infty}
\sum_{j_1=0}^{p_1} \sum_{j_2=0}^{p_2} C_{j_2j_1}\phi_{j_1}(t_1)
\phi_{j_2}(t_1)=\frac{1}{2}\psi_1(t_1)\psi_2(t_1)=K^{*}(t_1,t_1),\ \ \ 
t_1\in (t, T),
\end{equation}

\vspace{4mm}
\noindent
where $C_{j_2j_1}$ is defined by (\ref{333.40}) ($k=2)$.

This means that we can repeat the proof of Theorem 1 for 
the case $k=2$ and apply the Lebesgue's 
Dominated Convergence Theorem 
in the formula (\ref{udar1}), since Theorem 7 and (\ref{za1})
implies the convergence almost everywhere in $[t, T]^2$ and 
almost everywhere in $[t, T]$ $(t_1=t_2\in [t, T])$
of the multiple trigonometric
Fourier series 

\begin{equation}
\label{ziko456}
\lim\limits_{p\to\infty}
\sum_{j_1=0}^{p} \sum_{j_2=0}^{p} C_{j_2j_1}\phi_{j_1}(t_1)
\phi_{j_2}(t_2),\ \ \ t_1,t_2\in [t, T]^2
\end{equation}

\vspace{4mm}
\noindent
to the function 
$K^{*}(t_1,t_2)$ 
(the question of finding an integrable majorant
for Lebesgue's 
Dominated Convergence Theorem is omitted here).
So, we can obtain the particular case of Theorem~6.

Let us consider the another approach.
The following fact is well-known.

\vspace{2mm}

{\bf Proposition 1.} {\it Let 
$\bigl\{x_{n_1,\ldots,n_k}\bigr\}_{n_1,\ldots,n_k=1}^{\infty}$
be a multi-index sequence and let there exists the limit

$$
\lim\limits_{n_1,\ldots,n_k\to\infty}x_{n_1,\ldots,n_k}<\infty.
$$

\vspace{4mm}

\noindent
Moreover, let there exists the limit

\vspace{-1mm}
$$
\lim\limits_{n_k\to\infty}x_{n_1,\ldots,n_k}=y_{n_1,\ldots,n_{k-1}}
<\infty\ \ \ \hbox{for any}\ \ \ n_1,\ldots,n_{k-1}.
$$

\vspace{3mm}

\noindent
Then there exists the iterated limit

\vspace{-1mm}
$$
\lim\limits_{n_1,\ldots,n_{k-1}\to\infty}\ 
\lim\limits_{n_k\to\infty}x_{n_1,\ldots,n_k},
$$

\vspace{3mm}
\noindent
and moreover,

\vspace{-1mm}
$$
\lim\limits_{n_1,\ldots,n_{k-1}\to\infty}\ \lim\limits_{n_k\to\infty}
x_{n_1,\ldots,n_k}=
\lim\limits_{n_1,\ldots,n_k\to\infty}x_{n_1,\ldots,n_k}.
$$
}

\vspace{3mm}

Denote

$$
C_{j_s\ldots j_1}(t_{s+1},\ldots,t_k)=
\int\limits_{[t,T]^s}K(t_1,\ldots,t_k)
\prod_{l=1}^s \phi_{j_l}(t_l)dt_1\ldots dt_s\ \ \ (s=1,\ldots,k-1).
$$

\vspace{5mm}
\noindent
where $K(t_1,\ldots,t_k)$ has the form (\ref{ppp}).
For $s=k$ we suppose that $C_{j_k\ldots j_1}$
is defined by
(\ref{333.40}).

Consider the following Fourier series

\begin{equation}
\label{ww1}
\lim\limits_{p_1,p_2\to\infty}
\sum_{j_1=0}^{p_1} \sum_{j_2=0}^{p_2} C_{j_2j_1}(t_3,\ldots,t_k)
\phi_{j_1}(t_1)
\phi_{j_2}(t_2),
\end{equation}

\vspace{1mm}
\begin{equation}
\label{ww2}
\lim\limits_{p_1,p_2,p_3\to\infty}
\sum_{j_1=0}^{p_1}\sum_{j_1=0}^{p_1}
\sum_{j_3=0}^{p_3} C_{j_3j_2j_1}(t_4,\ldots,t_k)
\phi_{j_1}(t_1)
\phi_{j_2}(t_2)\phi_{j_3}(t_3),
\end{equation}

$$
\ldots
$$

\begin{equation}
\label{ww3}
\lim\limits_{p_1,\ldots,p_{k-1}\to\infty}
\sum_{j_1=0}^{p_1}\ldots
\sum_{j_{k-1}=0}^{p_{k-1}} C_{j_{k-1}\ldots j_1}(t_k)
\phi_{j_1}(t_1)\ldots
\phi_{j_{k-1}}(t_{k-1}),
\end{equation}

\vspace{1mm}

\begin{equation}
\label{ww4}
\lim\limits_{p_1,\ldots,p_{k}\to\infty}
\sum_{j_1=0}^{p_1}\ldots
\sum_{j_{k}=0}^{p_{k}} C_{j_{k}\ldots j_1}
\phi_{j_1}(t_1)\ldots
\phi_{j_{k}}(t_{k}),
\end{equation}

\vspace{5mm}
\noindent
where
$\{\phi_j(x)\}_{j=0}^{\infty}$ is a complete orthonormal
system of Legendre polynomials or trigonometric func\-ti\-ons
in the space $L_2([t, T])$. 

The author does not know the answere to the question
on existence of the limits (\ref{ww1})--(\ref{ww4}) even for
the case $p_1=\ldots=p_k$ and trigonometric Fourier series.
Obviously, at least 
for the case $k=2$ and $\psi_1(\tau),$ $\psi_2(\tau)\equiv 1$
the answere to the above question is positive
for the Fourier--Legendre series as well as 
for the trigonometric Fourier series.

If we suppose the 
existence of the limits (\ref{ww1})--(\ref{ww4}), then
combining Proposition 1 and the proof of Lemma 1 
we obtain

$$
K^{*}(t_1,\ldots,t_k)=
\sum_{j_1=0}^{\infty}C_{j_1}(t_2,\ldots,t_k)
\phi_{j_1}(t_1)=
$$

\vspace{3mm}
\begin{equation}
\label{ww10}
=\sum_{j_1=0}^{\infty}\sum_{j_2=0}^{\infty}C_{j_2j_1}(t_3,\ldots,t_k)
\phi_{j_1}(t_1)\phi_{j_2}(t_2)=
\end{equation}

\vspace{3mm}
$$
=
\lim\limits_{p_1,p_2\to\infty}
\sum_{j_1=0}^{p_1} \sum_{j_2=0}^{p_2} C_{j_2j_1}(t_3,\ldots,t_k)
\phi_{j_1}(t_1)
\phi_{j_2}(t_2)=
$$

\vspace{3mm}
$$
=
\lim\limits_{p_1,p_2\to\infty}
\sum_{j_1=0}^{p_1} \sum_{j_2=0}^{p_2} 
\sum_{j_3=0}^{\infty}C_{j_3j_2j_1}(t_4,\ldots,t_k)
\phi_{j_1}(t_1)
\phi_{j_2}(t_2)\phi_{j_3}(t_3)=
$$

\vspace{3mm}
\begin{equation}
\label{ww18}
=
\lim\limits_{p_1,p_2,p_3\to\infty}
\sum_{j_1=0}^{p_1} \sum_{j_2=0}^{p_2} 
\sum_{j_3=0}^{p_3}C_{j_3j_2j_1}(t_4,\ldots,t_k)
\phi_{j_1}(t_1)
\phi_{j_2}(t_2)\phi_{j_3}(t_3)=
\end{equation}

\vspace{3mm}
\begin{equation}
\label{ww181}
=\sum_{j_1=0}^{\infty}\sum_{j_2=0}^{\infty}\sum_{j_3=0}^{\infty}
C_{j_3j_2j_1}(t_4,\ldots,t_k)
\phi_{j_1}(t_1)\phi_{j_2}(t_2)\phi_{j_3}(t_3)=
\end{equation}

\vspace{3mm}
\begin{equation}
\label{ww182}
=
\lim\limits_{p_1,p_2,p_3\to\infty}
\sum_{j_1=0}^{p_1} \sum_{j_2=0}^{p_2}\sum_{j_3=0}^{p_3} 
\sum_{j_4=0}^{\infty}C_{j_4\ldots j_1}(t_5,\ldots,t_k)
\phi_{j_1}(t_1)
\phi_{j_3}(t_4)= 
\end{equation}

\vspace{2mm}
$$
=\ldots =
$$

\vspace{1mm}
\begin{equation}
\label{ww21}
=
\lim\limits_{p_1,\ldots,p_k\to\infty}
\sum_{j_1=0}^{p_1} \ldots
\sum_{j_k=0}^{p_k}C_{j_k\ldots j_1}
\phi_{j_1}(t_1)\ldots
\phi_{j_k}(t_k).
\end{equation}

\vspace{6mm}

Note that 
the transition from (\ref{ww18}) to (\ref{ww181})
is based on (\ref{ww10}) and the proof of Lemma 1.
The transition from (\ref{ww181}) to (\ref{ww182})
is based on (\ref{ww18}) and the proof of Lemma 1.

Using (\ref{ww21}) we could get the version of Theorem 1
with multiple series instead of iterated ones.

\vspace{5mm}

\section{Refinement of Theorems 1 and 2 for Iterated 
Stra\-to\-no\-vich Stochastic Integrals of Multiplicities $2$ and $3$
$(i_1,i_2,i_3=1,\ldots,m).$
The Case of Mean-Square Convergence}

\vspace{5mm}

In this section, it will be shown that the upper limits
in Theorems 1 and 2 (the cases $k=2,$ $k=3$ and $n=1$)
can be replaced by the usual limits.

\vspace{2mm}

{\bf Theorem 8} \cite{10a} (Sect.~2.4).\ 
{\it Suppose that every $\psi_l(\tau)$ $(l=1,2,3)$ is twice continuously
differentiable function at the interval
$[t, T]$ and
$\{\phi_j(x)\}_{j=0}^{\infty}$ is a complete
orthonormal system of 
trigonometric functions in the space $L_2([t, T])$. 
Then$,$ the iterated Stratonovich stochastic integrals 
$J^{*}[\psi^{(2)}]_{T,t}$ and $J^{*}[\psi^{(3)}]_{T,t}$ $(i_1,i_2,i_3=1,\ldots,m)$
defined by {\rm(\ref{str})}
are expanded into the 
conver\-ging 
in the mean-square sense 
iterated series

\vspace{-1mm}
\begin{equation}
\label{nov500}
\lim\limits_{p_1\to\infty}
\lim\limits_{p_2\to\infty}
{\sf M}\left\{\left(J^{*}[\psi^{(2)}]_{T,t}-
\sum_{j_1=0}^{p_1}\sum_{j_2=0}^{p_2}
C_{j_2j_1}\zeta^{(i_1)}_{j_1}\zeta^{(i_2)}_{j_2}\right)^{2}\right\}=0,
\end{equation}

\begin{equation}
\label{nov501}
\lim\limits_{p_2\to\infty}
\lim\limits_{p_1\to\infty}
{\sf M}\left\{\left(J^{*}[\psi^{(2)}]_{T,t}-
\sum_{j_2=0}^{p_2}\sum_{j_1=0}^{p_1}
C_{j_2j_1}\zeta^{(i_1)}_{j_1}\zeta^{(i_2)}_{j_2}\right)^{2}\right\}=0,
\end{equation}

\begin{equation}
\label{nov502}
\lim\limits_{p_1\to\infty}
\lim\limits_{p_2\to\infty}
\lim\limits_{p_3\to\infty}
{\sf M}\left\{\left(J^{*}[\psi^{(3)}]_{T,t}-
\sum_{j_1=0}^{p_1}\sum_{j_2=0}^{p_2}\sum_{j_3=0}^{p_3}
C_{j_3j_2j_1}\zeta^{(i_1)}_{j_1}\zeta^{(i_2)}_{j_2}\zeta^{(i_3)}_{j_3}\right)^{2}\right\}=0,
\end{equation}

\begin{equation}
\label{nov503}
\lim\limits_{p_3\to\infty}
\lim\limits_{p_2\to\infty}
\lim\limits_{p_1\to\infty}
{\sf M}\left\{\left(J^{*}[\psi^{(3)}]_{T,t}-
\sum_{j_3=0}^{p_3}\sum_{j_2=0}^{p_2}\sum_{j_1=0}^{p_1}
C_{j_3j_2j_1}\zeta^{(i_1)}_{j_1}\zeta^{(i_2)}_{j_2}\zeta^{(i_3)}_{j_3}\right)^{2}\right\}=0,
\end{equation}

\vspace{4mm}
\noindent
where 
$$
\zeta_{j}^{(i)}=
\int\limits_t^T \phi_{j}(s) d{\bf f}_s^{(i)}\ \ \ (i=1,\ldots,m,\ \ j=0, 1,\ldots )
$$ 

\vspace{3mm}
\noindent
are independent standard Gaussian random variables
for
various
$i$ or $j$ and 
$C_{j_2j_1},$ $C_{j_3j_2j_1}$ are defined by 
{\rm (\ref{333.40})}.}

\vspace{2mm}

{\bf Proof.}\ We will prove the equalities (\ref{nov500}) and (\ref{nov502})
(the equalities (\ref{nov501}) and (\ref{nov503}) can be proved similarly
using the expansion (\ref{otit3333}) instead of the expansion
(\ref{30.18})).

From (\ref{nov800}) we have
w.~p.~1 

$$
J^{*}[\psi^{(2)}]_{T,t}-
\sum_{j_1=0}^{p_1}\sum_{j_2=0}^{p_2}
C_{j_2j_1}\zeta^{(i_1)}_{j_1}\zeta^{(i_2)}_{j_2}=J[R_{p_1p_2}]_{T,t}^{(2)}=
$$

\vspace{1mm}
$$
=\int\limits_t^T\int\limits_t^{t_2}
R_{p_1p_2}(t_1,t_2)d{\bf f}_{t_1}^{(i_1)}d{\bf f}_{t_2}^{(i_2)}
+\int\limits_t^T\int\limits_t^{t_1}
R_{p_1p_2}(t_1,t_2)d{\bf f}_{t_2}^{(i_2)}d{\bf f}_{t_1}^{(i_1)}+
$$

\begin{equation}
\label{nov801}
+{\bf 1}_{\{i_1=i_2\}}
\int\limits_t^T R_{p_1p_2}(t_1,t_1)dt_1,
\end{equation}

\vspace{3mm}
\noindent
where we used the same notations as in (\ref{nov800}).

Uning (\ref{nov801}), we obtain

\vspace{-1mm}
$$
{\sf M}\left\{\left(J[R_{p_1p_2}]_{T,t}^{(2)}\right)^2\right\}=
\int\limits_t^T\int\limits_t^{t_2}
R_{p_1p_2}^2(t_1,t_2)dt_1 dt_2
+\int\limits_t^T\int\limits_t^{t_1}
R_{p_1p_2}^2(t_1,t_2)dt_2dt_1+
$$

\vspace{1mm}
$$
+{\bf 1}_{\{i_1=i_2\}} \left(2\int\limits_t^T\int\limits_t^{t_2}
R_{p_1p_2}(t_1,t_2)R_{p_1p_2}(t_2,t_1)dt_1dt_2+
\left(\int\limits_t^T R_{p_1p_2}(t_1,t_1)dt_1\right)^2\right)=
$$

\vspace{1mm}
$$
=\int\limits_t^T\int\limits_t^{t_2}
R_{p_1p_2}^2(t_1,t_2)dt_1 dt_2
+\int\limits_t^T\int\limits_{t_2}^{T}
R_{p_1p_2}^2(t_1,t_2)dt_1dt_2+
$$

\vspace{1mm}
$$
+{\bf 1}_{\{i_1=i_2\}} \left(\int\limits_t^T\int\limits_t^{t_2}
R_{p_1p_2}(t_1,t_2)R_{p_1p_2}(t_2,t_1)dt_1dt_2+\right.
$$

\vspace{1mm}
$$
+\left.
\int\limits_t^T\int\limits_{t_1}^T
R_{p_1p_2}(t_1,t_2)R_{p_1p_2}(t_2,t_1)dt_2dt_1\right)+
{\bf 1}_{\{i_1=i_2\}}\left(\int\limits_t^T R_{p_1p_2}(t_1,t_1)dt_1\right)^2=
$$

\vspace{1mm}
$$
=\int\limits_{[t,T]^2}
R_{p_1p_2}^2(t_1,t_2)dt_1 dt_2+
$$

\vspace{1mm}
$$
+{\bf 1}_{\{i_1=i_2\}} \left(\int\limits_t^T\int\limits_t^{t_2}
R_{p_1p_2}(t_1,t_2)R_{p_1p_2}(t_2,t_1)dt_1dt_2+\right.
$$

\vspace{1mm}
$$
+\left.
\int\limits_t^T\int\limits_{t_2}^T
R_{p_1p_2}(t_1,t_2)R_{p_1p_2}(t_2,t_1)dt_1dt_2\right)+
{\bf 1}_{\{i_1=i_2\}}\left(\int\limits_t^T R_{p_1p_2}(t_1,t_1)dt_1\right)^2=
$$

\vspace{1mm}
$$
=\int\limits_{[t,T]^2}
R_{p_1p_2}^2(t_1,t_2)dt_1 dt_2+
$$

\vspace{1mm}
\begin{equation}
\label{nov803}
+{\bf 1}_{\{i_1=i_2\}} \left(\int\limits_{[t,T]^2}
R_{p_1p_2}(t_1,t_2)R_{p_1p_2}(t_2,t_1)dt_1dt_2+
\left(\int\limits_t^T R_{p_1p_2}(t_1,t_1)dt_1\right)^2\right).
\end{equation}

\vspace{5mm}

Since the integrals on the right-hand side of (\ref{nov803}) 
exist as Riemann integrals, then they are equal to the 
corresponding Lebesgue integrals. 
Moreover,

$$
\lim\limits_{p_1\to\infty}\lim\limits_{p_2\to\infty}
R_{p_1p_2}(t_1,t_2)=0\ \ \ \hbox{when}\ \ \ (t_1,t_2)\in (t, T)^2,
$$

\vspace{3mm}
\noindent
where the left-hand side 
is bounded on $[t, T]^2$ (see (\ref{410})).

Then, applying two times (we mean an iterated passage to the limit
$\lim\limits_{p_1\to\infty}\lim\limits_{p_2\to\infty}$)
the Lebesgue's 
Dominated Convergence Theorem and taking into account
(\ref{leto8001}), (\ref{leto8002}), and (\ref{d2020}),
we obtain

\begin{equation}
\label{nov805}
\lim\limits_{p_1\to\infty}\lim\limits_{p_2\to\infty}
\int\limits_{[t, T]^2}
R_{p_1p_2}^2(t_1,t_2)dt_1 dt_2=0,
\end{equation}

\vspace{2mm}
\begin{equation}
\label{nov806}
\lim\limits_{p_1\to\infty}\lim\limits_{p_2\to\infty}
\int\limits_{[t, T]^2}
R_{p_1p_2}(t_1,t_2)R_{p_1p_2}(t_2,t_1)dt_1 dt_2=0,
\end{equation}

\vspace{1mm}
\begin{equation}
\label{nov807}
\lim\limits_{p_1\to\infty}\lim\limits_{p_2\to\infty}
\int\limits_t^T
R_{p_1p_2}(t_1,t_1)dt_1=0.
\end{equation}

\vspace{5mm}

The relations (\ref{nov803})--(\ref{nov807}) imply the following equality

\vspace{1mm}
$$
\lim\limits_{p_1\to\infty}\lim\limits_{p_2\to\infty}
{\sf M}\left\{\left(J[R_{p_1p_2}]_{T,t}^{(2)}\right)^2\right\}=0.
$$

\vspace{4mm}
\noindent
The relation (\ref{nov500}) is proved.

Let us prove the relation (\ref{nov502}).
Using (\ref{s1s}) and the
integration order replacement technique for 
iterated Ito stochastic integrals (see Chapter 3 in \cite{10a}-\cite{12aa-afterxxx}), we get 
w.~p.~1

\vspace{2mm}
$$
J^{*}[\psi^{(3)}]_{T,t}-
\sum_{j_1=0}^{p_1}\sum_{j_2=0}^{p_2}\sum_{j_3=0}^{p_3}
C_{j_3j_2j_1}\zeta^{(i_1)}_{j_1}\zeta^{(i_2)}_{j_2}\zeta^{(i_3)}_{j_3}
=J[R_{p_1p_2p_3}]_{T,t}^{(3)}=
$$

\vspace{1mm}
$$
=
\int\limits_t^T\int\limits_t^{t_3}\int\limits_t^{t_2}
R_{p_1 p_2 p_3}(t_1,t_2,t_3)
d{\bf f}_{t_1}^{(i_1)}
d{\bf f}_{t_2}^{(i_2)}
d{\bf f}_{t_3}^{(i_3)}+
$$

\vspace{1mm}
$$
+
\int\limits_t^T\int\limits_t^{t_3}\int\limits_t^{t_2}
R_{p_1 p_2 p_3}(t_1,t_3,t_2)
d{\bf f}_{t_1}^{(i_1)}
d{\bf f}_{t_2}^{(i_3)}
d{\bf f}_{t_3}^{(i_2)}+
$$

\vspace{1mm}
$$
+
\int\limits_t^T\int\limits_t^{t_3}\int\limits_t^{t_2}
R_{p_1 p_2 p_3}(t_2,t_1,t_3)
d{\bf f}_{t_1}^{(i_2)}
d{\bf f}_{t_2}^{(i_1)}
d{\bf f}_{t_3}^{(i_3)}+
$$

\vspace{1mm}
$$
+
\int\limits_t^T\int\limits_t^{t_3}\int\limits_t^{t_2}
R_{p_1 p_2 p_3}(t_2,t_3,t_1)
d{\bf f}_{t_1}^{(i_3)}
d{\bf f}_{t_2}^{(i_1)}
d{\bf f}_{t_3}^{(i_2)}+
$$

\vspace{1mm}
$$
+
\int\limits_t^T\int\limits_t^{t_3}\int\limits_t^{t_2}
R_{p_1 p_2 p_3}(t_3,t_2,t_1)
d{\bf f}_{t_1}^{(i_3)}
d{\bf f}_{t_2}^{(i_2)}
d{\bf f}_{t_3}^{(i_1)}+
$$

\vspace{1mm}
$$
+
\int\limits_t^T\int\limits_t^{t_3}\int\limits_t^{t_2}
R_{p_1 p_2 p_3}(t_3,t_1,t_2)
d{\bf f}_{t_1}^{(i_2)}
d{\bf f}_{t_2}^{(i_3)}
d{\bf f}_{t_3}^{(i_1)}+
$$

\vspace{1mm}
$$
+{\bf 1}_{\{i_1=i_2\}}
\int\limits_t^T\left(\int\limits_t^{T}
R_{p_1 p_2 p_3}(t_2,t_2,t_3)dt_2\right)
d{\bf f}_{t_3}^{(i_3)}+
$$

\vspace{1mm}
$$
+
{\bf 1}_{\{i_2=i_3\}}
\int\limits_t^T\left(\int\limits_t^{T}
R_{p_1 p_2 p_3}(t_1,t_2,t_2)dt_2\right)
d{\bf f}_{t_1}^{(i_1)}+
$$

\vspace{1mm}
\begin{equation}
\label{nov901}
+{\bf 1}_{\{i_1=i_3\}}
\int\limits_t^T\left(\int\limits_t^{T}
R_{p_1 p_2 p_3}(t_3,t_2,t_3)dt_3\right)
d{\bf f}_{t_2}^{(i_2)}.
\end{equation}

\vspace{5mm} 

Let us calculate the second moment of 
$J[R_{p_1p_2p_3}]_{T,t}^{(3)}$ using (\ref{nov901}).
We have

\vspace{2mm}
$$
{\sf M}\left\{\left(J[R_{p_1p_2p_3}]_{T,t}^{(3)}\right)^2\right\}=
$$

\vspace{1mm}
\begin{equation}
\label{novv1}
=
\int\limits_t^T\int\limits_t^{t_3}\int\limits_t^{t_2}
\left(\sum\limits_{(t_1,t_2,t_3)} R_{p_1 p_2 p_3}^2(t_1,t_2,t_3)\right)dt_1
dt_2dt_3+
\end{equation}

\vspace{1mm}
$$
+2\left({\bf 1}_{\{i_1=i_2\}}
\int\limits_t^T\int\limits_t^{t_3}\int\limits_t^{t_2}
G^{(1)}_{p_1 p_2 p_3}(t_1,t_2,t_3)dt_1dt_2dt_3+\right.
$$

\vspace{1mm}
$$
+{\bf 1}_{\{i_1=i_3\}}
\int\limits_t^T\int\limits_t^{t_3}\int\limits_t^{t_2}
G^{(2)}_{p_1 p_2 p_3}(t_1,t_2,t_3)dt_1dt_2dt_3+
$$

\vspace{1mm}
$$
+{\bf 1}_{\{i_2=i_3\}}
\int\limits_t^T\int\limits_t^{t_3}\int\limits_t^{t_2}
G^{(3)}_{p_1 p_2 p_3}(t_1,t_2,t_3)dt_1dt_2dt_3+
$$

\vspace{1mm}
$$
\left.+{\bf 1}_{\{i_1=i_2=i_3\}}
\int\limits_t^T\int\limits_t^{t_3}\int\limits_t^{t_2}
G^{(4)}_{p_1 p_2 p_3}(t_1,t_2,t_3)dt_1dt_2dt_3\right)+
$$

\vspace{1mm}
$$
+\int\limits_{[t,T]^3}
\biggl({\bf 1}_{\{i_1=i_2\}}
R_{p_1 p_2 p_3}(t_1,t_1,t_3)R_{p_1 p_2 p_3}(t_2,t_2,t_3)+\biggr.
$$

\vspace{1mm}
$$
+{\bf 1}_{\{i_2=i_3\}}
R_{p_1 p_2 p_3}(t_3,t_1,t_1)R_{p_1 p_2 p_3}(t_3,t_2,t_2)+
$$

\vspace{1mm}
$$
+{\bf 1}_{\{i_1=i_3\}}
R_{p_1 p_2 p_3}(t_1,t_3,t_1)R_{p_1 p_2 p_3}(t_2,t_3,t_2)+
$$

\vspace{1mm}
$$
+2\cdot {\bf 1}_{\{i_1=i_2=i_3\}}\biggl(
R_{p_1 p_2 p_3}(t_1,t_1,t_3)R_{p_1 p_2 p_3}(t_3,t_2,t_2)+\biggr.
$$

\vspace{1mm}
$$
+
R_{p_1 p_2 p_3}(t_1,t_1,t_3)R_{p_1 p_2 p_3}(t_2,t_3,t_2)+
$$

\vspace{1mm}
\begin{equation}
\label{nov980}
\biggl.\biggl.+
R_{p_1 p_2 p_3}(t_3,t_1,t_1)R_{p_1 p_2 p_3}(t_2,t_3,t_2)\biggr)\biggr)dt_1dt_2dt_3,
\end{equation}

\vspace{4mm}
\noindent
where permutation $(t_1,t_2,t_3)$ when summing in (\ref{novv1}) are 
performed only in the value $R_{p_1 p_2 p_3}^2(t_1,t_2,t_3)$ and 
the functions $G^{(i)}_{p_1 p_2 p_3}(t_1,t_2,t_3)$ $(i=1,\ldots,4)$
are defined by the following relations

\vspace{1mm}
$$
G^{(1)}_{p_1 p_2 p_3}(t_1,t_2,t_3)=
R_{p_1 p_2 p_3}(t_1,t_2,t_3)R_{p_1 p_2 p_3}(t_2,t_1,t_3)+
$$

\vspace{1mm}
$$
+R_{p_1 p_2 p_3}(t_1,t_3,t_2)R_{p_1 p_2 p_3}(t_3,t_1,t_2)+
$$

\vspace{1mm}
$$
+R_{p_1 p_2 p_3}(t_2,t_3,t_1)R_{p_1 p_2 p_3}(t_3,t_2,t_1),
$$

\vspace{3mm}
$$
G^{(2)}_{p_1 p_2 p_3}(t_1,t_2,t_3)=
R_{p_1 p_2 p_3}(t_1,t_2,t_3)R_{p_1 p_2 p_3}(t_3,t_2,t_1)+
$$

\vspace{1mm}
$$
+R_{p_1 p_2 p_3}(t_1,t_3,t_2)R_{p_1 p_2 p_3}(t_2,t_3,t_1)+
$$

\vspace{1mm}
$$
+R_{p_1 p_2 p_3}(t_2,t_1,t_3)R_{p_1 p_2 p_3}(t_3,t_1,t_2),
$$

\vspace{3mm}
$$
G^{(3)}_{p_1 p_2 p_3}(t_1,t_2,t_3)=
R_{p_1 p_2 p_3}(t_1,t_2,t_3)R_{p_1 p_2 p_3}(t_1,t_3,t_2)+
$$

\vspace{1mm}
$$
+R_{p_1 p_2 p_3}(t_2,t_1,t_3)R_{p_1 p_2 p_3}(t_2,t_3,t_1)+
$$

\vspace{1mm}
$$
+R_{p_1 p_2 p_3}(t_3,t_2,t_1)R_{p_1 p_2 p_3}(t_3,t_1,t_2),
$$

\vspace{3mm}
$$
G^{(4)}_{p_1 p_2 p_3}(t_1,t_2,t_3)=
R_{p_1 p_2 p_3}(t_1,t_2,t_3)R_{p_1 p_2 p_3}(t_2,t_3,t_1)+
$$

\vspace{1mm}
$$
+R_{p_1 p_2 p_3}(t_1,t_2,t_3)R_{p_1 p_2 p_3}(t_3,t_1,t_2)+
$$

\vspace{1mm}
$$
+R_{p_1 p_2 p_3}(t_1,t_3,t_2)R_{p_1 p_2 p_3}(t_2,t_1,t_3)+
$$

\vspace{1mm}
$$
+R_{p_1 p_2 p_3}(t_1,t_3,t_2)R_{p_1 p_2 p_3}(t_3,t_2,t_1)+
$$

\vspace{1mm}
$$
+R_{p_1 p_2 p_3}(t_2,t_1,t_3)R_{p_1 p_2 p_3}(t_3,t_2,t_1)+
$$

\vspace{1mm}
$$
+R_{p_1 p_2 p_3}(t_2,t_3,t_1)R_{p_1 p_2 p_3}(t_3,t_1,t_2).
$$

\vspace{5mm}

Further,

\vspace{-2mm}
$$
\int\limits_t^T\int\limits_t^{t_3}\int\limits_t^{t_2}
\left(\sum\limits_{(t_1,t_2,t_3)} R_{p_1 p_2 p_3}^2(t_1,t_2,t_3)\right)dt_1
dt_2dt_3=
$$

\vspace{2mm}
\begin{equation}
\label{novv60}
=\int\limits_{[t,T]^3}R_{p_1 p_2 p_3}^2(t_1,t_2,t_3)dt_1
dt_2dt_3.
\end{equation}

\vspace{5mm}

We will say that the function $\Phi(t_1,t_2,t_3)$ is symmetric if

\vspace{1mm}
$$
\Phi(t_1,t_2,t_3)=\Phi(t_1,t_3,t_2)=\Phi(t_2,t_1,t_3)=
\Phi(t_2,t_3,t_1)=
$$

\vspace{1mm}
$$
=\Phi(t_3,t_1,t_2)=\Phi(t_3,t_2,t_1).
$$

\vspace{5mm}

For the symmetric function $\Phi(t_1,t_2,t_3)$, we have

\vspace{1mm}
$$
\int\limits_t^T\int\limits_t^{t_3}\int\limits_t^{t_2}
\left(\sum\limits_{(t_1,t_2,t_3)} \Phi(t_1,t_2,t_3) \right)dt_1
dt_2dt_3= 
$$

\vspace{1mm}
$$
=6\int\limits_t^T\int\limits_t^{t_3}\int\limits_t^{t_2}
\Phi(t_1,t_2,t_3)dt_1dt_2dt_3= 
$$

\vspace{1mm}
\begin{equation}
\label{nov900}
=
\int\limits_{[t,T]^3}\Phi(t_1,t_2,t_3)dt_1
dt_2dt_3.
\end{equation}

\vspace{4mm} 

The relation (\ref{nov900}) implies that

\vspace{1mm}
\begin{equation}
\label{nov901a}
\int\limits_t^T\int\limits_t^{t_3}\int\limits_t^{t_2}
\Phi(t_1,t_2,t_3)dt_1dt_2dt_3=
\frac{1}{6}\int\limits_{[t,T]^3}\Phi(t_1,t_2,t_3)dt_1
dt_2dt_3.
\end{equation}

\vspace{4mm}

It is easy to check that the functions 
$G^{(i)}_{p_1 p_2 p_3}(t_1,t_2,t_3)$ $(i=1,\ldots,4)$
are symmetric. Using this property as well as 
(\ref{nov980}), (\ref{novv60}), and (\ref{nov901a}), we obtain

\vspace{1mm}
$$
{\sf M}\left\{\left(J[R_{p_1p_2p_3}]_{T,t}^{(3)}\right)^2\right\}=
\int\limits_{[t,T]^3}R_{p_1 p_2 p_3}^2(t_1,t_2,t_3)dt_1dt_2dt_3+
$$

\vspace{1mm}
$$
+\frac{1}{3}\int\limits_{[t,T]^3}\biggl({\bf 1}_{\{i_1=i_2\}}
G^{(1)}_{p_1 p_2 p_3}(t_1,t_2,t_3)dt_1dt_2dt_3+\biggr.
$$

\vspace{1mm}
$$
+{\bf 1}_{\{i_1=i_3\}}
G^{(2)}_{p_1 p_2 p_3}(t_1,t_2,t_3)dt_1dt_2dt_3+
$$

\vspace{1mm}
$$
+{\bf 1}_{\{i_2=i_3\}}
G^{(3)}_{p_1 p_2 p_3}(t_1,t_2,t_3)dt_1dt_2dt_3+
$$

\vspace{1mm}
$$
\biggl.+{\bf 1}_{\{i_1=i_2=i_3\}}
G^{(4)}_{p_1 p_2 p_3}(t_1,t_2,t_3)dt_1dt_2dt_3\biggr)dt_1dt_2dt_3+
$$

\vspace{1mm}
$$
+
\int\limits_{[t,T]^3}\biggl({\bf 1}_{\{i_1=i_2\}}
R_{p_1 p_2 p_3}(t_1,t_1,t_3)R_{p_1 p_2 p_3}(t_2,t_2,t_3)+\biggr.
$$

\vspace{1mm}
$$
+{\bf 1}_{\{i_2=i_3\}}
R_{p_1 p_2 p_3}(t_3,t_1,t_1)R_{p_1 p_2 p_3}(t_3,t_2,t_2)+
$$

\vspace{1mm}
$$
+{\bf 1}_{\{i_1=i_3\}}
R_{p_1 p_2 p_3}(t_1,t_3,t_1)R_{p_1 p_2 p_3}(t_2,t_3,t_2)+
$$

\vspace{1mm}
$$
+2\cdot {\bf 1}_{\{i_1=i_2=i_3\}}\biggl(
R_{p_1 p_2 p_3}(t_1,t_1,t_3)R_{p_1 p_2 p_3}(t_3,t_2,t_2)+\biggr.
$$

\vspace{1mm}
$$
+
R_{p_1 p_2 p_3}(t_1,t_1,t_3)R_{p_1 p_2 p_3}(t_2,t_3,t_2)+
$$

\vspace{-1mm}
\begin{equation}
\label{nov9801}
\biggl.\biggl.+
R_{p_1 p_2 p_3}(t_3,t_1,t_1)R_{p_1 p_2 p_3}(t_2,t_3,t_2)\biggr)\biggr)dt_1dt_2dt_3.
\end{equation}

\vspace{4mm} 

Since the integrals on the right-hand side of (\ref{nov9801}) 
exist as Riemann integrals, then they are equal to the 
corresponding Lebesgue integrals. 
Moreover, 

\vspace{1mm}
$$
\lim\limits_{p_1\to\infty}\lim\limits_{p_2\to\infty}\lim\limits_{p_3\to\infty}
R_{p_1 p_2 p_3}(t_1,t_2,t_3)=0\ \ \ \hbox{when}\ \ \ (t_1,t_2,t_3)\in (t, T)^3,
$$

\vspace{4mm}
\noindent
where the left-hand side 
is bounded on $[t, T]^3$
(see (\ref{410})).

Using (\ref{lab11}) and applying three times (we mean an iterated passage to the limit
$\lim\limits_{p_1\to\infty}\lim\limits_{p_2\to\infty}
\lim\limits_{p_3\to\infty}$)
the Lebesgue's Dominated Convergence Theorem in the equality (\ref{nov9801}), 
we obtain

\vspace{1mm}
$$
\lim\limits_{p_1\to\infty}\lim\limits_{p_2\to\infty}\lim\limits_{p_3\to\infty}
{\sf M}\left\{\left(J[R_{p_1p_2p_3}]_{T,t}^{(3)}\right)^2\right\}=0.
$$

\vspace{4mm}
\noindent
The relation (\ref{nov502}) is proved. Theorem 8 is proved.

Developing the approach used in the proof of 
Theorem 8, we can in principle prove 
the following formulas

\vspace{-1mm}
$$
\lim\limits_{p_1\to\infty}
\ldots
\lim\limits_{p_k\to\infty}
{\sf M}\left\{\left(J^{*}[\psi^{(k)}]_{T,t}-
\sum_{j_1=0}^{p_1}\ldots \sum_{j_k=0}^{p_k}
C_{j_k \ldots j_1}\zeta^{(i_1)}_{j_1}
\ldots \zeta^{(i_k)}_{j_k}\right)^{2}\right\}=0,
$$

\vspace{1mm}
$$
\lim\limits_{p_k\to\infty}
\ldots
\lim\limits_{p_1\to\infty}
{\sf M}\left\{\left(J^{*}[\psi^{(k)}]_{T,t}-
\sum_{j_k=0}^{p_k}\ldots \sum_{j_1=0}^{p_1}
C_{j_k \ldots j_1}\zeta^{(i_1)}_{j_1}
\ldots \zeta^{(i_k)}_{j_k}\right)^{2}\right\}=0,
$$

\vspace{5mm}
\noindent
which are correct under the conditions of Theorem 1.

\vspace{5mm}

\section{Expansion of Iterated Stratonovich Stochastic Integrals
of Multiplicity $k.$ The Case $i_1=\ldots=i_k\ne 0$ and
Different Weight Functions $\psi_1(\tau),\ldots,\psi_k(\tau)$}

\vspace{5mm}

In this section, we generalize the approach considered in \cite{17aa}
(also see \cite{10a}, Sect.~2.1.2) 
to the case $i_1=\ldots=i_k\ne 0$ and 
different weight functions $\psi_1(\tau),\ldots,\psi_k(\tau)$ $(k>2).$

Let us formulate the following theorem.

\vspace{2mm}

{\bf Theorem 9} \cite{10a} (Sect.~2.22).\ {\it Suppose that 
$\{\phi_j(x)\}_{j=0}^{\infty}$ is a complete orthonormal system of 
Legendre polynomials or trigonometric functions in the space $L_2([t, T]).$
Moreover, $\psi_1(\tau),\ldots,\psi_k(\tau)$ $(k\ge 2)$ are 
continuously differentiable nonrandom functions on $[t, T]$. Then, 
for the iterated Stratonovich stochastic integral

\vspace{-1mm}
$$
J^{*}[\psi^{(k)}]_{T,t}={\int\limits_t^{*}}^T\psi_k(t_k)\ldots
{\int\limits_t^{*}}^{t_2}\psi_1(t_1)d{\bf f}_{t_1}^{(i_1)}\ldots 
d{\bf f}_{t_k}^{(i_1)}\ \ \ (i_1=1,\ldots,m)
$$

\vspace{3mm}
\noindent
the following equality

\vspace{-1mm}
$$
\lim\limits_{p\to\infty}
{\sf M}\left\{\left(
J^{*}[\psi^{(k)}]_{T,t}-\sum_{j_1=0}^{p}\ldots \sum_{j_k=0}^{p}
C_{j_k \ldots j_1}\zeta_{j_1}^{(i_1)}\ldots \zeta_{j_k}^{(i_1)}\right)^{2n}\right\}=0
$$

\vspace{3mm}
\noindent
is valid, where $n\in\mathbb{N},$

\vspace{-1mm}
$$
C_{j_k \ldots j_1}=\int\limits_t^T\psi_k(t_k)\phi_{j_k}(t_k)
\ldots \int\limits_t^{t_2}\psi_1(t_1)\phi_{j_1}(t_1)dt_1\ldots dt_k
$$

\vspace{3mm}
\noindent
is the Fourier coefficient and

\vspace{-1mm}
$$
\zeta_{j}^{(i_1)}=
\int\limits_t^T \phi_{j}(\tau) d{\bf f}_{\tau}^{(i_1)}\ \ \ (i_1=1,\ldots,m)
$$ 

\vspace{3mm}
\noindent
are independent
standard Gaussian random variables for various 
$j$.}

\vspace{2mm}

{\bf Proof.}\ The case $k=2$ is proved in Theorem~6.
Consider the case $k>2$. First, consider the case $k=3$ in detail.
Define the auxiliary function

\vspace{2mm}
$$
K'(t_1,t_2,t_3)=\frac{1}{6}\left\{
\begin{matrix}
\psi_1(t_1)\psi_2(t_2)\psi_3(t_3),
\ \ t_1\le t_2 \le t_3
\cr\cr
\psi_1(t_1)\psi_2(t_3)\psi_3(t_2),
\ \ t_1\le t_3 \le t_2
\cr\cr
\psi_1(t_2)\psi_2(t_1)\psi_3(t_3),
\ \ t_2\le t_1 \le t_3
\cr\cr
\psi_1(t_2)\psi_2(t_3)\psi_3(t_1),
\ \ t_2\le t_3 \le t_1
\cr\cr
\psi_1(t_3)\psi_2(t_2)\psi_3(t_1),
\ \ t_3\le t_2 \le t_1
\cr\cr
\psi_1(t_3)\psi_2(t_1)\psi_3(t_2),
\ \ t_3\le t_1 \le t_2
\end{matrix}\right.,\ \ \ \ t_1,t_2,t_3\in[t,T].
$$

\vspace{6mm}

Using Lemma~3, Remark~1, and (\ref{30.4}), we obtain w.~p.~1

\vspace{1mm}
$$
J[K']_{T,t}^{(3)}=
\hbox{\vtop{\offinterlineskip\halign{
\hfil#\hfil\cr
{\rm l.i.m.}\cr
$\stackrel{}{{}_{N\to \infty}}$\cr
}} }\sum_{l_3=0}^{N-1}\sum_{l_2=0}^{N-1}
\sum_{l_1=0}^{N-1}
K'(\tau_{l_1},\tau_{l_2},\tau_{l_3})
\Delta{\bf f}_{\tau_{l_1}}^{(i_1)}
\Delta{\bf f}_{\tau_{l_2}}^{(i_1)}
\Delta{\bf f}_{\tau_{l_3}}^{(i_1)}=
$$

\vspace{2mm}
$$
=\hbox{\vtop{\offinterlineskip\halign{
\hfil#\hfil\cr
{\rm l.i.m.}\cr
$\stackrel{}{{}_{N\to \infty}}$\cr
}} }\left(\sum_{l_3=0}^{N-1}\sum_{l_2=0}^{l_3-1}
\sum_{l_1=0}^{l_2-1}
K'(\tau_{l_1},\tau_{l_2},\tau_{l_3})
\Delta{\bf f}_{\tau_{l_1}}^{(i_1)}
\Delta{\bf f}_{\tau_{l_2}}^{(i_1)}
\Delta{\bf f}_{\tau_{l_3}}^{(i_1)}+\right.
$$

\vspace{2mm}
$$
+\sum_{l_3=0}^{N-1}\sum_{l_1=0}^{l_3-1}
\sum_{l_2=0}^{l_1-1}
K'(\tau_{l_1},\tau_{l_2},\tau_{l_3})
\Delta{\bf f}_{\tau_{l_1}}^{(i_1)}
\Delta{\bf f}_{\tau_{l_2}}^{(i_1)}
\Delta{\bf f}_{\tau_{l_3}}^{(i_1)}+
$$

\vspace{2mm}
$$
+\sum_{l_2=0}^{N-1}\sum_{l_1=0}^{l_2-1}
\sum_{l_3=0}^{l_1-1}
K'(\tau_{l_1},\tau_{l_2},\tau_{l_3})
\Delta{\bf f}_{\tau_{l_1}}^{(i_1)}
\Delta{\bf f}_{\tau_{l_2}}^{(i_1)}
\Delta{\bf f}_{\tau_{l_3}}^{(i_1)}+
$$

\vspace{2mm}
$$
+\sum_{l_2=0}^{N-1}\sum_{l_3=0}^{l_2-1}
\sum_{l_1=0}^{l_3-1}
K'(\tau_{l_1},\tau_{l_2},\tau_{l_3})
\Delta{\bf f}_{\tau_{l_1}}^{(i_1)}
\Delta{\bf f}_{\tau_{l_2}}^{(i_1)}
\Delta{\bf f}_{\tau_{l_3}}^{(i_1)}+
$$

\vspace{2mm}
$$
+\sum_{l_1=0}^{N-1}\sum_{l_2=0}^{l_1-1}
\sum_{l_3=0}^{l_2-1}
K'(\tau_{l_1},\tau_{l_2},\tau_{l_3})
\Delta{\bf f}_{\tau_{l_1}}^{(i_1)}
\Delta{\bf f}_{\tau_{l_2}}^{(i_1)}
\Delta{\bf f}_{\tau_{l_3}}^{(i_1)}+
$$

\vspace{2mm}
$$
+\sum_{l_1=0}^{N-1}\sum_{l_3=0}^{l_1-1}
\sum_{l_2=0}^{l_3-1}
K'(\tau_{l_1},\tau_{l_2},\tau_{l_3})
\Delta{\bf f}_{\tau_{l_1}}^{(i_1)}
\Delta{\bf f}_{\tau_{l_2}}^{(i_1)}
\Delta{\bf f}_{\tau_{l_3}}^{(i_1)}+
$$

\vspace{2mm}
$$
+\sum_{l_2=0}^{N-1}\sum_{l_1=0}^{l_2-1}
K'(\tau_{l_1},\tau_{l_2},\tau_{l_1})
\left(\Delta{\bf f}_{\tau_{l_1}}^{(i_1)}\right)^2
\Delta{\bf f}_{\tau_{l_2}}^{(i_1)}+
\sum_{l_3=0}^{N-1}\sum_{l_1=0}^{l_3-1}
K'(\tau_{l_1},\tau_{l_3},\tau_{l_3})
\left(\Delta{\bf f}_{\tau_{l_3}}^{(i_1)}\right)^2
\Delta{\bf f}_{\tau_{l_1}}^{(i_1)}+
$$

\vspace{2mm}
$$
+\sum_{l_1=0}^{N-1}\sum_{l_2=0}^{l_1-1}
K'(\tau_{l_1},\tau_{l_2},\tau_{l_2})
\left(\Delta{\bf f}_{\tau_{l_2}}^{(i_1)}\right)^2
\Delta{\bf f}_{\tau_{l_1}}^{(i_1)}+
\sum_{l_3=0}^{N-1}\sum_{l_2=0}^{l_3-1}
K'(\tau_{l_3},\tau_{l_2},\tau_{l_3})
\left(\Delta{\bf f}_{\tau_{l_3}}^{(i_1)}\right)^2
\Delta{\bf f}_{\tau_{l_2}}^{(i_1)}+
$$

\vspace{2mm}
$$
+\sum_{l_3=0}^{N-1}\sum_{l_2=0}^{l_3-1}
K'(\tau_{l_2},\tau_{l_2},\tau_{l_3})
\left(\Delta{\bf f}_{\tau_{l_2}}^{(i_1)}\right)^2
\Delta{\bf f}_{\tau_{l_3}}^{(i_1)}\left.+
\sum_{l_2=0}^{N-1}\sum_{l_3=0}^{l_2-1}
K'(\tau_{l_2},\tau_{l_2},\tau_{l_3})
\left(\Delta{\bf f}_{\tau_{l_2}}^{(i_1)}\right)^2
\Delta{\bf f}_{\tau_{l_3}}^{(i_1)}\right)=
$$

\vspace{2mm}
$$
=\frac{1}{6}\left(
\int\limits_t^T\psi_3(t_3)
\int\limits_t^{t_3}\psi_2(t_2)
\int\limits_t^{t_2}\psi_1(t_1)
d{\bf f}_{t_1}^{(i_1)}
d{\bf f}_{t_2}^{(i_1)}d{\bf f}_{t_3}^{(i_1)}+\right.
\int\limits_t^T\psi_3(t_2)
\int\limits_t^{t_2}\psi_2(t_1)
\int\limits_t^{t_1}\psi_1(t_3)
d{\bf f}_{t_3}^{(i_1)}
d{\bf f}_{t_1}^{(i_1)}d{\bf f}_{t_2}^{(i_1)}+
$$

\vspace{2mm}
$$
+\int\limits_t^T\psi_3(t_2)
\int\limits_t^{t_2}\psi_2(t_3)
\int\limits_t^{t_3}\psi_1(t_1)
d{\bf f}_{t_1}^{(i_1)}
d{\bf f}_{t_3}^{(i_1)}d{\bf f}_{t_2}^{(i_1)}+
\int\limits_t^T\psi_3(t_3)
\int\limits_t^{t_3}\psi_2(t_1)
\int\limits_t^{t_1}\psi_1(t_2)
d{\bf f}_{t_2}^{(i_1)}
d{\bf f}_{t_1}^{(i_1)}d{\bf f}_{t_3}^{(i_1)}+
$$

\vspace{2mm}
$$
+\int\limits_t^T\psi_3(t_1)
\int\limits_t^{t_1}\psi_2(t_2)
\int\limits_t^{t_2}\psi_1(t_3)
d{\bf f}_{t_3}^{(i_1)}
d{\bf f}_{t_2}^{(i_1)}d{\bf f}_{t_1}^{(i_1)}+
\int\limits_t^T\psi_3(t_1)
\int\limits_t^{t_1}\psi_2(t_3)
\int\limits_t^{t_3}\psi_1(t_2)
d{\bf f}_{t_2}^{(i_1)}
d{\bf f}_{t_3}^{(i_1)}d{\bf f}_{t_1}^{(i_1)}+
$$

\vspace{2mm}
$$
+
\int\limits_t^T\psi_3(t_2)
\int\limits_t^{t_2}\psi_2(t_1)\psi_1(t_1)
dt_1d{\bf f}_{t_2}^{(i_1)}+
\int\limits_t^T\psi_3(t_1)
\int\limits_t^{t_1}\psi_2(t_2)\psi_1(t_2)
dt_2d{\bf f}_{t_1}^{(i_1)}+
$$

\vspace{2mm}
$$
+
\int\limits_t^T\psi_3(t_3)
\int\limits_t^{t_3}\psi_2(t_1)\psi_1(t_1)
dt_1d{\bf f}_{t_3}^{(i_1)}+
\int\limits_t^T\psi_3(t_3)\psi_2(t_3)
\int\limits_t^{t_3}\psi_1(t_1)
d{\bf f}_{t_1}^{(i_1)}
dt_3+
$$

\vspace{2mm}
$$
\left.+\int\limits_t^T\psi_3(t_3)\psi_2(t_3)
\int\limits_t^{t_3}\psi_1(t_2)
d{\bf f}_{t_2}^{(i_1)}
dt_3+
\int\limits_t^T\psi_3(t_2)\psi_2(t_2)
\int\limits_t^{t_2}\psi_1(t_3)
d{\bf f}_{t_3}^{(i_1)}
dt_2\right)=
$$

\vspace{2mm}
$$
=\int\limits_t^T\psi_3(t_3)
\int\limits_t^{t_3}\psi_2(t_2)
\int\limits_t^{t_2}\psi_1(t_1)
d{\bf f}_{t_1}^{(i_1)}
d{\bf f}_{t_2}^{(i_1)}d{\bf f}_{t_3}^{(i_1)}+
$$

\vspace{2mm}
$$
+\frac{1}{2}
\int\limits_t^T\psi_3(t_3)
\int\limits_t^{t_3}\psi_2(t_1)\psi_1(t_1)
dt_1d{\bf f}_{t_3}^{(i_1)}+
\frac{1}{2}
\int\limits_t^T\psi_3(t_3)\psi_2(t_3)
\int\limits_t^{t_3}\psi_1(t_1)
d{\bf f}_{t_1}^{(i_1)}
dt_3=
$$

\vspace{2mm}
\begin{equation}
\label{zikoxxx1}
=
{\int\limits_t^{*}}^T\psi_3(t_3)
{\int\limits_t^{*}}^{t_3}\psi_2(t_2)
{\int\limits_t^{*}}^{t_2}\psi_1(t_1)
d{\bf f}_{t_1}^{(i_1)}
d{\bf f}_{t_2}^{(i_1)}d{\bf f}_{t_3}^{(i_1)}\stackrel{\sf def}{=}J^{*}[\psi^{(3)}]_{T,t},
\end{equation}

\vspace{5mm}
\noindent
where the multiple stochastic integral $J[K']_{T,t}^{(3)}$
is defined by (\ref{30.34}) and 
$\left\{\tau_{j}\right\}_{j=0}^{N}$ is a partition of
$[t,T],$ which satisfies the condition {\rm (\ref{1111}).

\vspace{2mm}

{\it For each $\delta>0$ let us call the exact upper edge 
of difference $\left|f({\bf t}')-f({\bf t}'')\right|$ 
in the set
of all points ${\bf t}'$, ${\bf t}''$ which 
belong 
to the domain $D$ 
as the module of 
continuity of the function
$f({\bf t})$ $({\bf t}=(t_1,\ldots,t_k))$ in the 
$k$-dimentional domain
$D$ $(k\ge 1)$ 
if the distance between ${\bf t}',{\bf t}''$
satisfies the condition
$\rho\left({\bf t}',{\bf t}''\right)<\delta.$}

\vspace{2mm}

{\it We will say that the function 
of $k$ $(k\ge 1)$ variables  
$f({\bf t})$ $({\bf t}=(t_1,\ldots,t_k))$
belongs 
to the H\"{o}lder class with 
the parameter $\alpha\in (0, 1]$ $(f({\bf t})\in C^{\alpha}(D))$ 
in the domain $D$ 
if the module of 
continuity of the function
$f({\bf t})$ $({\bf t}=(t_1,\ldots,t_k))$ 
in the domain $D$ have the orders $o(\delta^{\alpha})$ $(\alpha \in (0, 1))$
and $O(\delta)$ $(\alpha=1)$.}

\vspace{2mm}

In 1967, Zhizhiashvili L.V.
proved
that the rectangular sums of multiple trigonometric Fourier series 
of the function of $k$ variables  
in the hypercube $[t,T]^k$ 
converge uniformly to this function in the hypercube $[t,T]^k$ if 
the function
belongs
to $C^{\alpha}([t,T]^k),$ $\alpha>0$ (definition
of the H\"{o}lder class with any parameter $\alpha>0$ can be found in 
the well known mathematical analysis tutorials \cite{IP}).

More precisely, the following statement is correct.

\vspace{2mm}

{\bf Theorem 10} \cite{IP}. {\it If the function
$f(x_1,\ldots,x_n)$ is periodic with period $2\pi$ with respect to each
variable and belongs in $\mathbb{R}^n$ to the H\"{o}lder class 
$C^{\alpha}(\mathbb{R}^n)$
for any $\alpha>0,$ then the rectangular partial
sums of multiple trigonometric Fourier series of the function
$f(x_1,\ldots,x_n)$ converge to this function uniformly in 
$\mathbb{R}^n$.}

\vspace{2mm}

In \cite{17aa} (also see \cite{10a}, Sect.~2.1.2)
it was shown that the following function 

\vspace{1mm} 
$$
K'(t_1,t_2)=\left\{
\begin{matrix}
\psi_1(t_1)\psi_2(t_2),\ \ t_1\le t_2\cr\cr\cr
\psi_1(t_2)\psi_2(t_1),\ \ t_2\le t_1
\end{matrix}
\right.,\ \ \ \ t_1,t_2\in[t,T]
$$

\vspace{5mm}
\noindent
belongs to the class $C^{1}([t,T]^2).$
Moreover, the following Fourier--Legendre expansion

\vspace{1mm}
$$
K'(t_1,t_2)=
\lim_{p\to\infty}
\sum_{j_1=0}^{p}\sum_{j_2=0}^{p}
\int\limits_t^T\int\limits_t^TK'(t_1,t_2)
\phi_{j_1}(t_1)\phi_{j_2}(t_2)dt_1 dt_2\cdot
\phi_{j_1}(t_1)\phi_{j_2}(t_2)=
$$

\vspace{2mm}
\begin{equation}
\label{334.ye}
=
\lim_{p\to\infty}\sum_{j_1=0}^{p}\sum_{j_2=0}^{p}\left(C_{j_2j_1}+
C_{j_1j_2}\right)
\phi_{j_1}(t_1)\phi_{j_2}(t_2)
\end{equation}

\vspace{4mm}
\noindent
is valid for $(t_1, t_2)\in (t, T)^2.$

Using Theorem~10 for $n=3$ 
and generalizing the Fourier--Legendre expansion (\ref{334.ye})
for the function $K'(t_1,t_2,t_3)$, we obtain 

\vspace{1mm}
$$
K'(t_1,t_2,t_3)=\lim_{p\to\infty}
\sum_{j_1=0}^{p}\sum_{j_2=0}^{p}
\sum_{j_3=0}^{p}\frac{1}{6}\biggl(C_{j_3j_2j_1}+
C_{j_3j_1j_2}+C_{j_2j_1j_3}+\biggr.
$$

\vspace{2mm}
\begin{equation}
\label{ziko9006}
\biggl.+C_{j_2j_3j_1}+C_{j_1j_2j_3}+C_{j_1j_3j_2}\biggr)
\phi_{j_1}(t_1)\phi_{j_2}(t_2)\phi_{j_3}(t_3),
\end{equation}

\vspace{3mm}
\noindent
where the multiple Fourier series (\ref{ziko9006}) converges to the function 
$K'(t_1,t_2,t_3)$ in $(t,T)^3$
and the partial sums of 
the series (\ref{ziko9006}) have an integrable majorant on $[t, T]^3$
that does not depend
on $p.$
For the trigonomertic case, the above statement follows from 
Theorem~10 (the proof that the function
$K'(t_1, t_2,t_3)$ belongs to the H\"{o}lder class 
with parameter $1$ in $[t, T]^3$ is omitted and 
can be carried out in the same way as for the function
$K'(t_1, t_2)$ in the two-dimensional case 
\cite{17aa}
(also see \cite{10a}, Sect.~2.1.2)).
The proof of generalization of the Fourier--Legendre expansion
(\ref{334.ye}) to the three-dimensional case (see (\ref{ziko9006}))
is omitted. The proof that the partial sums of 
the series (\ref{ziko9006}) have an integrable majorant on $[t, T]^3$
is also omitted.

Denote
$$
R'_{ppp}(t_1,t_2,t_3)=K'(t_1,t_2,t_3)-
\sum_{j_1=0}^{p}\sum_{j_2=0}^{p}
\sum_{j_3=0}^{p}\frac{1}{6}\biggl(C_{j_3j_2j_1}+
C_{j_3j_1j_2}+C_{j_2j_1j_3}+\biggr.
$$

\vspace{2mm}
$$
\biggl.+C_{j_2j_3j_1}+C_{j_1j_2j_3}+C_{j_1j_3j_2}\biggr)
\phi_{j_1}(t_1)\phi_{j_2}(t_2)\phi_{j_3}(t_3).
$$

\vspace{4mm}

Using Lemma~5 and (\ref{zikoxxx1}), we get w.~p.~1

\vspace{1mm}
$$
J^{*}[\psi^{(3)}]_{T,t}=J[K']_{T,t}^{(3)}=
\sum_{j_1=0}^{p}\sum_{j_2=0}^{p}
\sum_{j_3=0}^{p}\frac{1}{6}\biggl(C_{j_3j_2j_1}+
C_{j_3j_1j_2}+C_{j_2j_1j_3}+\biggr.
$$

\vspace{2mm}
$$
\biggl.+C_{j_2j_3j_1}+C_{j_1j_2j_3}+C_{j_1j_3j_2}\biggr)
\zeta_{j_1}^{(i_1)}\zeta_{j_2}^{(i_1)}\zeta_{j_3}^{(i_1)}+
J[R'_{ppp}]_{T,t}^{(3)}=
$$

\vspace{2mm}
$$
=\sum_{j_1=0}^{p}\sum_{j_2=0}^{p}
\sum_{j_3=0}^{p}C_{j_3j_2j_1}
\zeta_{j_1}^{(i_1)}\zeta_{j_2}^{(i_1)}\zeta_{j_3}^{(i_1)}+
J[R'_{ppp}]_{T,t}^{(3)}.
$$

\vspace{4mm}

Then 
$$
{\sf M}\left\{\left(J[R'_{ppp}]_{T,t}^{(3)}\right)^{2n}\right\}=
{\sf M}\left\{\left(J^{*}[\psi^{(3)}]_{T,t}-
\sum_{j_1=0}^{p}\sum_{j_2=0}^{p}
\sum_{j_3=0}^{p}C_{j_3j_2j_1}
\zeta_{j_1}^{(i_1)}\zeta_{j_2}^{(i_1)}\zeta_{j_3}^{(i_1)}\right)^{2n}\right\},
$$

\vspace{3mm}
\noindent
where $n\in \mathbb{N}.$

Applying (we mean here the passage to the limit
$\lim\limits_{p\to\infty}$)
the Lebesgue's Dominated Convergence Theorem to the integrals
on the right-hand side of (\ref{udar1}) for $k=3$ and 
$R'_{ppp}(t_1,t_2,t_3)$ instead of $R_{p_1p_2p_3}(t_1,t_2,t_3)$,
we obtain

\vspace{-1mm}
$$
\lim\limits_{p\to\infty}
{\sf M}\left\{\left(J[R'_{ppp}]_{T,t}^{(3)}\right)^{2n}\right\}=0.
$$

\vspace{3mm}
\noindent
Theorems 9 is proved for the case $k=3.$ 

To prove Theorem 9 for the case $k>3$, consider the auxiliary function

\vspace{1mm}
\begin{equation}
\label{zikoyyy1}
K'(t_1,\ldots,t_k)=\frac{1}{k!}\left\{
\begin{matrix}
\psi_1(t_1)\ldots\psi_k(t_k),
\ \ t_1\le \ldots \le t_k
\cr\cr
\ldots 
\cr\cr
\psi_1(t_{g_1})\ldots\psi_k(t_{g_k}),
\ \ t_{g_1}\le \ldots \le t_{g_k}
\cr\cr
\ldots
\cr\cr
\psi_1(t_k)\ldots\psi_k(t_1),
\ \ t_k\le \ldots \le t_1
\end{matrix}\right.,\ \ \ \ t_1,\ldots,t_k\in[t,T],
\end{equation}

\vspace{5mm}
\noindent
where $\{g_1,\ldots,g_k\}=\{1,\ldots,k\}$
and we take into account all
possible permutations $(g_1,\ldots,g_k)$ on the right-hand side of the formula
(\ref{zikoyyy1}).

Further, we have w.~p.~1

\vspace{-1mm}
\begin{equation}
\label{zikozzz1a}
J[K']_{T,t}^{(k)}=
J[\psi^{(k)}]_{T,t}+
\sum_{r=1}^{\left[k/2\right]}\frac{1}{2^r}
\sum_{(s_r,\ldots,s_1)\in {\rm A}_{k,r}}
J[\psi^{(k)}]_{T,t}^{s_r,\ldots,s_1},
\end{equation}

\vspace{4mm}
\noindent
where the function $K'(t_1,\ldots,t_k)$
is defined by (\ref{zikoyyy1}); another
notations are the same as in (\ref{30.1}) and 
Lemma~2 ($i_1=\ldots=i_k\ne 0$ in (\ref{30.1})).

From (\ref{zikozzz1a}) and Lemma~2 we obtain w.~p.~1

\begin{equation}
\label{zikozzz1}
J^{*}[\psi^{(k)}]_{T,t}=J[K']_{T,t}^{(k)},
\end{equation}

\vspace{4mm}
\noindent
where $i_1=\ldots=i_k\ne 0.$

Generalizing the above reasoning to the case $k>3$ and taking into account
(\ref{zikozzz1}), we get

\vspace{2mm}
$$
J^{*}[\psi^{(k)}]_{T,t}=
\sum_{j_1=0}^{p}\ldots
\sum_{j_k=0}^{p}\frac{1}{k!}\left(\sum\limits_{(j_1,\ldots,j_k)}
C_{j_k\ldots j_1}\right)
\zeta_{j_1}^{(i_1)}\ldots\zeta_{j_k}^{(i_1)}+
J[R'_{p\ldots p}]_{T,t}^{(k)}=
$$

\vspace{3mm}
$$
=\sum_{j_1=0}^{p}\ldots
\sum_{j_k=0}^{p}C_{j_k\ldots j_1}
\zeta_{j_1}^{(i_1)}\ldots \zeta_{j_k}^{(i_1)}+
J[R'_{p\ldots p}]_{T,t}^{(k)},
$$

\vspace{5mm}
\noindent
where

\vspace{-1mm}
$$
R'_{p\ldots p}(t_1,\ldots,t_k)\stackrel{\sf def}{=}K'(t_1,\ldots,t_k)-
$$

\vspace{2mm}
$$
-
\sum_{j_1=0}^{p}\ldots
\sum_{j_k=0}^{p}\frac{1}{k!}\left(\sum\limits_{(j_1,\ldots,j_k)}
C_{j_k\ldots j_1}\right)
\phi_{j_1}(t_1)\ldots\phi_{j_k}(t_k),
$$

\vspace{5mm}
\noindent
the expression

\vspace{-1mm}
$$
\sum\limits_{(j_1, \ldots, j_k)}
$$

\vspace{4mm}
\noindent
means the sum with respect to 
all possible  
permutations $(j_1, \ldots, j_k)$.

Further,

\vspace{-1mm}
$$
{\sf M}\left\{\left(J[R'_{p\ldots p}]_{T,t}^{(k)}\right)^{2n}\right\}=
{\sf M}\left\{\left(J^{*}[\psi^{(k)}]_{T,t}-
\sum_{j_1=0}^{p}\ldots
\sum_{j_k=0}^{p}C_{j_k\ldots j_1}
\zeta_{j_1}^{(i_1)}\ldots \zeta_{j_k}^{(i_1)}\right)^{2n}\right\},
$$

\vspace{4mm}
\noindent
where $n\in \mathbb{N}.$

Applying (we mean here the passage to the limit
$\lim\limits_{p\to\infty}$)
the Lebesgue's Dominated Convergence Theorem to the integrals
on the right-hand side of (\ref{udar1}) for the function
$R'_{p\ldots p}(t_1,\ldots,t_k)$ instead of the function $R_{p_1\ldots,p_k}(t_1,\ldots,t_k)$,
we obtain

$$
\lim\limits_{p\to\infty}
{\sf M}\left\{\left(J[R'_{p\ldots p}]_{T,t}^{(k)}\right)^{2n}\right\}=0.
$$

\vspace{4mm}
\noindent
Theorems 9 is proved.

\vspace{5mm}

\section{Recent Results on Expansion of Iterated Ito and Stratonovich
Stochastic Integrals}

\vspace{5mm}

Using (\ref{leto5008}), we can write (\ref{tyyy}) as

\vspace{1mm}

$$
J[\psi^{(k)}]_{T,t}=
\hbox{\vtop{\offinterlineskip\halign{
\hfil#\hfil\cr
{\rm l.i.m.}\cr
$\stackrel{}{{}_{p_1,\ldots,p_k\to \infty}}$\cr
}} }
\sum\limits_{j_1=0}^{p_1}\ldots
\sum\limits_{j_k=0}^{p_k}
C_{j_k\ldots j_1}\Biggl(
\prod_{l=1}^k\zeta_{j_l}^{(i_l)}+\sum\limits_{r=1}^{[k/2]}
(-1)^r \times
\Biggr.
$$

\vspace{2mm}
\begin{equation}
\label{leto6000hh}
\times
\sum_{\stackrel{(\{\{g_1, g_2\}, \ldots, 
\{g_{2r-1}, g_{2r}\}\}, \{q_1, \ldots, q_{k-2r}\})}
{{}_{\{g_1, g_2, \ldots, 
g_{2r-1}, g_{2r}, q_1, \ldots, q_{k-2r}\}=\{1, 2, \ldots, k\}}}}
\prod\limits_{s=1}^r
{\bf 1}_{\{i_{g_{{}_{2s-1}}}=~i_{g_{{}_{2s}}}\ne 0\}}
\Biggl.{\bf 1}_{\{j_{g_{{}_{2s-1}}}=~j_{g_{{}_{2s}}}\}}
\prod_{l=1}^{k-2r}\zeta_{j_{q_l}}^{(i_{q_l})}\Biggr),
\end{equation}

\vspace{5mm}
\noindent
where $[x]$ is an integer part of a real number $x,$
$\prod\limits_{\emptyset}
\stackrel{\sf def}{=}1,$ $\sum\limits_{\emptyset}
\stackrel{\sf def}{=}0;$
another notations are the same as in Theorem 4.

\vspace{2mm}

In particular, from (\ref{leto6000hh}) for $k=5$ we obtain

\vspace{3mm}

$$
J[\psi^{(5)}]_{T,t}=
\hbox{\vtop{\offinterlineskip\halign{
\hfil#\hfil\cr
{\rm l.i.m.}\cr
$\stackrel{}{{}_{p_1,\ldots,p_5\to \infty}}$\cr
}} }\sum_{j_1=0}^{p_1}\ldots\sum_{j_5=0}^{p_5}
C_{j_5\ldots j_1}\Biggl(
\prod_{l=1}^5\zeta_{j_l}^{(i_l)}-\Biggr.
$$

\vspace{2mm}
$$
-
\sum\limits_{\stackrel{(\{g_1, g_2\}, \{q_1, q_{2}, q_3\})}
{{}_{\{g_1, g_2, q_{1}, q_{2}, q_3\}=\{1, 2, 3, 4, 5\}}}}
{\bf 1}_{\{i_{g_{{}_{1}}}=~i_{g_{{}_{2}}}\ne 0\}}
{\bf 1}_{\{j_{g_{{}_{1}}}=~j_{g_{{}_{2}}}\}}
\prod_{l=1}^{3}\zeta_{j_{q_l}}^{(i_{q_l})}+
$$

\vspace{2mm}
$$
+
\sum_{\stackrel{(\{\{g_1, g_2\}, 
\{g_{3}, g_{4}\}\}, \{q_1\})}
{{}_{\{g_1, g_2, g_{3}, g_{4}, q_1\}=\{1, 2, 3, 4, 5\}}}}
{\bf 1}_{\{i_{g_{{}_{1}}}=~i_{g_{{}_{2}}}\ne 0\}}
{\bf 1}_{\{j_{g_{{}_{1}}}=~j_{g_{{}_{2}}}\}}
\Biggl.{\bf 1}_{\{i_{g_{{}_{3}}}=~i_{g_{{}_{4}}}\ne 0\}}
{\bf 1}_{\{j_{g_{{}_{3}}}=~j_{g_{{}_{4}}}\}}
\zeta_{j_{q_1}}^{(i_{q_1})}\Biggr).
$$

\vspace{7mm}
\noindent
The last equality obviously agrees with
(\ref{a5}).

Let us consider the generalization of Theorem 4 for the case
of an arbitrary complete orthonormal systems  
of functions in the space $L_2([t,T])$ 
and $\psi_1(\tau),\ldots,\psi_k(\tau)\in L_2([t, T]).$

\vspace{2mm}

{\bf Theorem~11}\ \cite{10a} (Sect.~1.11), \cite{11} (Sect.~15).
{\it Suppose that
$\psi_1(\tau),\ldots,\psi_k(\tau)\in L_2([t, T])$ and
$\{\phi_j(x)\}_{j=0}^{\infty}$ is an arbitrary complete orthonormal system  
of functions in the space $L_2([t,T]).$
Then the following expansion

$$
J[\psi^{(k)}]_{T,t}=
\hbox{\vtop{\offinterlineskip\halign{
\hfil#\hfil\cr
{\rm l.i.m.}\cr
$\stackrel{}{{}_{p_1,\ldots,p_k\to \infty}}$\cr
}} }
\sum\limits_{j_1=0}^{p_1}\ldots
\sum\limits_{j_k=0}^{p_k}
C_{j_k\ldots j_1}\Biggl(
\prod_{l=1}^k\zeta_{j_l}^{(i_l)}+\sum\limits_{r=1}^{[k/2]}
(-1)^r \times
\Biggr.
$$

\vspace{3mm}
$$
\times
\sum_{\stackrel{(\{\{g_1, g_2\}, \ldots, 
\{g_{2r-1}, g_{2r}\}\}, \{q_1, \ldots, q_{k-2r}\})}
{{}_{\{g_1, g_2, \ldots, 
g_{2r-1}, g_{2r}, q_1, \ldots, q_{k-2r}\}=\{1, 2, \ldots, k\}}}}
\prod\limits_{s=1}^r
{\bf 1}_{\{i_{g_{{}_{2s-1}}}=~i_{g_{{}_{2s}}}\ne 0\}}
\Biggl.{\bf 1}_{\{j_{g_{{}_{2s-1}}}=~j_{g_{{}_{2s}}}\}}
\prod_{l=1}^{k-2r}\zeta_{j_{q_l}}^{(i_{q_l})}\Biggr)
$$

\vspace{6mm}
\noindent
con\-verg\-ing in the mean-square sense is valid,
where $[x]$ is an integer part of a real number $x,$
$\prod\limits_{\emptyset}
\stackrel{\sf def}{=}1,$ $\sum\limits_{\emptyset}
\stackrel{\sf def}{=}0;$
another notations are the same as in Theorem~{\rm 4}.}

\vspace{2mm}

It should be noted that an analogue of Theorem 11 was considered 
in \cite{Rybakov1000}. 
Note that we use another notations 
\cite{10a} (Sect.~1.11), \cite{11} (Sect.~15)
in comparison with \cite{Rybakov1000}.
Moreover, the proof of an analogue of Theorem 11
from \cite{Rybakov1000} is different from the proof given in 
\cite{10a} (Sect.~1.11), \cite{11} (Sect.~15).

Recently, a new approach to the expansion and mean-square 
approximation of iterated Stratonovich stochastic integrals has been obtained
\cite{10a} (Sect.~2.10--2.16), \cite{12} (Sect.~13--19), 
\cite{15a} (Sect.~5--11), \cite{hhh111hhh} (Sect.~7--13), \cite{new-art-1-xxy}
(Sect.~4--9).
Let us formulate four theorems that were obtained using this approach.

\vspace{2mm}

{\bf Theorem 12}\ \cite{10a}, \cite{12}, \cite{15a}, \cite{hhh111hhh}, \cite{new-art-1-xxy}.\
{\it Suppose 
that $\{\phi_j(x)\}_{j=0}^{\infty}$ is a complete orthonormal system of 
Legendre polynomials or trigonometric functions in the space $L_2([t, T]).$
Furthermore, let $\psi_1(\tau), \psi_2(\tau),$ $\psi_3(\tau)$ are continuously dif\-ferentiable 
nonrandom functions on $[t, T].$ 
Then, for the 
iterated Stra\-to\-no\-vich stochastic integral of third multiplicity

$$
J^{*}[\psi^{(3)}]_{T,t}={\int\limits_t^{*}}^T\psi_3(t_3)
{\int\limits_t^{*}}^{t_3}\psi_2(t_2)
{\int\limits_t^{*}}^{t_2}\psi_1(t_1)
d{\bf w}_{t_1}^{(i_1)}
d{\bf w}_{t_2}^{(i_2)}d{\bf w}_{t_3}^{(i_3)}\ \ \ (i_1,i_2,i_3=0,1,\ldots,m)
$$

\vspace{4mm}
\noindent
the following 
relations

\vspace{-1mm}
\begin{equation}
\label{fin1}
J^{*}[\psi^{(3)}]_{T,t}
=\hbox{\vtop{\offinterlineskip\halign{
\hfil#\hfil\cr
{\rm l.i.m.}\cr
$\stackrel{}{{}_{p\to \infty}}$\cr
}} }
\sum\limits_{j_1, j_2, j_3=0}^{p}
C_{j_3 j_2 j_1}\zeta_{j_1}^{(i_1)}\zeta_{j_2}^{(i_2)}\zeta_{j_3}^{(i_3)},
\end{equation}

\vspace{3mm}
\begin{equation}
\label{fin2}
{\sf M}\left\{\left(
J^{*}[\psi^{(3)}]_{T,t}-
\sum\limits_{j_1, j_2, j_3=0}^{p}
C_{j_3 j_2 j_1}\zeta_{j_1}^{(i_1)}\zeta_{j_2}^{(i_2)}\zeta_{j_3}^{(i_3)}\right)^2\right\}
\le \frac{C}{p}
\end{equation}

\vspace{5mm}
\noindent
are fulfilled, where $i_1, i_2, i_3=0,1,\ldots,m$ in {\rm (\ref{fin1})} and 
$i_1, i_2, i_3=1,\ldots,m$ in {\rm (\ref{fin2})},
constant $C$ is independent of $p,$

$$
C_{j_3 j_2 j_1}=\int\limits_t^T\psi_3(t_3)\phi_{j_3}(t_3)
\int\limits_t^{t_3}\psi_2(t_2)\phi_{j_2}(t_2)
\int\limits_t^{t_2}\psi_1(t_1)\phi_{j_1}(t_1)dt_1dt_2dt_3
$$

\vspace{4mm}
\noindent
and
$$
\zeta_{j}^{(i)}=
\int\limits_t^T \phi_{j}(\tau) d{\bf f}_{\tau}^{(i)}
$$ 

\vspace{2mm}
\noindent
are independent standard Gaussian random variables for various 
$i$ or $j$ {\rm (}in the case when $i\ne 0${\rm );} 
another notations are the same as in Theorems~{\rm 4, 11}.}

\vspace{2mm}

{\bf Theorem 13}\ \cite{10a}, \cite{12}, \cite{15a}, \cite{hhh111hhh}, \cite{new-art-1-xxy}.\ 
{\it Let
$\{\phi_j(x)\}_{j=0}^{\infty}$ be a complete orthonormal system of 
Legendre polynomials or trigonometric functions in the space $L_2([t, T]).$
Furthermore, let $\psi_1(\tau), \ldots, \psi_4(\tau)$ be continuously dif\-ferentiable 
nonrandom functions on $[t, T].$ 
Then, for the 
iterated Stra\-to\-no\-vich stochastic integral of fourth multiplicity

\begin{equation}
\label{fin0}
J^{*}[\psi^{(4)}]_{T,t}={\int\limits_t^{*}}^T\psi_4(t_4)
{\int\limits_t^{*}}^{t_4}\psi_3(t_3)
{\int\limits_t^{*}}^{t_3}\psi_2(t_2)
{\int\limits_t^{*}}^{t_2}\psi_1(t_1)
d{\bf w}_{t_1}^{(i_1)}
d{\bf w}_{t_2}^{(i_2)}d{\bf w}_{t_3}^{(i_3)}d{\bf w}_{t_4}^{(i_4)}
\end{equation}

\vspace{4mm}
\noindent
the following 
relations

\begin{equation}
\label{fin3}
J^{*}[\psi^{(4)}]_{T,t}
=\hbox{\vtop{\offinterlineskip\halign{
\hfil#\hfil\cr
{\rm l.i.m.}\cr
$\stackrel{}{{}_{p\to \infty}}$\cr
}} }
\sum\limits_{j_1, j_2, j_3,j_4=0}^{p}
C_{j_4j_3 j_2 j_1}\zeta_{j_1}^{(i_1)}\zeta_{j_2}^{(i_2)}\zeta_{j_3}^{(i_3)}\zeta_{j_4}^{(i_4)},
\end{equation}

\vspace{3mm}

\begin{equation}
\label{fin4}
{\sf M}\left\{\left(
J^{*}[\psi^{(4)}]_{T,t}-
\sum\limits_{j_1, j_2, j_3, j_4=0}^{p}
C_{j_4 j_3 j_2 j_1}\zeta_{j_1}^{(i_1)}\zeta_{j_2}^{(i_2)}\zeta_{j_3}^{(i_3)}
\zeta_{j_4}^{(i_4)}
\right)^2\right\}
\le \frac{C}{p^{1-\varepsilon}}
\end{equation}

\vspace{5mm}
\noindent
are fulfilled, where $i_1, \ldots , i_4=0,1,\ldots,m$ in {\rm (\ref{fin0}),} {\rm (\ref{fin3})} 
and $i_1, \ldots, i_4=1,\ldots,m$ in {\rm (\ref{fin4}),}
constant $C$ does not depend on $p,$
$\varepsilon$ is an arbitrary
small positive real number 
for the case of complete orthonormal system of 
Legendre polynomials in the space $L_2([t, T])$
and $\varepsilon=0$ for the case of
complete orthonormal system of 
trigonometric functions in the space $L_2([t, T]),$

$$
C_{j_4 j_3 j_2 j_1}=
$$

$$
=
\int\limits_t^T\psi_4(t_4)\phi_{j_4}(t_4)
\int\limits_t^{t_4}\psi_3(t_3)\phi_{j_3}(t_3)
\int\limits_t^{t_3}\psi_2(t_2)\phi_{j_2}(t_2)
\int\limits_t^{t_2}\psi_1(t_1)\phi_{j_1}(t_1)dt_1dt_2dt_3dt_4;
$$

\vspace{4mm}
\noindent
another notations are the same as in Theorem~{\rm 12}.}

\vspace{2mm}

{\bf Theorem 14}\ \cite{10a}, \cite{12}, \cite{15a}, \cite{hhh111hhh}, \cite{new-art-1-xxy}.\
{\it Assume 
that $\{\phi_j(x)\}_{j=0}^{\infty}$ is a complete orthonormal system of 
Legendre polynomials or trigonometric functions in the space $L_2([t, T])$
and $\psi_1(\tau), \ldots, \psi_5(\tau)$ are continuously dif\-ferentiable 
nonrandom functions on $[t, T].$ 
Then, for the 
iterated Stra\-to\-no\-vich stochastic integral of fifth multiplicity

\begin{equation}
\label{fin7}
J^{*}[\psi^{(5)}]_{T,t}={\int\limits_t^{*}}^T\psi_5(t_5)
\ldots
{\int\limits_t^{*}}^{t_2}\psi_1(t_1)
d{\bf w}_{t_1}^{(i_1)}
\ldots d{\bf w}_{t_5}^{(i_5)}
\end{equation}

\vspace{4mm}
\noindent
the following 
relations

\begin{equation}
\label{fin8}
J^{*}[\psi^{(5)}]_{T,t}
=\hbox{\vtop{\offinterlineskip\halign{
\hfil#\hfil\cr
{\rm l.i.m.}\cr
$\stackrel{}{{}_{p\to \infty}}$\cr
}} }
\sum\limits_{j_1,\ldots,j_5=0}^{p}
C_{j_5 \ldots j_1}\zeta_{j_1}^{(i_1)}\ldots \zeta_{j_5}^{(i_5)},
\end{equation}

\vspace{3mm}

\begin{equation}
\label{fin9}
{\sf M}\left\{\left(
J^{*}[\psi^{(5)}]_{T,t}-
\sum\limits_{j_1, \ldots, j_5=0}^{p}
C_{j_5 \ldots j_1}\zeta_{j_1}^{(i_1)}\ldots
\zeta_{j_5}^{(i_5)}
\right)^2\right\}
\le \frac{C}{p^{1-\varepsilon}}
\end{equation}

\vspace{5mm}
\noindent
are fulfilled, where $i_1, \ldots , i_5=0,1,\ldots,m$ in {\rm (\ref{fin7}),} {\rm (\ref{fin8})} 
and $i_1, \ldots, i_5=1,\ldots,m$ in {\rm (\ref{fin9}),}
constant $C$ is independent of $p,$
$\varepsilon$ is an arbitrary
small positive real number 
for the case of complete orthonormal system of 
Legendre polynomials in the space $L_2([t, T])$
and $\varepsilon=0$ for the case of
complete orthonormal system of 
trigonometric functions in the space $L_2([t, T]),$

\vspace{-1mm}
$$
C_{j_5 \ldots j_1}=
\int\limits_t^T\psi_5(t_5)\phi_{j_5}(t_5)\ldots
\int\limits_t^{t_2}\psi_1(t_1)\phi_{j_1}(t_1)dt_1\ldots dt_5;
$$

\vspace{3mm}
\noindent
another notations are the same as in Theorems~{\rm 12, 13}.}

\vspace{2mm}

{\bf Theorem 15}\ \cite{10a}, \cite{12}, \cite{15a}, \cite{hhh111hhh}.\
{\it Suppose that 
$\{\phi_j(x)\}_{j=0}^{\infty}$ is a complete orthonormal system of 
Legendre polynomials or trigonometric functions in the space $L_2([t, T]).$
Then, for the 
iterated Stratonovich stochastic integral of sixth multiplicity

\begin{equation}
\label{after10001qu1}
J_{T,t}^{*(i_1\ldots i_6)}={\int\limits_t^{*}}^T
\ldots
{\int\limits_t^{*}}^{t_2}
d{\bf w}_{t_1}^{(i_1)}
\ldots d{\bf w}_{t_6}^{(i_6)}
\end{equation}

\vspace{3mm}
\noindent
the following 
expansion 

\vspace{-1mm}
$$
J_{T,t}^{*(i_1\ldots i_6)}
=\hbox{\vtop{\offinterlineskip\halign{
\hfil#\hfil\cr
{\rm l.i.m.}\cr
$\stackrel{}{{}_{p\to \infty}}$\cr
}} }
\sum\limits_{j_1, \ldots, j_6=0}^{p}
C_{j_6 \ldots j_1}\zeta_{j_1}^{(i_1)}\ldots
\zeta_{j_6}^{(i_6)}
$$

\vspace{4mm}
\noindent
that converges in the mean-square sense is valid, where
$i_1, \ldots, i_6=0, 1,\ldots,m,$

$$
C_{j_6 \ldots j_1}=
\int\limits_t^T\phi_{j_6}(t_6)\ldots
\int\limits_t^{t_2}\phi_{j_1}(t_1)dt_1\ldots dt_6;
$$

\vspace{3mm}
\noindent
another notations are the same as in Theorems~{\rm 12--14}.}

The results of Theorems~12--15 were developed in 
\cite{10a} (Chapter~2).
In particular, analogues of Theorem~15 for iterated Stratonovich stochastic
integrals of multiplicities 7 and 8 were obtained in \cite{10a} (Sect.~2.36, 2.37).
In addition, the variants of Theorems 12--15
were obtained
for the case when $\{\phi_j(x)\}_{j=0}^{\infty}$ is an arbitrary complete orthonormal system
of functions in $L_2([t, T])$ \cite{10a} (Sect.~2.1.4, 2.23, 2.24, 2.31--2.34).

\vspace{5mm}

\section{Theorems 3--5, 12--15 from Point
of View of the Wong--Zakai Approximation}

\vspace{5mm}

The iterated Ito stochastic integrals and solutions
of Ito SDEs are complex and important func\-ti\-onals
from the independent components ${\bf f}_{s}^{(i)},$
$i=1,\ldots,m$ of the multidimensional
Wiener process ${\bf f}_{s},$ $s\in[0, T].$
Let ${\bf f}_{s}^{(i)p},$ $p\in\mathbb{N}$ 
be some approximation of
${\bf f}_{s}^{(i)},$
$i=1,\ldots,m$.
Suppose that 
${\bf f}_{s}^{(i)p}$
converges to
${\bf f}_{s}^{(i)},$
$i=1,\ldots,m$ if $p\to\infty$ in some sense and has
differentiable sample trajectories.

A natural question arises: if we replace 
${\bf f}_{s}^{(i)}$
by ${\bf f}_{s}^{(i)p},$
$i=1,\ldots,m$ in the functionals
mentioned above, will the resulting
functionals converge to the original
functionals from the components 
${\bf f}_{s}^{(i)},$
$i=1,\ldots,m$ of the multidimentional
Wiener process ${\bf f}_{s}$?
The answere to this question is negative 
in the general case. However, 
in the pioneering works of Wong E. and Zakai M. \cite{W-Z-1},
\cite{W-Z-2},
it was shown that under the special conditions and 
for some types of approximations 
of the Wiener process the answere is affirmative
with one peculiarity: the convergence takes place 
to the iterated Stratonovich stochastic integrals
and solutions of Stratonovich SDEs and not to iterated 
Ito stochastic integrals and solutions
of Ito SDEs.
The piecewise 
linear approximation 
as well as the regularization by convolution 
\cite{W-Z-1}-\cite{Watanabe} relate the 
mentioned types of approximations
of the Wiener process. The above approximation 
of stochastic integrals and solutions of SDEs 
is often called the Wong--Zakai approximation.

Let ${\bf w}_{\tau},$ $\tau\in[0, T]$ is a random vector with 
an $m+1$ components: ${\bf w}_{\tau}^{(i)}={\bf f}_{\tau}^{(i)}$ 
for $i=1,\ldots,m$ and 
${\bf w}_{\tau}^{(0)}=\tau,$\ 
${\bf f}_{\tau}^{(i)}$ $(i=1,\ldots,m)$
are independent standard Wiener processes.

It is well known that the following representation 
takes place \cite{Lipt}, \cite{7e}

\begin{equation}
\label{um1x}
{\bf w}_{\tau}^{(i)}-{\bf w}_{t}^{(i)}=
\sum_{j=0}^{\infty}\int\limits_t^{\tau}
\phi_j(s)ds\ \zeta_j^{(i)},\ \ \ \zeta_j^{(i)}=
\int\limits_t^T \phi_j(s)d{\bf w}_s^{(i)},
\end{equation}

\vspace{4mm}
\noindent
where $\tau\in[t, T],$ $t\ge 0,$
$\{\phi_j(x)\}_{j=0}^{\infty}$ is an arbitrary complete 
orthonormal system of functions in the space $L_2([t, T]),$ and
$\zeta_j^{(i)}$ are independent standard Gaussian 
random variables for various $i$ or $j.$
Moreover, the series (\ref{um1x}) converges for any $\tau\in [t, T]$
in the mean-square sense.

Let ${\bf w}_{\tau}^{(i)p}-{\bf w}_{t}^{(i)p}$ be 
the mean-square approximation of the process
${\bf w}_{\tau}^{(i)}-{\bf w}_{t}^{(i)},$
which has the following form

\vspace{-3mm}
\begin{equation}
\label{um1xx}
{\bf w}_{\tau}^{(i)p}-{\bf w}_{t}^{(i)p}=
\sum_{j=0}^{p}\int\limits_t^{\tau}
\phi_j(s)ds\ \zeta_j^{(i)}.
\end{equation}

\vspace{3mm}

From (\ref{um1xx}) we obtain

\vspace{-4mm}
\begin{equation}
\label{um1xxx}
d{\bf w}_{\tau}^{(i)p}=
\sum_{j=0}^{p}
\phi_j(\tau)\zeta_j^{(i)} d\tau.
\end{equation}

\vspace{4mm}

Consider the following iterated Riemann--Stieltjes
integral

\begin{equation}
\label{um1xxxx}
\int\limits_t^T
\psi_k(t_k)\ldots \int\limits_t^{t_2}\psi_1(t_1)
d{\bf w}_{t_1}^{(i_1)p_1}\ldots d{\bf w}_{t_k}^{(i_k)p_k},
\end{equation}

\vspace{4mm}
\noindent
where $i_1,\ldots,i_k=0,1,\ldots,m,$\ \ $p_1,\ldots,p_k\in \mathbb{N},$ 

\begin{equation}
\label{um1xxx1}
d{\bf w}_{\tau}^{(i)p}=
\left\{\begin{matrix}
d{\bf f}_{\tau}^{(i)p}\ &\hbox{\rm for}\ \ \ i=1,\ldots,m\cr\cr\cr
d\tau^p\ &\hbox{\rm for}\ \ \ i=0
\end{matrix}
,\right.
\end{equation}

\vspace{4mm}
\noindent
and $d{\bf f}_{\tau}^{(i)p},$ $d\tau^p$ are defined by the relation (\ref{um1xxx}).

Let us substitute (\ref{um1xxx}) into (\ref{um1xxxx})

\begin{equation}
\label{um1xxxx1}
\int\limits_t^T
\psi_k(t_k)\ldots \int\limits_t^{t_2}\psi_1(t_1)
d{\bf w}_{t_1}^{(i_1)p_1}\ldots d{\bf w}_{t_k}^{(i_k)p_k}=
\sum\limits_{j_1=0}^{p_1} \ldots \sum\limits_{j_k=0}^{p_k}
C_{j_k \ldots j_1}\prod\limits_{l=1}^k \zeta_{j_l}^{(i_l)},
\end{equation}

\vspace{4mm}
\noindent
where 
$$
\zeta_j^{(i)}=\int\limits_t^T \phi_j(s)d{\bf w}_s^{(i)}
$$ 

\vspace{2mm}
\noindent
are independent standard Gaussian random variables for various 
$i$ or $j$ (in the case when $i\ne 0$),
${\bf w}_{s}^{(i)}={\bf f}_{s}^{(i)}$ for
$i=1,\ldots,m$ and 
${\bf w}_{s}^{(0)}=s,$

$$
C_{j_k \ldots j_1}=\int\limits_t^T\psi_k(t_k)\phi_{j_k}(t_k)\ldots
\int\limits_t^{t_2}
\psi_1(t_1)\phi_{j_1}(t_1)
dt_1\ldots dt_k
$$

\vspace{4mm}
\noindent
is the Fourier coefficient.

To best of our knowledge \cite{W-Z-1}-\cite{Watanabe}
the approximations of the Wiener process
in the Wong--Zakai approximation must satisfy fairly strong
restrictions
\cite{Watanabe}
(see Definition 7.1, pp.~480--481).
Moreover, approximations of the Wiener process that are
similar to (\ref{um1xx})
were not considered in \cite{W-Z-1}, \cite{W-Z-2}
(also see \cite{Watanabe}, Theorems 7.1, 7.2).
Therefore, the proof of analogs of Theorems 7.1 and 7.2 \cite{Watanabe}
for approximations of the Wiener 
process based on its series expansion (\ref{um1x})
should be carried out separately.
Thus, the mean-square convergence of the right-hand side
of (\ref{um1xxxx1}) to the iterated Stratonovich stochastic integral 
(\ref{str})
does not follow from the results of the papers
\cite{W-Z-1}, \cite{W-Z-2} (also see \cite{Watanabe},
Theorems 7.1, 7.2).

Nevertheless, the authors of the works
\cite{1}
(Sect.~5.8, pp.~202--204), 
\cite{Zapad-2} (pp.~438-439), \cite{Zapad-4} (pp.~82-84),  
\cite{Zapad-9} (pp.~263-264) use 
the Wong--Zakai approximation 
\cite{W-Z-1}-\cite{Watanabe} (without rigorous proof) within the frames
of the method of expansion of iterated stochastic integrals
based on the trigonometric series expansion 
of the Brownian bridge process (version
of the so-called Karhunen--Loeve expansion).

From the other hand, Theorems 3--5, 12--15 from this 
paper can be considered as the proof of the
Wong--Zakai approximation based on the iterated 
Riemann--Stieltjes integrals (\ref{um1xxxx}) of multiplicities 1 to 6
and the Wiener process approximation (\ref{um1xx}) 
on the base of its series expansion.
At that, the mentioned Riemann--Stieltjes integrals converge
(according to Theorems 3--5, 12--15)
to the appropriate Stratonovich 
stochastic integrals (\ref{str}) of multiplicities 1 to 6. Recall that
$\{\phi_j(x)\}_{j=0}^{\infty}$ (see (\ref{um1x}), (\ref{um1xx}), and
Theorems 3, 12--15)
is a complete 
orthonormal system of Legendre polynomials or 
trigonometric functions 
in the space $L_2([t, T])$.

To illustrate the above reasoning, 
consider two examples for the case $k=2,$
$\psi_1(s),$ $\psi_2(s)\equiv 1;$ $i_1, i_2=1,\ldots,m.$

The first example relates to the piecewise linear approximation
of the multidimensional Wiener process (these approximations 
were considered in \cite{W-Z-1}-\cite{Watanabe}).

Let ${\bf b}_{\Delta}^{(i)}(t),$ $t\in[0, T]$ be the piecewise
linear approximation of the $i$th component ${\bf f}_t^{(i)}$
of the multidimensional standard Wiener process ${\bf f}_t,$
$t\in [0, T]$ with independent components
${\bf f}_t^{(i)},$ $i=1,\ldots,m,$ i.e.

$$
{\bf b}_{\Delta}^{(i)}(t)={\bf f}_{k\Delta}^{(i)}+
\frac{t-k\Delta}{\Delta}\Delta{\bf f}_{k\Delta}^{(i)},
$$

\vspace{3mm}
\noindent
where 

\vspace{-2mm}
$$
\Delta{\bf f}_{k\Delta}^{(i)}={\bf f}_{(k+1)\Delta}^{(i)}-
{\bf f}_{k\Delta}^{(i)},\ \ \
t\in[k\Delta, (k+1)\Delta),\ \ \ k=0, 1,\ldots, N-1.
$$

\vspace{4mm}

Note that w.~p.~1

\vspace{-1mm}
\begin{equation}
\label{pridum}
\frac{d{\bf b}_{\Delta}^{(i)}}{dt}(t)=
\frac{\Delta{\bf f}_{k\Delta}^{(i)}}{\Delta},\ \ \
t\in[k\Delta, (k+1)\Delta),\ \ \ k=0, 1,\ldots, N-1.
\end{equation}

\vspace{4mm}

Consider the following iterated Riemann--Stieltjes
integral

\vspace{1mm}
$$
\int\limits_0^T
\int\limits_0^{s}
d{\bf b}_{\Delta}^{(i_1)}(\tau)d{\bf b}_{\Delta}^{(i_2)}(s),\ \ \ 
i_1,i_2=1,\ldots,m.
$$

\vspace{4mm}

Using (\ref{pridum}) and additive property of Riemann--Stieltjes integrals, 
we can write w.~p.~1

\vspace{2mm}
$$
\int\limits_0^T
\int\limits_0^{s}
d{\bf b}_{\Delta}^{(i_1)}(\tau)d{\bf b}_{\Delta}^{(i_2)}(s)=
\int\limits_0^T
\int\limits_0^{s}
\frac{d{\bf b}_{\Delta}^{(i_1)}}{d\tau}(\tau)d\tau
\frac{d {\bf b}_{\Delta}^{(i_2)}}{d s}(s)
ds =
$$

\vspace{3mm}
$$
=
\sum\limits_{l=0}^{N-1}\int\limits_{l\Delta}^{(l+1)\Delta}
\left(
\sum\limits_{q=0}^{l-1}\int\limits_{q\Delta}^{(q+1)\Delta}
\frac{\Delta{\bf f}_{q\Delta}^{(i_1)}}{\Delta}d\tau+
\int\limits_{l\Delta}^{s}
\frac{\Delta{\bf f}_{l\Delta}^{(i_1)}}{\Delta}d\tau\right)
\frac{\Delta{\bf f}_{l\Delta}^{(i_2)}}{\Delta}ds=
$$

\vspace{3mm}
$$
=\sum\limits_{l=0}^{N-1}\sum\limits_{q=0}^{l-1}
\Delta{\bf f}_{q\Delta}^{(i_1)}
\Delta{\bf f}_{l\Delta}^{(i_2)}+
\frac{1}{\Delta^2}\sum\limits_{l=0}^{N-1}
\Delta{\bf f}_{l\Delta}^{(i_1)}
\Delta{\bf f}_{l\Delta}^{(i_2)}
\int\limits_{l\Delta}^{(l+1)\Delta}
\int\limits_{l\Delta}^{s}d\tau ds=
$$

\vspace{3mm}
\begin{equation}
\label{oh-ty}
=\sum\limits_{l=0}^{N-1}\sum\limits_{q=0}^{l-1}
\Delta{\bf f}_{q\Delta}^{(i_1)}
\Delta{\bf f}_{l\Delta}^{(i_2)}+
\frac{1}{2}\sum\limits_{l=0}^{N-1}
\Delta{\bf f}_{l\Delta}^{(i_1)}
\Delta{\bf f}_{l\Delta}^{(i_2)}.
\end{equation}

\vspace{6mm}

Using (\ref{oh-ty}) and standard relation between Stratonovich
and Ito stochastic integrals, it 
is not difficult to show 
that

\vspace{1mm}
$$
\hbox{\vtop{\offinterlineskip\halign{
\hfil#\hfil\cr
{\rm l.i.m.}\cr
$\stackrel{}{{}_{N\to \infty}}$\cr
}} }
\int\limits_0^T
\int\limits_0^{s}
d{\bf b}_{\Delta}^{(i_1)}(\tau)d{\bf b}_{\Delta}^{(i_2)}(s)=
\int\limits_0^T
\int\limits_0^{s}
d{\bf f}_{\tau}^{(i_1)}d{\bf f}_{s}^{(i_2)}+
\frac{1}{2}{\bf 1}_{\{i_1=i_2\}}\int\limits_0^T ds=
$$

\vspace{3mm}
\begin{equation}
\label{uh-111}
=
{\int\limits_0^{*}}^T
{\int\limits_0^{*}}^s
d{\bf f}_{\tau}^{(i_1)}d{\bf f}_{s}^{(i_2)},
\end{equation}

\vspace{5mm}
\noindent
where $\Delta\to 0$ if $N\to\infty$ ($N\Delta=T$).

Obviously, (\ref{uh-111}) agrees with Theorem 7.1 (see \cite{Watanabe},
p.~486).

The next example relates to the approximation
of the Wiener process based on its series expansion
(\ref{um1x}) for $t=0$, where
$\{\phi_j(x)\}_{j=0}^{\infty}$ 
is a complete 
orthonormal system of Legendre polynomials or 
trigonometric functions 
in the space $L_2([0, T])$.

Consider the following iterated Riemann--Stieltjes
integral

\vspace{-1mm}
\begin{equation}
\label{abcd1}
\int\limits_0^T
\int\limits_0^{s}
d{\bf f}_{\tau}^{(i_1)p}d{\bf f}_{s}^{(i_2)p},\ \ \ 
i_1,i_2=1,\ldots,m,
\end{equation}

\vspace{3mm}
\noindent
where $d{\bf f}_{\tau}^{(i)p}$ is defined by the
relation
(\ref{um1xxx}).

Let us substitute (\ref{um1xxx}) into (\ref{abcd1}) 

\vspace{-1mm}
\begin{equation}
\label{set18}
\int\limits_0^T
\int\limits_0^{s}
d{\bf f}_{\tau}^{(i_1)p}d{\bf f}_{s}^{(i_2)p}=
\sum\limits_{j_1,j_2=0}^p
C_{j_2 j_1} \zeta_{j_1}^{(i_1)}\zeta_{j_2}^{(i_2)},
\end{equation}

\vspace{3mm}
\noindent
where 
$$
C_{j_2 j_1}=
\int\limits_0^T \phi_{j_2}(s)\int\limits_0^s
\phi_{j_1}(\tau)d\tau ds
$$

\vspace{3mm}
\noindent
is the Fourier coefficient; another notations 
are the same as in (\ref{um1xxxx1}).

As we noted above, approximations of the Wiener process that are
similar to (\ref{um1xx})
were not considered in \cite{W-Z-1}, \cite{W-Z-2}
(also see Theorems 7.1, 7.2 in \cite{Watanabe}).
Furthermore, the extension of the results of Theorems 7.1 and 7.2
\cite{Watanabe} to the case under consideration is
not obvious.               

On the other hand, we can apply the theory built in Chapters 1 and 2
of the monographs \cite{10a}-\cite{12aa-afterxxx}. More precisely, 
using 
Theorems 3, 5 from this paper,
we obtain from (\ref{set18}) the desired result

\vspace{-1mm}
$$
\hbox{\vtop{\offinterlineskip\halign{
\hfil#\hfil\cr
{\rm l.i.m.}\cr
$\stackrel{}{{}_{p\to \infty}}$\cr
}} }
\int\limits_0^T
\int\limits_0^{s}
d{\bf f}_{\tau}^{(i_1)p}d{\bf f}_{s}^{(i_2)p}=
\hbox{\vtop{\offinterlineskip\halign{
\hfil#\hfil\cr
{\rm l.i.m.}\cr
$\stackrel{}{{}_{p\to \infty}}$\cr
}} }
\sum\limits_{j_1,j_2=0}^p
C_{j_2 j_1} \zeta_{j_1}^{(i_1)}\zeta_{j_2}^{(i_2)}=
$$

\vspace{2mm}
\begin{equation}
\label{umen-bl}
=
{\int\limits_0^{*}}^T
{\int\limits_0^{*}}^s
d{\bf f}_{\tau}^{(i_1)}d{\bf f}_{s}^{(i_2)}.
\end{equation}

\vspace{5mm}

From the other hand, by Theorem 4
(see (\ref{leto5001})) for the case
$k=2$ we obtain from (\ref{set18}) the following relation

\vspace{-2mm}
$$
\hbox{\vtop{\offinterlineskip\halign{
\hfil#\hfil\cr
{\rm l.i.m.}\cr
$\stackrel{}{{}_{p\to \infty}}$\cr
}} }
\int\limits_0^T
\int\limits_0^{s}
d{\bf f}_{\tau}^{(i_1)p}d{\bf f}_{s}^{(i_2)p}=
\hbox{\vtop{\offinterlineskip\halign{
\hfil#\hfil\cr
{\rm l.i.m.}\cr
$\stackrel{}{{}_{p\to \infty}}$\cr
}} }
\sum\limits_{j_1,j_2=0}^p
C_{j_2 j_1} \zeta_{j_1}^{(i_1)}\zeta_{j_2}^{(i_2)}=
$$

\vspace{2mm}
$$
=
\hbox{\vtop{\offinterlineskip\halign{
\hfil#\hfil\cr
{\rm l.i.m.}\cr
$\stackrel{}{{}_{p\to \infty}}$\cr
}} }
\sum\limits_{j_1,j_2=0}^p
C_{j_2 j_1} \biggl(\zeta_{j_1}^{(i_1)}\zeta_{j_2}^{(i_2)}-
{\bf 1}_{\{i_1=i_2\}}{\bf 1}_{\{j_1=j_2\}}\biggr)+
{\bf 1}_{\{i_1=i_2\}}\sum\limits_{j_1=0}^{\infty}
C_{j_1 j_1}=
$$

\vspace{2mm}
\begin{equation}
\label{umen-blx}
=
\int\limits_0^T
\int\limits_0^{s}
d{\bf f}_{\tau}^{(i_1)}d{\bf f}_{s}^{(i_2)}+
{\bf 1}_{\{i_1=i_2\}}\sum\limits_{j_1=0}^{\infty}
C_{j_1 j_1}.
\end{equation}

\vspace{5mm}

Since
$$
\sum\limits_{j_1=0}^{\infty}
C_{j_1 j_1}=\frac{1}{2}\sum\limits_{j_1=0}^{\infty}
\left(\int\limits_0^T \phi_j(\tau)d\tau\right)^2
=\frac{1}{2}
\left(\int\limits_0^T \phi_0(\tau)d\tau\right)^2=\frac{1}{2}
\int\limits_0^T ds,
$$

\vspace{4mm}
\noindent
then from (\ref{umen-blx}) we obtain (\ref{umen-bl}).

\end{document}